\documentclass[11pt]{article}
\usepackage[a4paper,margin=1in]{geometry}
\usepackage{amsmath,amssymb,amsthm,mathtools,bm}
\usepackage{geometry}
\usepackage{enumitem}
\usepackage{booktabs}
\usepackage{subfigure}
\usepackage{hyperref}
\usepackage[nameinlink,noabbrev]{cleveref}

\numberwithin{equation}{section}

\counterwithin{table}{section}
\counterwithin{figure}{section}
\newtheorem{definition}{Definition}[section]
\newtheorem{theorem}{Theorem}[section]
\newtheorem{proposition}{Proposition}[section]

\numberwithin{equation}{section}
\makeatletter
\renewcommand*{\p@subfigure}{\thefigure}
\makeatother

\title{Triply conservative schemes for degenerate convection-dispersion equations}
\author{
	Kaichang Yu\thanks{Center for Interdisciplinary Applied Mathematics,
		Zhejiang University, Hangzhou 310000, China.
		Email: \texttt{kaichangyu@zju.edu.cn}.}
	\and
	Juan Cheng\thanks{Corresponding author. Academy for Multidisciplinary Studies, Capital Normal University, Beijing 100048, China.
Email: \texttt{jcheng@cnu.edu.cn}. Research partially supported by Science Challenge Project TZ2025007 and NSFC grant 12571432.}
	\and
	Chi-Wang Shu\thanks{Division of Applied Mathematics, Brown University, Providence, RI 02912, USA.
		Email: \texttt{chi-wang\_shu@brown.edu}.}
}

\date{}

\begin{document}

\maketitle	
\begin{abstract}
It is well known that for numerically solving hyperbolic conservation laws, even a convergent numerical scheme still may not always converge to a weak solution. The remedy is the famous Lax-Wendroff theorem. 
In a recent work on degenerate convection-diffusion equations,
an analogous phenomenon was identified and doubly
conservative schemes were introduced.
In this paper, we address the same issue for solving the following nonlinear degenerate dispersive or convection-dispersion equations:
\begin{equation*} 
u_t+f(u)_x+A(u)_{xxx}=0.
\end{equation*}
The resulting equation is regarded as degenerate if $A'(u)$ is allowed to vanish. An important example is the Rosenau-Hyman $K(m,n)$ equation. In contrast to classical KdV solitons possessing exponentially decaying tails, $K(m,n)$ equations can admit compactly supported solitary waves, known as compactons, whose regularity can deteriorate at the boundary of their support. We first give an example in which both a classical conservative and doubly conservative scheme converges to a limit that fails to satisfy the ultra-weak formulation. To address this issue, we introduce the notion of a triply conservative (TriC) scheme and prove a Lax-Wendroff type theorem for consistent TriC schemes. Under suitable assumptions, a convergent sequence generated by a consistent TriC scheme will converge to an ultra-weak solution. Several existing high-order schemes are shown to have a TriC structure. We further introduce a new high-order TriC finite volume scheme based on a connection between triple cell averages and quadratic B-splines. The theory and the scheme are also suited for two-dimensional anisotropic problems and non-uniform meshes. Numerical experiments in one and two dimensions confirm the high-order accuracy on smooth solutions and illustrate the behavior of the schemes for compactons and traveling fronts.
\end{abstract}


\section{Introduction}
We consider the following nonlinear dispersive equation:
\begin{equation}\label{1D TC model}
u_{t} + f(u)_{x} + A(u)_{xxx}=0,
\end{equation}
and its two-dimensional anisotropic generalization.
The equation may be equipped with different initial and boundary conditions. $f(u)$ and $A(u)$ in \eqref{1D TC model} are Lipschitz continuous functions of $u$. 
If $A'(u)$ may be zero somewhere, then the equation \eqref{1D TC model} is called a degenerate dispersive equation.
The class of equations in the form \eqref{1D TC model}, together with its two-dimensional extension, includes the classical Korteweg-de Vries (KdV) equation \cite{KdV_original}, the Zakharov–Kuznetsov equation \cite{ZETF_1974_ZK-equation} and the compacton equations proposed in \cite{PRL_1993_Rosenau_Hyman}. The KdV equation has nonlinear convection and linear dispersion, which may generate soliton solutions. 
These solitons are localized and have exponentially decaying tails but are not compactly supported.
In 1993, to study the nonlinear dispersion in pattern formation, Rosenau and Hyman introduced the $K(m,n)$ equation which means that $f(u)=u^m$ and $A(u)=u^n$ with  non-negative integers $m,n$ and $n>1$. 
The $K(m,n)$ equation is degenerate if $u=0$ somewhere. This will generate dispersive solutions with low-regularity, compact support and sharp moving interfaces. 
A representative example is the compacton, which is also introduced in \cite{PRL_1993_Rosenau_Hyman} and recognized as a soliton with compact support.
In addition, $K(0,n)$ type equations may also possess one-sided traveling-front solutions \cite{StAM_2024_Anco_Gandarias}, which vanish on one side of the interface and are nonzero on the other side. 
These two types of solutions can have limited regularity at their interfaces, so we consider them in the ultra-weak sense.
For analyses of low-regularity solutions to degenerate dispersive equations, see \cite{StAM_2024_Anco_Gandarias,CPAM_2019_Germain_Benjamin_Jeremy, Nonlinearity_2018_Zilburg_Rosenau}.

The low regularity also poses a challenge for numerical approximation.
A standard conservative discretization may converge to a limit that does not satisfy the ultra-weak formulation.
The problem setting is shown in Section~\ref{compare_1D} and here we only present the comparison between a numerical solution and the exact weak solution. For this problem, the solution is computed on $400$ cells and also 
numerically converges as the mesh refines. This result is also the same regardless of whether backward Euler or explicit third-order strong-stability-preserving Runge-Kutta (SSPRK3) \cite{SIREV_2001_Gottlieb_Shu_Tadmor} time discretization is used.
\begin{figure}[htbp]
	\centering
	\includegraphics[width=.35\textwidth]{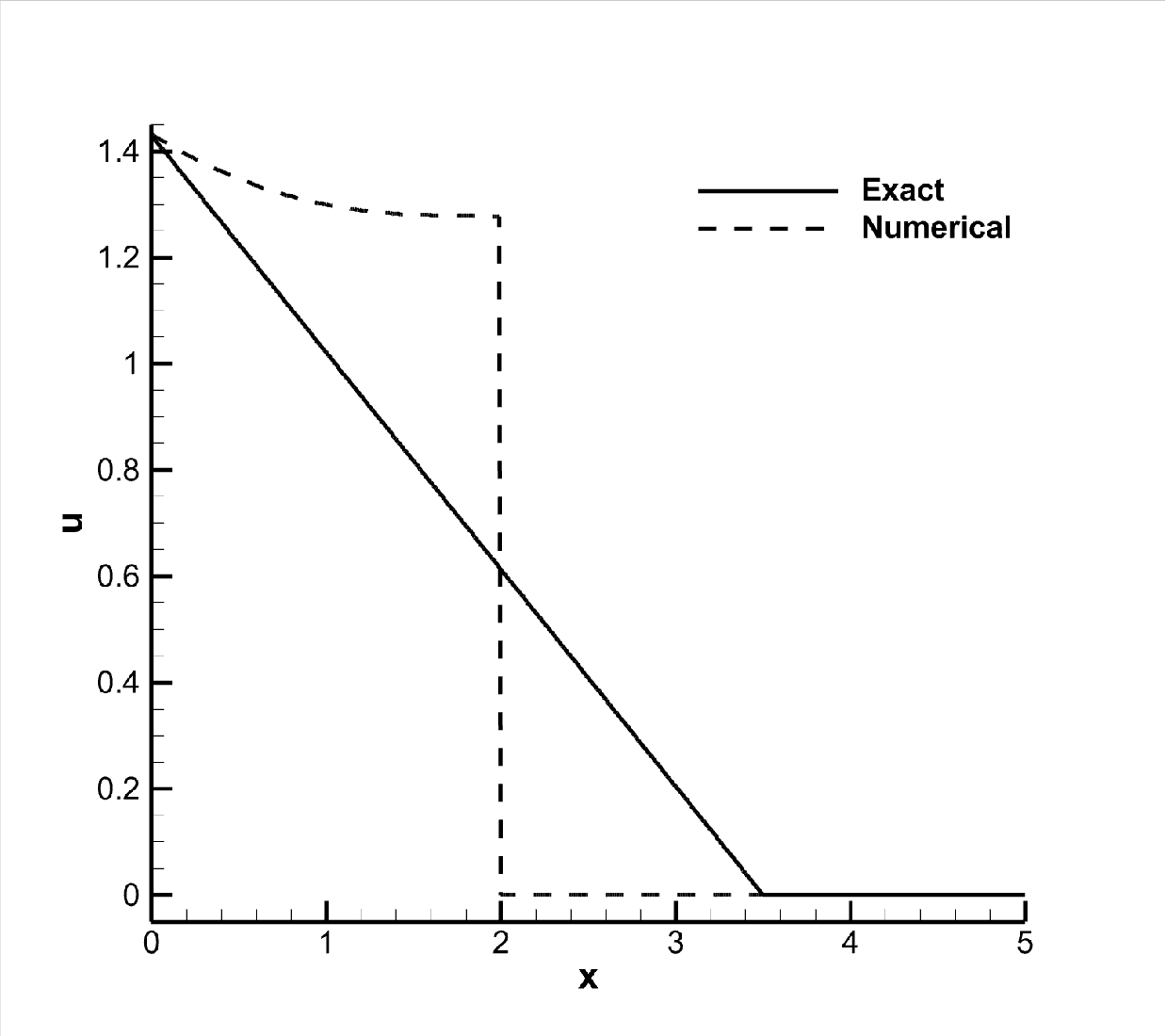}
	\caption{Comparison between an exact solution of \eqref{1D TC model} and a conservative numerical solution.}
\end{figure}

This phenomenon is similar to the case in the numerical approximation of hyperbolic conservation laws \cite{Lax_Wendroff_1960} and degenerate convection-diffusion equations \cite{MCOM_2026_Yu_Cheng_Liu_Shu}.
In hyperbolic conservation laws, the classical Lax-Wendroff theorem \cite{Lax_Wendroff_1960} provides a fundamental framework ensuring that a conservative numerical scheme, under suitable assumptions, when convergent, always converges to a weak solution.
This theorem also provides a principal theoretical motivation for constructing numerical schemes in conservative form. 
More recently, an analogous issue has been investigated in the approximation of degenerate convection-diffusion equations, 
for which classical conservative schemes are not enough to ensure that the numerical solutions converge to an ultra-weak solution. Therefore, this phenomenon leads to the construction of doubly conservative schemes in \cite{MCOM_2026_Yu_Cheng_Liu_Shu}.

It is worth mentioning that the numerical scheme used in the above example is in both the conservative and the
doubly conservative form. Therefore, for the third-order degenerate dispersive equation, neither the result in \cite{Lax_Wendroff_1960} nor that of \cite{MCOM_2026_Yu_Cheng_Liu_Shu} suffices. This motivates us to develop a similar theory for the degenerate dispersive equation, namely, the theory of triply conservative (TriC) schemes. Furthermore, compared with the doubly conservative theory, the TriC construction links triple cell averages to quadratic B-spline moments, and our Lax-Wendroff type convergence proof uses a auxiliary linear functional of the test function.

The literature contains many numerical schemes for solving the degenerate dispersive equation \eqref{1D TC model}, such as finite difference methods \cite{MCS_2023_Ahmat_Qiu, JCP_1995_deFrutos_Lopez_Sanz}, discontinuous Galerkin methods \cite{JCP_2004_Levy_Shu_Yan, CAMC_2026_Zheng_Xu}, the particle method \cite{JCP_2001_Chertock_Levy}, Pad\'e approximation \cite{MCS_2007_Rus_Villatoro},  the operator splitting method \cite{JCP_1999_Holden_Karlsen_Risebro} and so on. 
In particular, in \cite{MCS_2023_Ahmat_Qiu}, the authors constructed a finite difference scheme based on the triple cell average on uniform meshes, but they did not discuss whether the scheme might converge to an incorrect limit. In \cite{JCP_2004_Levy_Shu_Yan}, the authors constructed a special energy stable local discontinuous Galerkin (LDG) scheme for the $K(m,m)$ equation with an odd $m$. In \cite{JCP_1999_Holden_Karlsen_Risebro}, the authors discussed the temporal semi-discretization and proved a Lax-Wendroff type theorem in the time direction. 
Unlike the work in \cite{JCP_1999_Holden_Karlsen_Risebro}, we mainly consider the spatial discretization, especially for the nonlinear degenerate dispersion term. In our view, the problem presented above is mainly caused by the spatial discretization and we will give a Lax-Wendroff type theorem for a fully discrete scheme which is also suited for high-order numerical schemes. Furthermore, some existing schemes belong to our newly defined TriC category, including \cite{MCS_2023_Ahmat_Qiu,MCOM_2008_ultraweak,JCP_2004_Levy_Shu_Yan}.

This paper is organized as follows. In Section~\ref{section 2}, we first introduce the definition of TriC schemes and the corresponding consistency in one dimension. Then we prove a Lax-Wendroff type theorem for the TriC scheme and construct a new high-order TriC finite volume (TriC-FV) scheme for solving dispersive equations. When constructing the scheme, we also establish the connection between the triple cell average and the quadratic B-spline, which is a key observation for constructing the scheme on non-uniform meshes. Finally we also present some high-order schemes in the literature as examples.
In Section~\ref{Sec3}, we generalize all the results in Section~\ref{section 2} to two dimensions. Then in Section~\ref{Sec4}, we will present the numerical results of the newly constructed TriC-FV scheme. We test the optimal accuracy on non-uniform meshes, compare the TriC scheme and non-TriC scheme for solving the degenerate dispersive equation, compare the high-order TriC-FV scheme and the low-order scheme at the equal cells and the equal freedoms, and also simulate the dynamics of compactons. Both one-dimensional and two-dimensional results are presented. Finally, a brief conclusion is given in Section~\ref{Sec5}.

\section{Triply conservative schemes: one-dimensional case}\label{section 2}
Following the spirit of the pioneering work of Lax and Wendroff \cite{Lax_Wendroff_1960} for hyperbolic conservation laws, we will present the definition of the triply conservative schemes and the corresponding notion of consistency in one dimension. Then we will prove a Lax-Wendroff type theorem of the degenerate dispersive equation \eqref{1D TC model} and present some examples of TriC schemes.
Before we introduce the TriC scheme, using the same definition introduced in \cite{JCP_1999_Holden_Karlsen_Risebro,StAM_2024_Anco_Gandarias}, we first define the ultra-weak solution of the equation \eqref{1D TC model}.
\begin{definition}
	A function $u:\mathbb{R}^{2}\to\mathbb{R}$ is an ultra-weak solution of \eqref{1D TC model} if it satisfies that, for $u_0\in L^1_{loc}(\mathbb{R})$, $u,\ f(u),\ A(u)\in L^1_{loc}\bigl(\mathbb{R}\times\mathbb{R}_{+}\bigr)$, the following equation
	\begin{equation}\label{eq:ultra-weak-formulation}
		\int_{\mathbb{R}_{+}}\int_{\mathbb{R}}
		\left(
		u \partial_t\phi
		+f(u) \partial_x\phi
		+A(u) \partial_{xxx}\phi
		\right) dx dt
		+\int_{\mathbb{R}}u_0(x)\phi(x,0) dx=0,
	\end{equation}
	holds for every $\phi\in C_c^\infty(\mathbb{R}^{2})$.
\end{definition}

\subsection{Definition of triply conservative schemes}
In this section, we will give the definition of TriC schemes and the corresponding consistency. Decompose the domain $\Omega=[x_{L},x_{R}]$ into $N$ cells on a  non-uniform mesh $x_{L}=x_{\frac{1}{2}}<x_{\frac{3}{2}}<\cdots<x_{N+\frac{1}{2}}=x_{R}$. Denote $\Delta x_{j} = x_{j+\frac{1}{2}} - x_{j-\frac{1}{2}}$, $I_{j} = [x_{j-\frac{1}{2}},x_{j+\frac{1}{2}}]$ and $x_{j}$ is the midpoint of $I_{j}$. 
The temporal domain is discretized into $0 = t_{0}<t_{1}<\cdots$.  
We define the maximal mesh size $h=\max_{j\in\mathbb{Z}}|I_{j}|$, where $|I_j|$ is the length of $I_j$, and assume that $|I_{j}|>c_{0}h$ holds for any $j\in\mathbb{Z}$, where $c_{0}>0$ is independent of the mesh size.
The time step size is defined as $\Delta t_{s} = t_{s+1} - t_{s}$ and $\tau = \max_{s}\Delta t_{s}$.
Since the numerical solution is only defined on temporal slices $\mathbb{R}\times \{t_{s}\}_{s=0}^{ \infty }$, to generalize the numerical solution $u_h(x,t)$ to a function over the whole domain $\mathbb{R}\times \mathbb{R}^{+} $, we extend it to a piecewise constant function over temporal intervals: $u_h(x,t) \equiv u_h(x,t_{s}),\ \forall t\in [t_{s},t_{s+1})$.
For simplicity, we denote $\sum_{s=0}^{+\infty}$ as $\sum_{s}$ and $\sum_{j\in\mathbb{Z}}$ as $\sum_{j}$.
As spatial discretization is our primary focus, we restrict our consideration to spatial formulation and only the forward Euler temporal discretization will be used.
The definition of TriC schemes of \eqref{1D TC model} is given as follows.

\begin{definition}[Triply Conservative Scheme]\label{1D def tric}
	A numerical scheme for \eqref{1D TC model} is called a TriC scheme if it can be written in the following form:
	\begin{equation}\label{1D TriC Def}
		\small
	\begin{aligned}
		\frac{\widetilde{u}_{j}^{s+1} - \widetilde{u}_{j}^{s}}{\Delta t_{s}} =  - \frac{\widetilde{f}_{j+\frac{1}{2}}^{s} - \widetilde{f}_{j-\frac{1}{2}}^{s}}{V_{j}} 
		- \frac{1}{V_j} \bigg( & \frac{2}{\Delta x_{j}+\Delta x_{j+1}} \bigg( \frac{\widetilde{A}^{s}_{j+\frac{3}{2}} - \widetilde{A}^{s}_{j+\frac{1}{2}}}{\Delta x_{j+1}} - \frac{\widetilde{A}^{s}_{j+\frac{1}{2}} - \widetilde{A}^{s}_{j-\frac{1}{2}}}{\Delta x_{j}} \bigg) - \\ &\frac{2}{\Delta x_{j-1}+\Delta x_{j}} \bigg( \frac{\widetilde{A}_{j+\frac{1}{2}}^{s} - \widetilde{A}^{s}_{j-\frac{1}{2}}}{\Delta x_{j}} - \frac{\widetilde{A}_{j-\frac{1}{2}}^{s} - \widetilde{A}^{s}_{j-\frac{3}{2}}}{\Delta x_{j-1}} \bigg)\bigg),
		\end{aligned}
	\end{equation}
	where $\widetilde{u}_{j}^{s}$ is the generalized locally conserved variable defined in \cite{C&F_2018_shi_shu}, $\widetilde{f}_{j+\frac{1}{2}}^{s}$ and $\widetilde{A}_{j+\frac{1}{2}}^{s}$ are numerical fluxes. $V_{j}$ is the local volume corresponding to the locally conserved variable $\widetilde{u}_{j}^{s}$. All these tilde-quantities are functions locally depending on the numerical solution $u_{h}^{s}(B_{j+\frac{1}{2}})$, where $B_{j+\frac{1}{2}} = \{x\in \mathbb{R}:|x-x_{j+\frac{1}{2}}|\leq ch \}$, with $c>0$ being  a positive integer independent of the mesh sizes.
\end{definition}
For different schemes, the locally conserved variable could be viewed as a different local spatial average of numerical solutions $u_{h}^{s}$. These different averages could be viewed as different non-negative mass one
kernels $K_{h,j}$ in $\widetilde{u}_{j}(t) = \int_{\mathbb{R}}u_{h}(x,t) K_{h,j}\left(x\right) dx$. If 
$K_{h,j}(x) = \delta(x-x_{j})$, then $\widetilde{u}_{j}$ is the point value; if $K_{h,j}\left(x\right)= \frac{1}{|I_{j}|} 
I_{I_{j}}(x)$ is the indicator function of $I_{j}$ divided by the length $|I_{j}|$, then $\widetilde{u}_{j} = \bar{u}_{j}$ is the usual cell average, and if $K_{h,j}(x) = 
\frac{2}{\Delta x_{j} + \Delta x_{j+1}}\begin{cases}
	\small 
	\dfrac{x-x_{j-\frac{1}{2}}}{\Delta x_j},
	& x\in [x_{j-\frac{1}{2}},x_{j+\frac{1}{2}}], \\
	\dfrac{x_{j+\frac{3}{2}}-x}{\Delta x_{j+1}},
	& x\in [x_{j+\frac{1}{2}},x_{j+\frac{3}{2}}], 
\end{cases}$, then $\widetilde{u}_{j}$ is the double cell average \cite{SISC_2012_ZLS_fvm}.  Similar discussion can also be seen in \cite{MCOM_2026_Yu_Cheng_Liu_Shu}. In this paper, we will take the kernel as the normalized quadratic B-spline and the definition will be given in Section~\ref{subsec:tric-fv 1D}.

Furthermore, the tilde quantities in \eqref{1D TriC Def} are also required to be consistent. The definition of consistency is presented as follows and is similar to the discussion for hyperbolic conservation laws in \cite{C&F_2018_shi_shu} and degenerate convection-diffusion equations in \cite{MCOM_2026_Yu_Cheng_Liu_Shu}. A TriC scheme is consistent if its tilde quantities are consistent.  
\begin{definition}[Consistency]\label{1D def consistency}
	A triply conservative scheme is consistent if the tilde quantities $\widetilde{u}_{j}^{s}$, $\widetilde{f}_{j+\frac{1}{2}}^{s}$, $\widetilde{A}_{j+\frac{1}{2}}^{s}$ satisfy the following conditions:
	\begin{enumerate}
		\item If $u_{h}^{s}(x)\equiv u$ for all $x\in B_{j+\frac{1}{2}}$, we have:
		\begin{equation}\label{1D consistency 1}
			\widetilde{u}_{j}^{s} = u,\quad \widetilde{f}_{j+\frac{1}{2}}^{s} = f(u),\quad \widetilde{A}_{j+\frac{1}{2}}^{s} = A(u).
		\end{equation}
		
		\item The tilde functions are all Lipschitz continuous with respect to the numerical solution, which means that there is a positive constant $C>0$ such that
		\begin{equation}\label{1D consistency 2}
			\left| \widetilde{u}_{j}^{s} - \widetilde{v}_{j}^{s} \right| + \left|\widetilde{f}_{j+\frac{1}{2}}^{s} (u_{h}^{s}) - \widetilde{f}_{j+\frac{1}{2}}^{s}(v_{h}^{s}) \right| +\left|\widetilde{A}_{j+\frac{1}{2}}^{s} (u_{h}^{s}) - \widetilde{A}_{j+\frac{1}{2}}^{s}(v_{h}^{s}) \right| \leq C\| u_{h}^{s} - v_{h}^{s} \|_{L^{\infty}(B_{j+\frac{1}{2}})} ,
		\end{equation}
		for any two functions $u_{h}^{s}$ and $v_{h}^{s}$ in the space of numerical solutions and for all integers $j$ and $s>0$.

		\item The locally conservative quantity $\widetilde{u}_{j}^{s}$ and local volume $V_{j}$ could recover the global conservation:
		\begin{equation}\label{1D consistency 3}
			\sum_{j}V_{j} \widetilde{u}_{j}^{s}= \int_{\mathbb{R}}u_{h}^{0}(x)dx, \quad \forall s\geq 0.
		\end{equation}
	\end{enumerate}
\end{definition}

\subsection{Lax-Wendroff type theorem}
In this section, we prove a new Lax-Wendroff type theorem for the TriC scheme for solving \eqref{1D TC model}. First, assume that the approximation of the initial value is accurate enough in the following sense: 
\begin{equation}
	\int_{\mathbb{R}}(u_{0} - u_{h}^{0})\phi(x,0) dx \to 0,\ \text{as}\ h\to 0,\ \forall \phi\in C_{0}^{\infty}(\mathbb{R}^{2}). 
\end{equation}
We also need the bounded total variation assumption, i.e.,
\begin{equation}
	\text{TV}(u_{h}^{s})<C,\quad \forall s\geq0,
\end{equation}
where
\begin{equation}
	\text{TV}(u_{h}^{s}) \triangleq \mathop{\text{sup}}\limits_{\Xi}\sum_{j = -\infty}^{\infty}\left|u_{h}^{s}(x_{j+\frac{1}{2}}) - u_{h}^{s}(x_{j-\frac{1}{2}})\right|,
\end{equation}
where $\Xi$ denotes the set of all partitions of the real line $\mathbb{R}$. Based on this condition, we will get an important property for our proof, i.e.,
\begin{equation}\label{tv bounded}
	h\sum_{j}\sup_{x\in B_{j+\frac{1}{2}}}|u_{h}^{s}(x) - u_{h}^{s}(x_{j+\frac{1}{2}})| \to 0,\ \text{as}\ h\to 0,\forall s\geq 0.
\end{equation}

\begin{theorem}[Lax--Wendroff type theorem]\label{Lax-Wendroff thm}
Let $\{u_{h}\}$ be the numerical solutions generated by a consistent triply conservative scheme satisfying the above assumptions. If $u_h$ converges almost everywhere to $u$ as $\tau,h\to 0$, then $u$ is an ultra-weak solution of \eqref{1D TC model}.
\end{theorem}

\begin{proof}
Define $\phi_{j}^{s}=\phi(x_{j},t_{s})$ with $\phi\in C_{0}^{\infty}(\mathbb{R}^{2})$. In this proof, we shall omit the superscript $s$ and the variable $t$ in $\phi(x,t)$ for simplicity whenever no confusion can arise.
As in the proofs of previous Lax-Wendroff type theorems, see \cite{Lax_Wendroff_1960,C&F_2018_shi_shu,MCOM_2026_Yu_Cheng_Liu_Shu}, we multiply \eqref{1D TriC Def} by $V_{j} \Delta t_{s}\phi_{j}$ and take the summation by parts to transfer the difference structure onto $\phi$, and finally use the regularity of $\phi$ to approximate the weak formulation.

For the temporal and the convective terms, the proof is basically the same as those in \cite{Lax_Wendroff_1960, C&F_2018_shi_shu, MCOM_2026_Yu_Cheng_Liu_Shu}. 
Now we mainly consider the dispersive term.
Define a auxiliary linear functional 
\begin{equation}\label{lambdaj}
	\lambda_{j} = 2\phi_{j} - \frac{1}{2}\left(\phi_{j+\frac{1}{2}} + \phi_{j-\frac{1}{2}}\right).
\end{equation}
Notice that by Taylor expansion, $\lambda_{j} = \phi_{j} + \mathcal{O}(h^{2})$.
We will first prove that:
 \begin{equation}\label{eq:2.11}
 	\small
 	\begin{aligned} 
 	&\sum_{j}\left( \frac{2}{\Delta x_{j}+\Delta x_{j+1}} \left( \frac{\widetilde{A}_{j+\frac{3}{2}} - \widetilde{A}_{j+\frac{1}{2}}}{\Delta x_{j+1}} - \frac{\widetilde{A}_{j+\frac{1}{2}} - \widetilde{A}_{j-\frac{1}{2}}}{\Delta x_{j}} \right) - \frac{2}{\Delta x_{j-1}+\Delta x_{j}} \left( \frac{\widetilde{A}_{j+\frac{1}{2}} - \widetilde{A}_{j-\frac{1}{2}}}{\Delta x_{j}} - \frac{\widetilde{A}_{j-\frac{1}{2}} - \widetilde{A}_{j-\frac{3}{2}}}{\Delta x_{j-1}} \right) \right) \lambda_{j}\\
 	&\to -\int_{\mathbb{R}} A(u)\phi_{xxx} dx.
 	\end{aligned} 
 \end{equation}
Then we will prove that 
 \begin{equation}\label{eq:2.12}
	\small
	\begin{aligned} 
		&\sum_{j}\left( \frac{2}{\Delta x_{j}+\Delta x_{j+1}} \left( \frac{\widetilde{A}_{j+\frac{3}{2}} - \widetilde{A}_{j+\frac{1}{2}}}{\Delta x_{j+1}} - \frac{\widetilde{A}_{j+\frac{1}{2}} - \widetilde{A}_{j-\frac{1}{2}}}{\Delta x_{j}} \right) - \frac{2}{\Delta x_{j-1}+\Delta x_{j}} \left( \frac{\widetilde{A}_{j+\frac{1}{2}} - \widetilde{A}_{j-\frac{1}{2}}}{\Delta x_{j}} - \frac{\widetilde{A}_{j-\frac{1}{2}} - \widetilde{A}_{j-\frac{3}{2}}}{\Delta x_{j-1}} \right) \right) (\phi_{j} - \lambda_{j})\\
		&\to 0.
	\end{aligned} 
\end{equation}
By summation by parts in \eqref{eq:2.11}, we have
 \begin{equation}\label{SBP 1D}
 	\small
 	\begin{aligned}
 			&\sum_{j}\frac{2}{\Delta x_{j}+\Delta x_{j+1}} \left( \frac{\widetilde{A}_{j+\frac{3}{2}} - \widetilde{A}_{j+\frac{1}{2}}}{\Delta x_{j+1}} - \frac{\widetilde{A}_{j+\frac{1}{2}} - \widetilde{A}_{j-\frac{1}{2}}}{\Delta x_{j}} \right) - \frac{2}{\Delta x_{j-1}+\Delta x_{j}} \left( \frac{\widetilde{A}_{j+\frac{1}{2}} - \widetilde{A}_{j-\frac{1}{2}}}{\Delta x_{j}} - \frac{\widetilde{A}_{j-\frac{1}{2}} - \widetilde{A}_{j-\frac{3}{2}}}{\Delta x_{j-1}} \right) \lambda_{j}\\
 			=&-\sum_{j}\widetilde{A}_{j+\frac{1}{2}} \bigg( \frac{1}{\Delta x_{j+1}}
 			\left(
 			\frac{2}{\Delta x_j+\Delta x_{j+1}}(\lambda_j-\lambda_{j+1})
 			-
 			\frac{2}{\Delta x_{j+1}+\Delta x_{j+2}}(\lambda_{j+1}-\lambda_{j+2})
 			\right) \\
 			&- \frac{1}{\Delta x_{j}}
 			\left(
 			\frac{2}{\Delta x_{j-1}+\Delta x_{j}}(\lambda_{j-1}-\lambda_{j})
 			-
 			\frac{2}{\Delta x_{j}+\Delta x_{j+1}}(\lambda_{j}-\lambda_{j+1})
 			\right)  \bigg) .
 	\end{aligned}
 \end{equation}
 Substituting \eqref{lambdaj} into \eqref{SBP 1D} and applying routine Taylor expansion, we obtain the following expression for the coefficient of $\widetilde{A}_{j+\frac{1}{2}}$:
 \begin{equation}\label{eq:2.13}
 	\small
 	\begin{aligned}
 		&\bigg( \frac{1}{\Delta x_{j+1}}
 		\left(
 		\frac{2}{\Delta x_j+\Delta x_{j+1}}(\lambda_j-\lambda_{j+1})
 		-
 		\frac{2}{\Delta x_{j+1}+\Delta x_{j+2}}(\lambda_{j+1}-\lambda_{j+2})
 		\right) \\
 		&- \frac{1}{\Delta x_{j}}
 		\left(
 		\frac{2}{\Delta x_{j-1}+\Delta x_{j}}(\lambda_{j-1}-\lambda_{j})
 		-
 		\frac{2}{\Delta x_{j}+\Delta x_{j+1}}(\lambda_{j}-\lambda_{j+1})
 		\right)  \bigg)= \omega_{j+\frac{1}{2}} \phi_{xxx}(x_{j+\frac{1}{2}})+ o(1),
 	\end{aligned}
 \end{equation} 
where 
\begin{equation}\label{1Deq:omega}
	\small
		\begin{aligned}
			\omega_{j+\frac{1}{2}}
			=&
			\frac{\Delta x_{j-1}}{12}-\frac{(\Delta x_{j-1})^2}{12\Delta x_j}
			+\frac{(\Delta x_j)^2}{12\Delta x_{j+1}}
			+\frac{5\Delta x_j}{12}
			+\frac{5\Delta x_{j+1}}{12}
			+\frac{\Delta x_{j+2}}{12}-
			\frac{(\Delta x_{j+2})^2}{12\Delta x_{j+1}}
			+
			\frac{(\Delta x_{j+1})^2}{12\Delta x_j}\\
			= & \frac{\Delta x_j+\Delta x_{j+1}}2+\eta_{j+1}-\eta_j ,
	\end{aligned}
\end{equation}
with
\begin{equation}\label{1Deq:eta-def}
	\small
		\eta_j \triangleq \frac{(\Delta x_{j-1})^2-(\Delta x_{j+1})^2}{12\Delta x_j}
		+\frac{\Delta x_{j+1}-\Delta x_{j-1}}{12} = \mathcal{O}(h).
\end{equation}
Together with \eqref{SBP 1D}--\eqref{1Deq:eta-def}, we have
\begin{equation}\label{eq2.16}
	\small
	\begin{aligned}
		&\sum_{j}\widetilde{A}_{j+\frac{1}{2}} \bigg( \frac{1}{\Delta x_{j+1}}
		\left(
		\frac{2}{\Delta x_j+\Delta x_{j+1}}(\lambda_j-\lambda_{j+1})
		-
		\frac{2}{\Delta x_{j+1}+\Delta x_{j+2}}(\lambda_{j+1}-\lambda_{j+2})
		\right) \\
		&- \frac{1}{\Delta x_{j}}
		\left(
		\frac{2}{\Delta x_{j-1}+\Delta x_{j}}(\lambda_{j-1}-\lambda_{j})
		-
		\frac{2}{\Delta x_{j}+\Delta x_{j+1}}(\lambda_{j}-\lambda_{j+1})
		\right)  \bigg) \\
		=&\sum_{j}\widetilde{A}_{j+\frac{1}{2}}\omega_{j+\frac{1}{2}}\phi_{xxx}(x_{j+\frac{1}{2}})=\sum_{j}\widetilde{A}_{j+\frac{1}{2}}\phi_{xxx}(x_{j+\frac{1}{2}})\left( \frac{\Delta x_{j}+\Delta x_{j+1}}{2} + \eta_{j+1} - \eta_{j} \right)\\
		=&\sum_{j}\widetilde{A}_{j+\frac{1}{2}}\phi_{xxx}(x_{j+\frac{1}{2}})\left( \frac{\Delta x_{j}+\Delta x_{j+1}}{2} \right)  + \widetilde{A}_{j+\frac{1}{2}}\phi_{xxx}(x_{j+\frac{1}{2}}) \left( \eta_{j+1} - \eta_{j} \right)\\
		=&\sum_{j}\widetilde{A}_{j+\frac{1}{2}}\phi_{xxx}(x_{j+\frac{1}{2}})\left( \frac{\Delta x_{j}+\Delta x_{j+1}}{2} \right) - \left(\widetilde{A}_{j+\frac{1}{2}} \phi_{xxx}(x_{j+\frac{1}{2}}) - \widetilde{A}_{j-\frac{1}{2}} \phi_{xxx}(x_{j-\frac{1}{2}}) \right) \eta_{j}.
	\end{aligned}
\end{equation}
Notice that $\eta_{j} = \mathcal{O}(h)$ and 
$$
\begin{aligned} 
\widetilde{A}_{j+\frac{1}{2}} \phi_{xxx}(x_{j+\frac{1}{2}}) - \widetilde{A}_{j-\frac{1}{2}} \phi_{xxx}(x_{j-\frac{1}{2}}) = \left( \widetilde{A}_{j+\frac{1}{2}} - \widetilde{A}_{j-\frac{1}{2}}  \right) \phi_{xxx}(x_{j+\frac{1}{2}}) + \widetilde{A}_{j-\frac{1}{2}} \left( \phi_{xxx}(x_{j+\frac{1}{2}}) - \phi_{xxx}(x_{j-\frac{1}{2}}) \right),
\end{aligned}
$$
by using the fact that $\phi_{xxx}(x_{j+\frac{1}{2}}) - \phi_{xxx}(x_{j-\frac{1}{2}}) = \mathcal{O}(h)$ and the bounded variation assumption \eqref{tv bounded}, we can conclude the magnitude of the last term of the last equation in \eqref{eq2.16} is $\mathcal{O}(h)$. Therefore we can conclude that 
$$
	\small
		\begin{aligned}
		&\sum_{j} \bigg(  \widetilde{A}_{j+\frac{1}{2}} \bigg( \frac{1}{\Delta x_{j+1}}
		\left(
		\frac{2}{\Delta x_j+\Delta x_{j+1}}(\lambda_j-\lambda_{j+1})
		-
		\frac{2}{\Delta x_{j+1}+\Delta x_{j+2}}(\lambda_{j+1}-\lambda_{j+2})
		\right) \\
		&- \frac{1}{\Delta x_{j}}
		\left(
		\frac{2}{\Delta x_{j-1}+\Delta x_{j}}(\lambda_{j-1}-\lambda_{j})
		-
		\frac{2}{\Delta x_{j}+\Delta x_{j+1}}(\lambda_{j}-\lambda_{j+1})
		\right)  \bigg) \bigg) \\
		=&\sum_{j}\widetilde{A}_{j+\frac{1}{2}}\phi_{xxx}(x_{j+\frac{1}{2}})\left( \frac{\Delta x_{j}+\Delta x_{j+1}}{2} \right) + \mathcal{O}(h).
	\end{aligned}
$$
Then we can conclude that $\sum_{j}\widetilde{A}_{j+\frac{1}{2}}\phi_{xxx}(x_{j+\frac{1}{2}})\left( \frac{\Delta x_{j}+\Delta x_{j+1}}{2} \right) \to \int_{\mathbb{R}} A(u)\phi_{xxx}(x) dx $ as $h,\tau \to 0$. 

Next, we will prove \eqref{eq:2.12}. 
By the discrete scheme and summation by parts in time, we obtain
By \eqref{1D TriC Def},
\begin{equation}\label{eq:2.18}
	\begin{aligned}
		&\Delta t_s\sum_j\Bigg[
		\frac{2}{\Delta x_j+\Delta x_{j+1}}
		\left(\frac{\widetilde A_{j+\frac32}^s-\widetilde A_{j+\frac12}^s}{\Delta x_{j+1}}
		-\frac{\widetilde A_{j+\frac12}^s-\widetilde A_{j-\frac12}^s}{\Delta x_j}\right)
		\\
		&		-\frac{2}{\Delta x_{j-1}+\Delta x_j}
		\left(\frac{\widetilde A_{j+\frac12}^s-\widetilde A_{j-\frac12}^s}{\Delta x_j}
		-\frac{\widetilde A_{j-\frac12}^s-\widetilde A_{j-\frac32}^s}{\Delta x_{j-1}}\right)
		\Bigg](\phi_j^s-\lambda_j^s)
		\\
		=&-\sum_jV_j(\widetilde u_j^{s+1}-\widetilde u_j^s)
		(\phi_j^s-\lambda_j^s)-\Delta t_s\sum_j
		(\widetilde f_{j+\frac12}^s-\widetilde f_{j-\frac12}^s)
		(\phi_j^s-\lambda_j^s).
	\end{aligned}
\end{equation}
Combined with the bounded variation assumption \eqref{tv bounded}, $\lambda_{j}-\phi_{j} = \mathcal{O}(h^{2})$, under the standard analysis procedure in the literature \cite{Lax_Wendroff_1960,C&F_2018_shi_shu,MCOM_2026_Yu_Cheng_Liu_Shu}, we finish the proof.

\end{proof}

We next illustrate how to construct or identify a TriC scheme. We first develop a new TriC finite volume method based on the notion of a ``triple cell average'', which will be introduced below. We then show that several existing methods in the literature, including the finite difference scheme \cite{MCS_2023_Ahmat_Qiu} and the discontinuous Galerkin schemes \cite{MCOM_2008_ultraweak,JCP_2004_Levy_Shu_Yan}, admit a TriC structure.

\subsection{A triply conservative finite volume framework}\label{subsec:tric-fv 1D}
The TriC structure can serve as a design principle for constructing a new, high-order, nonstandard finite volume method. The key issue is the choice of the locally conserved variable to be evolved. For hyperbolic conservation laws, the standard cell average arises naturally from a single local integration. For convection-diffusion equations, the double cell average is used \cite{SISC_2012_ZLS_fvm, MCOM_2026_Yu_Cheng_Liu_Shu}. For the third-order dispersive equation \eqref{1D TC model}, the analogous construction requires the triple cell average. This leads us to integrate \eqref{1D TC model} three times, which also reveals an equivalent representation in terms of an inner product with the quadratic B-spline function \cite{Book_de_Boor}.
These two definitions lead to equivalent numerical formulations on uniform meshes. 
On non-uniform meshes, the inner-product formulation provides a more convenient definition for defining the triple average and constructing the TriC-FV scheme.

\subsubsection{Triple cell averages and quadratic B-splines} 
We will illustrate the motivation for introducing the triple cell average. We start from a uniform mesh of size $\Delta x$ and show the triple-integral form of triple cell averages. 
Integrating the equation \eqref{1D TC model} three times in space, we have
\begin{equation}\label{1D Tric uniform}
	\begin{aligned} 
	&\frac{d}{dt}\int_{x_{j-\frac{1}{2}}}^{x_{j+\frac{1}{2}}} \int_{\xi-\frac{\Delta x}{2}}^{\xi+\frac{\Delta x}{2}}\int_{\eta-\frac{\Delta x}{2}}^{\eta+\frac{\Delta x}{2}} u(x,t)dx d\eta d\xi + \int_{x_{j-\frac{1}{2}}}^{x_{j+\frac{1}{2}}} \int_{\xi-\frac{\Delta x}{2}}^{\xi+\frac{\Delta x}{2}} f(u(\eta+\frac{\Delta x}{2},t)) - f(u(\eta-\frac{\Delta x}{2},t))d\eta d\xi \\
	=& -\left(  A(u(x_{j+\frac{3}{2}},t)) - 3 A(u(x_{j+\frac{1}{2}},t)) + 3 A(u(x_{j-\frac{1}{2}},t)) - A(u(x_{j-\frac{3}{2}},t)) \right).
\end{aligned}
\end{equation}
The above equation apparently lies in the TriC form if we multiply both sides by the scale factor $\frac{1}{(\Delta x)^{3}}$. 
The triple-integral form of the triple cell average in $I_j$ is then defined as follows naturally:
\begin{equation}\label{1D triple average}
	\bar{\bar{\bar{u}}}_{j}(t)
	=
	\frac{1}{(\Delta x)^3}
	\int_{x_{j-\frac{1}{2}}}^{x_{j+\frac{1}{2}}} \int_{\xi-\frac{\Delta x}{2}}^{\xi+\frac{\Delta x}{2}}\int_{\eta-\frac{\Delta x}{2}}^{\eta+\frac{\Delta x}{2}} u(x,t)dx d\eta d\xi .
\end{equation}
We can therefore construct a nonstandard finite volume scheme based on \eqref{1D Tric uniform} and \eqref{1D triple average}.
However, the above procedure relies heavily on the assumption of the uniform mesh. 
Inspired by the doubly conservative schemes for degenerate convection-diffusion equations in the recent works \cite{MCOM_2026_Yu_Cheng_Liu_Shu,arxiv_2026_Yu_Cheng_Shu}, through exchanging the order of the threefold integral \eqref{1D triple average}, we obtain the following equivalence between the triple cell average and the quadratic B-spline.
\begin{proposition}
	Let $x_{j}=(x_{j-\frac{1}{2}}+x_{j+\frac{1}{2}})/2$ be the midpoint of the interval $I_{j}$. For a locally integrable function $q$, define its triple cell average of $I_j$ by \eqref{1D triple average}. Then we have
\begin{equation}\label{1D triple-average-bspline} 
		\bar{\bar{\bar{q}}}_j = \frac{1}{\Delta x} \int_{\mathbb{R}} q(x) B_2\!\left(\frac{x-x_j}{\Delta x}\right) dx, 
\end{equation} 
where $B_2(\hat{x})$ is the centered quadratic B-spline which is defined on the reference interval by
\begin{equation}\label{eq:quadratic-bspline}
	\small
	B_2(\hat{x})
	=
	\begin{cases}
		\displaystyle
		\frac{1}{2}\left(\frac{3}{2}+\hat{x} \right)^2,
		& -\frac{3}{2}\le \hat{x}<-\frac{1}{2},
		\\
		\frac{3}{4}-\hat{x}^2,
		& -\frac{1}{2}\le \hat{x}\le\frac{1}{2},
		\\
		\displaystyle
		\frac{1}{2}\left(\frac{3}{2}-\hat{x}\right)^2,
		& \frac{1}{2} < \hat{x}\le\frac{3}{2},
		\\
		0,
		& \text{otherwise}.
	\end{cases}
\end{equation}
\end{proposition}
 
\begin{proof} 
	Since $q\in L^1_{\mathrm{loc}}(\mathbb{R})$, exchanging the order of integration in \eqref{1D triple average}, we obtain
	\begin{align}
		(\Delta x)^3 \bar{\bar{\bar{q}}}_j
		=&
		\int_{x_{j-\frac{3}{2}}}^{x_{j-\frac{1}{2}}}
		q(x)
		\int_{x_{j-\frac{1}{2}}}^{x+\Delta x}
		\int_{\eta-\frac{\Delta x}{2}}^{x+\frac{\Delta x}{2}}
		d\xi d\eta dx
		+
		\int_{x_{j-\frac{1}{2}}}^{x_{j+\frac{1}{2}}}
		q(x)
		\Bigg[
		\int_{x_{j-\frac{1}{2}}}^{x}
		\int_{x-\frac{\Delta x}{2}}^{\eta+\frac{\Delta x}{2}}
		d\xi d\eta
		+
		\int_{x}^{x_{j+\frac{1}{2}}}
		\int_{\eta-\frac{\Delta x}{2}}^{x+\frac{\Delta x}{2}}
		d\xi d\eta
		\Bigg]dx
		\nonumber\\
		&+
		\int_{x_{j+\frac{1}{2}}}^{x_{j+\frac{3}{2}}}
		q(x)
		\int_{x-\Delta x}^{x_{j+\frac{1}{2}}}
		\int_{x-\frac{\Delta x}{2}}^{\eta+\frac{\Delta x}{2}}
		d\xi d\eta dx.
		\label{eq:triple-average-reordered}
	\end{align}
	Evaluating the inner two integrals in the three regions gives
	\begin{equation}\label{eq:triple-average-piecewise}
		\begin{aligned}
			(\Delta x)^3 \bar{\bar{\bar{q}}}_j
			=&
			\int_{x_{j-\frac{3}{2}}}^{x_{j-\frac{1}{2}}}
			q(x)
			\frac{1}{2}
			\left(
			x-x_j+\frac{3\Delta x}{2}
			\right)^2 dx
			+
			\int_{x_{j-\frac{1}{2}}}^{x_{j+\frac{1}{2}}}
			q(x)
			\left[
			\frac{3}{4}(\Delta x)^2-(x-x_j)^2
			\right] dx
			\\
			&+
			\int_{x_{j+\frac{1}{2}}}^{x_{j+\frac{3}{2}}}
			q(x) 
			\frac{1}{2}
			\left(
			x_j+\frac{3\Delta x}{2}-x
			\right)^2 dx.
		\end{aligned}
	\end{equation}
	Dividing by $(\Delta x)^3$ and comparing the three piecewise weights with
	\eqref{eq:quadratic-bspline} yields
	\[
	\bar{\bar{\bar{q}}}_j
	=
	\frac{1}{\Delta x}
	\int_{\mathbb{R}}
	q(x)
	B_2 \left(\frac{x-x_j}{\Delta x}\right) dx,
	\]
	which proves \eqref{1D triple-average-bspline}.
\end{proof}

Therefore, on uniform meshes, the triple cell average is exactly a weighted average of $u$ with respect to the normalized quadratic B-spline kernel $\frac{1}{\Delta x}B_{2}\left(\frac{x-x_{j}}{\Delta x}\right)$.
For non-uniform meshes, repeating integration in \eqref{1D triple average} no longer works.
Therefore we define triple cell averages through the inner form \eqref{1D triple-average-bspline}.
Denote $B_{2,j}(x)$ as the quadratic spline function based on the intervals $\overline{I_{j-1}\cup I_j\cup I_{j+1}}$.
The expression of $B_{2,j}$ is as follows:
\begin{equation}\label{non-uniform-B2-piecewise}
	B_{2,j}(x)=
	\begin{cases}
		\displaystyle
		\frac{(x-x_{j-\frac{3}{2}})^2}
		{\Delta x_{j-1}(\Delta x_{j-1}+\Delta x_{j})},
		&x\in I_{j-1},\\[2ex]
		\displaystyle
		\frac{(x-x_{j-\frac{3}{2}})(x_{j+\frac{1}{2}}-x)}{(\Delta x_{j-1}+\Delta x_{j})\Delta x_{j}}
		+
		\frac{(x_{j+\frac{3}{2}}-x)(x-x_{j-\frac{1}{2}})}{(\Delta x_{j}+\Delta x_{j+1})\Delta x_{j}},
		&x\in I_j,\\[2ex]
		\displaystyle
		\frac{(x_{j+\frac{3}{2}}-x)^2}{\Delta x_{j+1}(\Delta x_{j}+\Delta x_{j+1})},
		&x\in I_{j+1},\\[2ex]
		\displaystyle
		0,&\text{otherwise}.
	\end{cases}
\end{equation}
It is clear that the B-splines are nonnegative and satisfy the identity
\begin{equation}\label{1D:B2-partition-unity}
	\sum_j B_{2,j}(x)=1.
\end{equation}
Moreover, $V_j	\triangleq 	\int_{\mathbb{R}}B_{2,j}(x) dx = \left(\Delta x_{j-1}+\Delta x_j+\Delta x_{j+1}\right)/3$.
We therefore define the non-uniform triple cell average by
\begin{equation}\label{1D:non-uniform-triple-average}
	\bar{\bar{\bar{u}}}_{j}(t)
	=
	\frac{1}{V_j}
	\int_{\mathbb{R}}u(x,t)B_{2,j}(x) dx.
\end{equation}
When the mesh is uniform, $V_j=\Delta x$ and
\eqref{1D:non-uniform-triple-average} reduces exactly to
\eqref{1D triple average}. In addition, \eqref{1D:B2-partition-unity} gives
\begin{equation}
	\sum_j V_j \bar{\bar{\bar{u}}}_j(t)
	=
	\int_{\mathbb{R}}u(x,t) dx.
\end{equation}

\subsubsection{The TriC-FV scheme}\label{subsubsec:tric-fv-scheme}
We now construct the TriC-FV scheme on non-uniform meshes that evolves the triple cell average based on its inner-form. The construction follows naturally from the quadratic B-spline representation of the triple cell average. Additionally, the first derivative of $B_{2,j}$ produces a conservative difference associated with the convective term, whereas its third derivative produces exactly the third-order conservative difference required by the TriC structure.
Let $B_{1,j+\frac{1}{2}}$ denote the linear $B$-spline supported on $I_j\cup I_{j+1}$, defined by
\begin{equation}
	B_{1,j+\frac{1}{2}}(x)
	=
	\begin{cases}
		\displaystyle
		\frac{x-x_{j-\frac{1}{2}}}{\Delta x_j},
		& x\in I_j,\\[2ex]
		\displaystyle
		\frac{x_{j+\frac{3}{2}}-x}{\Delta x_{j+1}},
		& x\in I_{j+1},\\[2ex]
		0,
		& \text{otherwise}.
	\end{cases}
\end{equation}
A direct differentiation of $B_{2,j}$ in \eqref{non-uniform-B2-piecewise} gives
\begin{equation}\label{1D:B2-first-derivative}
	B_{2,j}'(x)
	=
	\frac{2}{\Delta x_{j-1}+\Delta x_j}
	B_{1,j-\frac{1}{2}}(x)
	-
	\frac{2}{\Delta x_j+\Delta x_{j+1}}
	B_{1,j+\frac{1}{2}}(x).
\end{equation}
Since $B_{2,j}$ is $C^1$ and piecewise quadratic, its third derivative is understood in the distributional sense.

We first derive the corresponding semi-discrete formulation. Multiplying \eqref{1D TC model} by $B_{2,j}$, integrating over $\mathbb{R}$, and formally integrating by parts yield
\begin{equation}\label{1D:B2-tested-balance}
	\frac{d}{dt}
	\int_{\mathbb{R}}
	u(x,t)B_{2,j}(x) dx
	-
	\int_{\mathbb{R}}
	f(u(x,t))B_{2,j}'(x) dx
	-
	\left\langle
	A(u(\cdot,t)),D_x^3B_{2,j}
	\right\rangle
	=0.
\end{equation}
The brackets denote inner product in the distributional sense. 
No boundary terms arise since $B_{2,j}$ is compactly supported. By the definition of the triple cell average in \eqref{1D:non-uniform-triple-average},
\begin{equation}
	\frac{d}{dt}
	\int_{\mathbb{R}}
	u(x,t)B_{2,j}(x) dx
	=
	V_j
	\frac{d}{dt}
	\bar{\bar{\bar{u}}}_j(t).
\end{equation}
For the convective term, substituting \eqref{1D:B2-first-derivative} into \eqref{1D:B2-tested-balance} gives
\begin{equation}
	\begin{aligned}
		\int_{\mathbb{R}}
		f(u)B_{2,j}'dx
		=&
		\frac{2}{\Delta x_{j-1}+\Delta x_j}
		\int_{I_{j-1}\cup I_j}
		f(u)B_{1,j-\frac{1}{2}}dx -\frac{2}{\Delta x_j+\Delta x_{j+1}}
		\int_{I_j\cup I_{j+1}}
		f(u)B_{1,j+\frac{1}{2}} dx .
	\end{aligned}
\end{equation}
Therefore, the natural generalized flux associated with the interface $x_{j+\frac{1}{2}}$ is
\begin{equation}\label{1D:tric-convective-flux-continuous}
	\widetilde{f}_{j+\frac{1}{2}}(t) = \frac{2}{\Delta x_j+\Delta x_{j+1}}
	\int_{I_j\cup I_{j+1}}
	f(u(x,t))B_{1,j+\frac{1}{2}}(x)dx .
\end{equation}
For the nonlinear dispersive term, the third distributional derivative of $B_{2,j}$ gives
\begin{equation}\label{1D:B2-third-derivative-pairing}
	\begin{aligned}
		-
		\left\langle
		A(u),D_x^3B_{2,j}
		\right\rangle 
		=&
		\frac{2}{\Delta x_j+\Delta x_{j+1}}
		\left(
		\frac{
			A(u(x_{j+\frac{3}{2}},t))
			-
			A(u(x_{j+\frac{1}{2}},t))
		}{\Delta x_{j+1}}
		-
		\frac{
			A(u(x_{j+\frac{1}{2}},t))
			-
			A(u(x_{j-\frac{1}{2}},t))
		}{\Delta x_j}
		\right)
		\\
		&-
		\frac{2}{\Delta x_{j-1}+\Delta x_j}
		\left(
		\frac{
			A(u(x_{j+\frac{1}{2}},t))
			-
			A(u(x_{j-\frac{1}{2}},t))
		}{\Delta x_j}
		-
		\frac{
			A(u(x_{j-\frac{1}{2}},t))
			-
			A(u(x_{j-\frac{3}{2}},t))
		}{\Delta x_{j-1}}
		\right).
	\end{aligned}
\end{equation}
Hence, the third-order spatial derivative is converted exactly into the needed TriC formulation.

To get a fully-discrete scheme, as we mainly consider the spatial discretization in this paper, we use the forward Euler scheme to discretize the temporal derivative. Furthermore, at the time level $t=t_{s}$, we will reconstruct a high-order numerical solution based on the triple cell average and introduce the numerical fluxes $\widetilde{f}_{j+\frac{1}{2}}^{s}$ and $\widetilde{A}_{j+\frac{1}{2}}^{s}$ to approximate the right-hand terms in \eqref{1D:tric-convective-flux-continuous}-- \eqref{1D:B2-third-derivative-pairing}. Their concrete reconstruction procedure will be specified later. 
Then the resulting fully-discrete TriC-FV scheme is
\begin{equation}\label{1D:tric-fv-forward-euler}
	\begin{aligned}
		\frac{
			\bar{\bar{\bar{u}}}_j^{s+1}
			-
			\bar{\bar{\bar{u}}}_j^s
		}
		{\Delta t_s}
		= &
		-
		\frac{
			\widetilde{f}_{j+\frac{1}{2}}^s
			-
			\widetilde{f}_{j-\frac{1}{2}}^s
		}
		{V_j} - \frac{1}{V_j}
		\frac{2}{\Delta x_j+\Delta x_{j+1}}
		\left(	\frac{
			\widetilde{A}_{j+\frac{3}{2}}^s
			-\widetilde{A}_{j+\frac{1}{2}}^s
		}{\Delta x_{j+1}}
		-	\frac{	\widetilde{A}_{j+\frac{1}{2}}^s
			-\widetilde{A}_{j-\frac{1}{2}}^s		}{\Delta x_j}
		\right)
		\\
		&+
		\frac{1}{V_j}
		\frac{2}{\Delta x_{j-1}+\Delta x_j}
		\left(
		\frac{
			\widetilde{A}_{j+\frac{1}{2}}^s
			-
			\widetilde{A}_{j-\frac{1}{2}}^s
		}{\Delta x_j}
		-
		\frac{
			\widetilde{A}_{j-\frac{1}{2}}^s
			-
			\widetilde{A}_{j-\frac{3}{2}}^s
		}{\Delta x_{j-1}}
		\right),
	\end{aligned}
\end{equation}
where  $\bar{\bar{\bar{u}}}_j^s$ is the triple cell average of $I_{j}$ at time $t_{s}$.
Comparing \eqref{1D:tric-fv-forward-euler} with \eqref{1D TriC Def}, we see that the generalized locally conserved variable is precisely the triple cell average,
\begin{equation*}
	\widetilde{u}_j^s
	=
	\bar{\bar{\bar{u}}}_j^s
	=
	\frac{1}{V_j}
	\int_{\mathbb{R}}
	u_h^s(x)B_{2,j}(x) dx,
\end{equation*}
with the corresponding local volume
\begin{equation}\label{1D:tric-fv-local-volume}
	V_j
	=
	\frac{
		\Delta x_{j-1}
		+
		\Delta x_j
		+
		\Delta x_{j+1}
	}{3}.
\end{equation}
Therefore, once $\widetilde{f}_{j+\frac{1}{2}}^s$ and $\widetilde{A}_{j+\frac{1}{2}}^s$ locally depend on the numerical solution, the scheme \eqref{1D:tric-fv-forward-euler} is a TriC scheme.

We now specify the reconstruction and implementation of the tilde functions. Unlike the traditional finite volume schemes, the reconstruction proposed here is  based on the triple cell average $\bar{\bar{\bar{u}}}_{j}^{s}$ rather than the standard cell averages. More precisely, we construct a local polynomial $u_{j}^{s}(x)$ whose triple cell averages match the given data on $I_j$ and on a fixed number of neighboring cells. Hence this reconstructed polynomial depends only on local information. The detailed construction will be given in the next subsection.

For the convective term, the generalized flux $\tilde{f}_{j+\frac{1}{2}}(t)$ in \eqref{1D:tric-convective-flux-continuous} involves an integral over $I_j\cup I_{j+1}$ and the numerical flux $\tilde{f}_{j+\frac{1}{2}}^{s}$ will be computed by numerical quadrature. Let $x_{j,\alpha}$ and $\omega_{j,\alpha}$ denote the quadrature points and the corresponding quadrature weights on a cell $I_{j}$. We define
\begin{equation}\label{1D:tric-convective-quadrature}
	\begin{aligned}
		\widetilde{f}_{j+\frac{1}{2}}^s
		=
		\frac{2}{\Delta x_j+\Delta x_{j+1}}
		\bigg(
		&
		\sum_{\alpha}
		\omega_{j,\alpha}
		\hat{f}_{j+\frac{1}{2}}^{s}\left(u_{j}^{s}(x_{j,\alpha}),u_{j+1}^{s}(x_{j,\alpha})\right)
		B_{1,j+\frac{1}{2}}(x_{j,\alpha})\\
		+&
		\sum_{\alpha}
		\omega_{j+1,\alpha}
		\hat{f}_{j+\frac{1}{2}}^{s}\left(u_{j}^{s}(x_{j+1,\alpha}),u_{j+1}^{s}(x_{j+1,\alpha})\right)
		B_{1,j+\frac{1}{2}}(x_{j+1,\alpha})
		\bigg),
	\end{aligned}
\end{equation}
where $\hat{f}_{j+\frac{1}{2}}^{s}$ is the classical consistent, monotone and locally Lipschitz numerical flux for stability (see \cite{Lax_Wendroff_1960,Book_Leveque}), such as the Lax-Friedrichs flux:
\begin{equation}\label{eq: 1D LF flux}
	\hat{f}_{j+\frac{1}{2}}^{s}(u^{-},u^{+}) = \frac{1}{2} \left(f(u^{-})+f(u^{+})\right) - \max_{u}|f'(u)|(u^{+} - u^{-}).
\end{equation}
The quadrature rule $\{ \omega_{j,\alpha}, x_{j,\alpha} \}$ is chosen sufficiently accurate so that the desired order of accuracy of the spatial discretization is preserved.

For the dispersive term, no numerical integration is needed since \eqref{1D:B2-third-derivative-pairing} only involves point values at cell interfaces. We therefore set
\begin{equation}\label{1D:tric-dispersive-flux}
	\widetilde{A}_{j+\frac{1}{2}}^s
	=
	\hat{A}_{j+\frac{1}{2}}^{s}\left(u_{j}^{s}(x_{j+\frac{1}{2}}), u_{j+1}^{s}(x_{j+\frac{1}{2}})
	\right).
\end{equation}
For stability, $\hat{A}_{j+\frac{1}{2}}^{s}$ is chosen as a classical consistent, reverse-monotone, and locally Lipschitz continuous numerical flux (see \cite{JCP_2004_Levy_Shu_Yan}). Here, reverse monotonicity means that $\hat{A}_{j+\frac{1}{2}}^{s}(\cdot,\cdot)$ is non-increasing in its first argument and non-decreasing in its second argument. 
We adopt the reverse Godunov flux, defined by
\begin{equation}\label{1D-FV-dispersive-point-flux}
	\widehat A^{\mathrm G}(v,w)=
	\begin{cases}
		\displaystyle\max_{z\in[v,w]}A(z),&v\leq w,\\[1mm]
		\displaystyle\min_{z\in[w,v]}A(z),&v>w.
	\end{cases}
\end{equation}
All the quadrature values appearing in \eqref{1D:tric-convective-quadrature} and the interface values appearing in \eqref{1D:tric-dispersive-flux} are obtained from the same reconstruction based on triple cell averages. This reconstruction will be described in detail in the next subsection.

With the scheme defined as above, the reconstruction is local and the numerical fluxes are locally Lipschitz continuous. Together with \eqref{1D:non-uniform-triple-average}, these properties imply that the new finite volume scheme is triply conservative and consistent in the sense of Def.~\ref{1D def tric} and Def.~\ref{1D def consistency}. This yields the following proposition.

\begin{proposition}
	Assume that the reconstruction from triple cell averages is local, constant-preserving, and locally Lipschitz with respect to its input data. Assume also that the quadrature rule in \eqref{1D:tric-convective-quadrature} is exact for linear functions on each cell. Then the scheme \eqref{1D:tric-fv-forward-euler}, together with \eqref{1D:tric-convective-quadrature} and \eqref{1D:tric-dispersive-flux}, is a consistent triply conservative scheme in the sense of \eqref{1D TriC Def} and \eqref{1D def consistency}.
\end{proposition}

\subsubsection{High-order reconstruction based on triple cell averages}
In this paper, we focus on the construction of a fifth-order triply conservative finite volume scheme, for which a polynomial of degree four is needed. 
We note that the reconstruction procedure developed below can be extended to higher-order approximations.

Given triple cell averages $\left\{\bar{\bar{\bar{u}}}_{j}^{s}\right\}_{j}$, we reconstruct a fourth-order polynomial $u^s_{j}(x)$ on each cell $I_{j}$. 
Various reconstruction methods are available, such as \cite{SIREV_2009_Shu,JCP_1989_Shu_Osher,JCP_2018_Zhu_Shu_MR-WENO}, which employ nonlinear reconstruction to achieve high-order accuracy while suppressing spurious numerical oscillations. 
Since this paper focuses on resolving the correct ultra-weak solution, we find that linear reconstruction is sufficient for all the examples considered. 
Therefore, we only describe the linear reconstruction here.

Recall that the inner-product-form triple cell average associated with $B_{2,\ell}(x)$ is defined in \eqref{1D:non-uniform-triple-average}.
At a fixed time level $t_s$, write
\begin{equation}
	\mathcal{M}_\ell[q]=\frac{1}{V_\ell}
	\int_{\mathbb{R}}q(x)B_{2,\ell}(x) dx,
	\qquad
	V_\ell=\frac{\Delta x_{\ell-1}+\Delta x_\ell+\Delta x_{\ell+1}}{3}.
\end{equation}
Here $B_{2,\ell}$ is the quadratic B-spline defined in \eqref{non-uniform-B2-piecewise}. In particular, $\mathcal{M}_\ell[1]=1$. For the target cell $I_j$, we use the five-cell stencil
\begin{equation}\label{eq:large-stencil}
	\mathcal{S}_j^5=\{I_{j-2},I_{j-1},I_j,I_{j+1},I_{j+2}\}
\end{equation}
and determine $u_j^s(x)\in \mathbb{P}^4$ by
\begin{equation}
	\mathcal{M}_{j+r}[u_j^s]=\bar{\bar{\bar{u}}}_{j+r}^s,
	\qquad r=-2,-1,0,1,2.
\end{equation}
Finally, the polynomials $u_j^s$ and $u_{j+1}^s$ are used directly to obtain the point values in \eqref{1D:tric-convective-quadrature} and \eqref{1D:tric-dispersive-flux}.
It is clear that the locality, constant preservation, and Lipschitz continuity provide the consistency requirement in \eqref{1D consistency 1}-\eqref{1D consistency 3}. 

\subsection{Triply conservative structure of existing high-order schemes}\label{subsec:existing-tric}
The structure also exists in several existing high-order schemes, including the successful LDG-I scheme in \cite{JCP_2004_Levy_Shu_Yan}, the finite difference scheme in \cite{MCS_2023_Ahmat_Qiu} and the ultra-weak DG (UWDG) scheme \cite{MCOM_2008_ultraweak}.
This observation provides independent evidence that the TriC definition captures a natural property of discretizations for third-order dispersive equations rather than only an artifact of the finite volume construction.
We will show how to recognize these schemes as TriC schemes.

\subsubsection{Finite difference WENO scheme}
Consider the fifth-order direct WENO scheme of Ahmat and Qiu \cite{MCS_2023_Ahmat_Qiu}. 
For this scheme, the linear reconstruction preserves the triple conservation.
On a uniform mesh, for the linear reconstruction, its semi-discrete formulation is
\begin{equation}\label{1D-MCS-semi}
	\frac{d u_j}{dt}
	=-\frac{\widehat{F}_{j+\frac{1}{2}}-\widehat{F}_{j-\frac{1}{2}}}{\Delta x}
	-\frac{\widehat{G}_{j+\frac{1}{2}}-\widehat{G}_{j-\frac{1}{2}}}{\Delta x},
\end{equation}
where $u_j$ denotes the nodal value, and $\widehat{F}$ and $\widehat{G}$ are the convective and dispersive numerical fluxes, respectively.
The linear reconstruction is based on a triple-integral representation of the dispersive flux. Specifically, an auxiliary function $ p$ is introduced through
\begin{equation}
	A(u(x)) = \frac{1}{(\Delta x)^3}
	\int_{x-\frac{\Delta x}{2}}^{x+\frac{\Delta x}{2}}
	\int_{\eta-\frac{\Delta x}{2}}^{\eta+\frac{\Delta x}{2}}
	\int_{\xi-\frac{\Delta x}{2}}^{\xi+\frac{\Delta x}{2}}
	p(\theta) d\theta d\xi d\eta.
\end{equation}
The reconstruction polynomial $p_j(x)$ is determined by matching its normalized triple averages centered at the stencil nodes $x_\ell$ to $A(u(x_\ell))$, with the same construction applied to each split flux.
Thus the dispersive term $A(u)_{xxx}$ is approximated as:
\begin{equation}
	A(u)_{xxx} = \frac{p_{j}(x_{j+\frac{3}{2}}) - 3p_{j}(x_{j+\frac{1}{2}}) + 3p_{j}(x_{j-\frac{1}{2}}) - p_{j}(x_{j-\frac{3}{2}})}{(\Delta x)^{3}},
\end{equation}
where the corresponding $\hat{G}_{j+\frac{1}{2}}$ in \eqref{1D-MCS-semi} is taken as $G_{j+\frac{1}{2}}^{L}$ in the linear reconstruction with
\begin{equation}
	G_{j+\frac{1}{2}}^{L} = \frac{p_{j}(x_{j+\frac{3}{2}}) - 2p_{j}(x_{j+\frac{1}{2}}) + p_{j}(x_{j-\frac{1}{2}})}{(\Delta x)^{2}}.
\end{equation}
For the nonlinear weighted procedure, the situation changes. The nonlinearly reconstructed numerical flux $\hat{G}_{j+\frac{1}{2}}^{N}$ is formed by a nonlinear weighted combination of different linearly reconstructed fluxes $G_{j+\frac{1}{2}}^{L,\ell} , \ \ell=1,2,3$ on different stencils.
This nonlinear weighting breaks the linear structure underlying the three successive summation-by-parts steps, so the triply conservative form is no longer preserved.
However, a TriC structure can be obtained by applying nonlinear reconstruction directly to the numerical solution, as in \cite{JCP_2018_Zhu_Shu_MR-WENO}. We leave this extension for the future work.
Furthermore, the constant $c$ in the Def.~\ref{1D def tric} is the length of the reconstruction stencil.
Thus, the linear reconstruction in \cite{MCS_2023_Ahmat_Qiu} is TriC and consistent.

\subsubsection{Ultra-weak discontinuous Galerkin scheme}
For the UWDG method proposed in \cite{MCOM_2008_ultraweak}, the triply conservative structure follows naturally by testing the scheme with the quadratic B-spline $B_{2,j}$ in \eqref{non-uniform-B2-piecewise}, since the test functions are defined independently on each cell and $B_{2,j}$ is itself a piecewise polynomial. Summing the resulting local formulations over all cells then yields the triply conservative form of the UWDG scheme.
The details are given below.

The semidiscrete UWDG scheme for \eqref{1D TC model} seeks $u_h\in \mathbb{P}^{k}\ (k\geq 2)$, such that the following equation holds for any $v_{h}\in\mathbb{P}^{k}$ on each cell $I_{j}$:
\begin{equation}
	\begin{aligned}
		&\int_{I_j} (u_h)_t v_h   dx
		- \int_{I_j} f(u_h) (v_h)_x   dx
		- \int_{I_j} A(u_h) (v_h)_{xxx}   dx
		+ \left( \widehat{f}(u_h^-,u_h^+) v_h^- \right)_{j+\frac{1}{2}}
		- \left( \widehat{f}(u_h^-,u_h^+) v_h^+ \right)_{j-\frac{1}{2}}
		\\
		&\quad
		+ \left( \widehat{A}_h (v_h)_{xx}^- \right)_{j+\frac{1}{2}}
		- \left( \widehat{A}_h (v_h)_{xx}^+ \right)_{j-\frac{1}{2}}
		- \left( \widetilde{(A_h)_x} (v_h)_x^- \right)_{j+\frac{1}{2}}
		+ \left( \widetilde{(A_h)_x} (v_h)_x^+ \right)_{j-\frac{1}{2}}
		\\
		&\quad
		+ \left( \widetilde{(A_h)_{xx}} v_h^- \right)_{j+\frac{1}{2}}
		- \left( \widetilde{(A_h)_{xx}} v_h^+ \right)_{j-\frac{1}{2}}
		= 0 .
	\end{aligned}
\end{equation}
where $\widehat{A}_h,\widetilde{(A_h)_x}, \widetilde{(A_h)_{xx}}$ are numerical fluxes. In particular, $\widehat{A}_{h}$ is Lipschitz continuous and satisfies the constant-preserving property.
Taking the test function $v_{h}=B_{2,j}$ piecewisely, applying the forward Euler discretization in time, taking the summation of all of them, we find that the terms such as $\widetilde{(A_h)_x}, \widetilde{(A_h)_{xx}}$ are canceled. 
The constant $c$ of the UWDG scheme in the Definition~\ref{1D def tric} is $5$.
We can directly get the scheme verified as a consistent TriC scheme. 

\subsubsection{Local discontinuous Galerkin scheme}
Here we discuss the triple conservation of a LDG scheme, namely the LDG-I method which is proposed in \cite{JCP_2004_Levy_Shu_Yan} for the $K(n,n)$ equation. Similar to the UWDG scheme, the test function can be taken piecewise as $B_{2,j}$.
Summing all the cells, then we will get the TriC structure of the $\mathbb{P}^{k}$ ($k\geq2$) LDG-I scheme. The $P^k$ LDG-I scheme $(k\geq2)$ is shown as follows:
\begin{equation}
	\begin{aligned}
		&\int_{I_j} u_t \psi dx
		-\int_{I_j} (v+p)\psi_x dx
		+\left(\widetilde{v}_{j+\frac{1}{2}}+\widehat{p}_{j+\frac{1}{2}}\right)
		\psi_{j+\frac{1}{2}}^{-}
		-\left(\widetilde{v}_{j-\frac{1}{2}}+\widehat{p}_{j-\frac{1}{2}}\right)
		\psi_{j-\frac{1}{2}}^{+}
		=0,
		\\
		&\int_{I_j} p l dx
		+\int_{I_j} q l_x dx
		-\widehat{q}_{j+\frac{1}{2}}l_{j+\frac{1}{2}}^{-}
		+\widehat{q}_{j-\frac{1}{2}}l_{j-\frac{1}{2}}^{+}
		=0,
		\\
		&\int_{I_j} q w dx
		+\int_{I_j} v w_x dx
		-\widehat{v}_{j+\frac{1}{2}}w_{j+\frac{1}{2}}^{-}
		+\widehat{v}_{j-\frac{1}{2}}w_{j-\frac{1}{2}}^{+}
		=0,
		\\
		&\int_{I_j} v z dx
		-\int_{I_j} u^n z dx
		=0.
	\end{aligned}
\end{equation}
The hat quantities are numerical fluxes and specifically, $\hat{v}$ satisfies the constant-preserving property and is locally Lipschitz continuous. Taking $\psi = B_{2,j}, l=\partial_{x}B_{2,j}, w = -\partial_{xx}B_{2,j}$ piecewisely, summing the resulting equations, and applying the forward Euler discretization in time, we obtain the triply conservative form \eqref{1D TriC Def} with $\tilde{A}_{j+\frac{1}{2}}^{s} = \hat{v}_{j+\frac{1}{2}}^{s}$. 
By the last equation of the LDG system, $v=\Pi((u_h^s)^n)$, where $\Pi$ denotes the local $L^2$-projection. Since $\Pi$ preserves 
constants and is Lipschitz continuous of $u_h$ if $u_h$ is bounded \cite{Book_Brenner_Scott}, $\widehat v$ satisfies the consistency in Def.~\ref{1D def consistency}. 
The constant $c$ of the LDG-I scheme in the Def.~\ref{1D def tric} is $5$.
Hence, the LDG-I scheme is triply conservative and consistent.

\section{Triply conservative schemes: two-dimensional case}\label{Sec3}
In this section, we extend the one-dimensional results in
Section~\ref{section 2} to two dimensions.
We consider the following nonlinear degenerate dispersive equation:
\begin{equation}\label{2D TC model}
	u_t
	+ \partial_x f_1(u)
	+ \partial_y f_2(u)
	+ \partial_{xxx} A_1(u)
	+ \partial_{xxy} A_2(u)
	+ \partial_{xyy} A_3(u)
	+ \partial_{yyy} A_4(u)
	= 0,
\end{equation}
where $(x,y)\in\mathbb{R}^2$ and $t>0$.
Here, $f_1(u),f_2(u)$ are the convective fluxes, while
$A_\ell(u)$, $\ell=1,\cdots,4$ are the third-order dispersive terms.
The equation is degenerate dispersive if the derivatives of the dispersive fluxes are allowed to vanish for certain values of $u$.
Following the one-dimensional case, we first introduce the ultra-weak solution of \eqref{2D TC model}.
\begin{definition}
\label{def:2D-ultra-weak}
A function $u(x,y,t) \in L_{loc}^{1}(\mathbb{R}^{2}\times\mathbb{R}_{+})$ is called an ultra-weak solution of \eqref{2D TC model} with initial data $u_0 \in L_{loc}^{1}(\mathbb{R}^{2})$ if it satisfies that: for $f_{i}(u),A_{\ell}(u)\in L_{loc}^{1}(\mathbb{R}^{2}\times\mathbb{R}_{+}),\ i=1,2,\ \ell=1,\ldots,4$, the following equation 
\begin{equation}\label{eq:2D-ultra-weak}
	\begin{aligned}
		&\int_{\mathbb{R}_{+}} \int_{\mathbb{R}^2}
		\Bigl[
		u\phi_t
		+ f_1(u)\phi_x
		+ f_2(u)\phi_y
		+ A_1(u)\phi_{xxx}
		+ A_2(u)\phi_{xxy}
		+ A_3(u)\phi_{xyy}
		+ A_4(u)\phi_{yyy}
		\Bigr] dx dy dt
		\\
		=& -\int_{\mathbb{R}^2}
		u_0(x,y)\phi(x,y,0)dx dy,
	\end{aligned}
\end{equation}
holds for any $\phi\in C_c^\infty(\mathbb{R}^3)$.
\end{definition}
 We next introduce the definition of TriC schemes of the equation \eqref{2D TC model} in two dimensions and establish
a Lax-Wendroff type theorem.

\subsection{Definitions of TriC schemes in two dimensions}
We give the definition of TriC schemes on two-dimensional Cartesian meshes. Let $\mathcal{T}_h=\{I_{ij}\}_{i,j\in\mathbb{Z}}$ be a possibly non-uniform Cartesian mesh of $\mathbb{R}^2$, with the rectangular cells $I_{ij}=I_i^{x}\times I_j^{y}$, where $I_i^{x}=[x_{i-\frac{1}{2}},x_{i+\frac{1}{2}}]$ and
$I_j^{y}=[y_{j-\frac{1}{2}},y_{j+\frac{1}{2}}]$. 
The midpoint of $I_{ij}$ is $(x_{i},y_{j})$ with $x_{i} = (x_{i-\frac{1}{2}}+x_{i+\frac{1}{2}})/2$ and $y_{j} = (y_{j-\frac{1}{2}}+y_{j+\frac{1}{2}})/2$.
We denote the cell widths and areas by $\Delta x_i=x_{i+\frac{1}{2}}-x_{i-\frac{1}{2}},
\ \Delta y_j=y_{j+\frac{1}{2}}-y_{j-\frac{1}{2}}, \ |I_{ij}|=\Delta x_i\Delta y_j$.
The maximal edge length is denoted by $h=\sup_{i,j}\max\{\Delta x_i,\Delta y_j\}$. 
Now we give the TriC definition on two-dimensional Cartesian meshes, which is the tensor product version of \eqref{1D TriC Def}.

\begin{definition}[Triply conservative scheme]
	\label{def:2D-TriC}
	A numerical scheme for \eqref{2D TC model} is TriC if it can be written in the following form:
	\begin{equation} \label{2D-TriC-full-form}
	\begin{aligned}
		&\frac{\widetilde{u}_{ij}^{s+1}-\widetilde{u}_{ij}^{s}}{\Delta t_s}
		+\frac{\widetilde{f}_{1,i+\frac{1}{2},j}^{s}
			-\widetilde{f}_{1,i-\frac{1}{2},j}^{s}}{V_i^x}
		+\frac{\widetilde{f}_{2,i,j+\frac{1}{2}}^{s}
			-\widetilde{f}_{2,i,j-\frac{1}{2}}^{s}}{V_j^y}
		\\
		&+\frac{1}{V_i^x}\Biggl[
		\frac{2}{\Delta x_i+\Delta x_{i+1}}
		\left(
		\frac{\widetilde{A}_{1,i+\frac{3}{2},j}^{s}
			-\widetilde{A}_{1,i+\frac{1}{2},j}^{s}}{\Delta x_{i+1}}
		-
		\frac{\widetilde{A}_{1,i+\frac{1}{2},j}^{s}
			-\widetilde{A}_{1,i-\frac{1}{2},j}^{s}}{\Delta x_i}
		\right)
		\\
		&\qquad\quad-
		\frac{2}{\Delta x_{i-1}+\Delta x_i}
		\left(
		\frac{\widetilde{A}_{1,i+\frac{1}{2},j}^{s}
			-\widetilde{A}_{1,i-\frac{1}{2},j}^{s}}{\Delta x_i}
		-
		\frac{\widetilde{A}_{1,i-\frac{1}{2},j}^{s}
			-\widetilde{A}_{1,i-\frac{3}{2},j}^{s}}{\Delta x_{i-1}}
		\right)\Biggr]
		\\
		&+\frac{1}{V_i^xV_j^y}\Biggl[
		\frac{2}{\Delta x_i+\Delta x_{i+1}}
		\Bigl(
		\widetilde{A}_{2,i+1,j+\frac{1}{2}}^{s}
		-\widetilde{A}_{2,i+1,j-\frac{1}{2}}^{s}
		-\widetilde{A}_{2,i,j+\frac{1}{2}}^{s}
		+\widetilde{A}_{2,i,j-\frac{1}{2}}^{s}
		\Bigr)
		\\
		&\qquad\quad-
		\frac{2}{\Delta x_{i-1}+\Delta x_i}
		\Bigl(
		\widetilde{A}_{2,i,j+\frac{1}{2}}^{s}
		-\widetilde{A}_{2,i,j-\frac{1}{2}}^{s}
		-\widetilde{A}_{2,i-1,j+\frac{1}{2}}^{s}
		+\widetilde{A}_{2,i-1,j-\frac{1}{2}}^{s}
		\Bigr)\Biggr]
		\\
		&+\frac{1}{V_i^xV_j^y}\Biggl[
		\frac{2}{\Delta y_j+\Delta y_{j+1}}
		\Bigl(
		\widetilde{A}_{3,i+\frac{1}{2},j+1}^{s}
		-\widetilde{A}_{3,i-\frac{1}{2},j+1}^{s}
		-\widetilde{A}_{3,i+\frac{1}{2},j}^{s}
		+\widetilde{A}_{3,i-\frac{1}{2},j}^{s}
		\Bigr)
		\\
		&\qquad\quad-
		\frac{2}{\Delta y_{j-1}+\Delta y_j}
		\Bigl(
		\widetilde{A}_{3,i+\frac{1}{2},j}^{s}
		-\widetilde{A}_{3,i-\frac{1}{2},j}^{s}
		-\widetilde{A}_{3,i+\frac{1}{2},j-1}^{s}
		+\widetilde{A}_{3,i-\frac{1}{2},j-1}^{s}
		\Bigr)\Biggr]
		\\
		&+\frac{1}{V_j^y}\Biggl[
		\frac{2}{\Delta y_j+\Delta y_{j+1}}
		\left(
		\frac{\widetilde{A}_{4,i,j+\frac{3}{2}}^{s}
			-\widetilde{A}_{4,i,j+\frac{1}{2}}^{s}}{\Delta y_{j+1}}
		-
		\frac{\widetilde{A}_{4,i,j+\frac{1}{2}}^{s}
			-\widetilde{A}_{4,i,j-\frac{1}{2}}^{s}}{\Delta y_j}
		\right)
		\\
		&\qquad\quad-
		\frac{2}{\Delta y_{j-1}+\Delta y_j}
		\left(
		\frac{\widetilde{A}_{4,i,j+\frac{1}{2}}^{s}
			-\widetilde{A}_{4,i,j-\frac{1}{2}}^{s}}{\Delta y_j}
		-
		\frac{\widetilde{A}_{4,i,j-\frac{1}{2}}^{s}
			-\widetilde{A}_{4,i,j-\frac{3}{2}}^{s}}{\Delta y_{j-1}}
		\right)\Biggr]
		=0,
	\end{aligned}
\end{equation}
where $V_{i}^{x} = \left(\Delta x_{i-1} + \Delta x_{i} + \Delta x_{i+1}\right)/3$ and $V_{j}^{y} = \left(\Delta y_{j-1} + \Delta y_{j} + \Delta y_{j+1}\right)/3$.	As in \cite{C&F_2018_shi_shu,MCOM_2026_Yu_Cheng_Liu_Shu}, the quantity $\widetilde{u}_{ij}^{s}$ is the generalized conserved variable and the quantities $\widetilde{f}_{1,i+\frac{1}{2},j}^{s},\widetilde{f}_{2,i,j+\frac{1}{2}}^{s},\widetilde{A}_{1,i+\frac{1}{2},j}^{s},$ $\dots,\widetilde{A}_{4,i,j+\frac{1}{2}}^{s}$ are generalized numerical fluxes. Each tilde quantity depends locally on the numerical solution $u_{h}^{s}(B_{ij})$ where $B_{ij} = \{\mathbf{x}\in\mathbb{R}^{2}:|\mathbf{x}-\mathbf{x}_{ij}|\leq ch\}$, where $\mathbf{x}_{ij}$ is the center of $I_{ij}$ and $c$ is a positive integer independent of the mesh size.
\end{definition}

We also require additional consistency.

\begin{definition}[Consistency]\label{def:2D-consistency}
A TriC scheme for \eqref{2D TC model} is consistent if:
\begin{enumerate}[label=(\roman*)]
\item If $u_h^s\equiv u$ on $B_{ij}$, then
\begin{equation}\label{eq:2D-consistency-constants}
\begin{aligned}
&\widetilde{u}_{ij}^s=u,\ 
\widetilde{f}_{1,i+\frac{1}{2},j}^s=f_{1}(u),\
\widetilde{f}_{2,i,j+\frac{1}{2}}^s=f_{2}(u),\\ &\widetilde{A}_{\ell_{1},i+\frac{1}{2},j}^s=A_{\ell_{1}}(u),\ \ell_{1}=1,3,\  \widetilde{A}_{\ell_{2},i,j+\frac{1}{2}}^s=A_{\ell_{2}}(u),\ \ell_{2}=2,4.
\end{aligned}
\end{equation}
\item The tilde functions are all Lipschitz
\begin{equation}\label{eq:2D-consistency-Lipschitz}
\begin{aligned}
&\left|\widetilde{u}_{ij}^s - \widetilde{v}_{ij}^s\right| + \left|\widetilde{f}_{1,i+\frac{1}{2},j}^s(u_h)
-\widetilde{f}_{1,i+\frac{1}{2},j}^s(v_h)\right|  
 + \left| \widetilde{f}_{2,i,j+\frac{1}{2}}^s(u_h)
- \widetilde{f}_{2,i,j+\frac{1}{2}}^s(v_h) \right|\\
+ & \left| \widetilde{A}_{1,i+\frac{1}{2},j}^s(u_h)
- \widetilde{A}_{1,i+\frac{1}{2},j}^s(v_h) \right|
+ \left| \widetilde{A}_{2,i,j+\frac{1}{2}}^s(u_h)
 - \widetilde{A}_{2,i,j+\frac{1}{2}}^s(v_h) \right| \\
+ & \left| \widetilde{A}_{3,i+\frac{1}{2},j}^s(u_h)
 - \widetilde{A}_{3,i+\frac{1}{2},j}^s(v_h)\right| 
 + \left| \widetilde{A}_{4,i,j+\frac{1}{2}}^s(u_h)
 - \widetilde{A}_{4,i,j+\frac{1}{2}}^s(v_h) \right|\\
\leq & C \|u_h-v_h\|_{L^\infty(B_{ij})},
\end{aligned}
\end{equation}
for any two numerical solutions $u_h,\ v_h$ and any indices $i,j$.

\item The locally conserved quantity could recover the global conservation:
\begin{equation}\label{eq:2D-consistency-mass}
\sum_{i,j}V_i^xV_j^y\widetilde{u}_{ij}^s
=\int u_h^s(x,y) dx dy,
\quad s\geq0.
\end{equation}
\end{enumerate}
\end{definition}

\subsection{Lax-Wendroff type theorem}
\label{subsec:2D-LW}
In this subsection, we prove the two-dimensional Lax-Wendroff type
theorem. For brevity, we mainly present the argument for the dispersive
terms; the convective terms and locally conserved quantity are handled as before \cite{C&F_2018_shi_shu, MCOM_2026_Yu_Cheng_Liu_Shu}.
We consider Cartesian meshes satisfying
$c_0h\leq\Delta x_i,\Delta y_j\leq h$, with $c_0>0$
independent of the mesh sizes. As in the one-dimensional case,
we require the initial approximation to satisfy
\begin{equation}\label{eq:2D-LW-initial}
	\int_{\mathbb{R}^2}(u_0-u_h^0)\phi(x,y,0) dx dy\longrightarrow0,
	\qquad h,\tau\to0,\qquad
	\forall\phi\in C_c^\infty(\mathbb{R}^3).
\end{equation}
We also assume the following local bounded variation condition:
\begin{equation}\label{eq:2D-LW-oscillation}
    h^2\sum_{i,j}\sup_{(x,y)\in B_{ij}}
	|u_h^s(x,y)-u_h^s(x_i,y_j)|
	\to 0,\quad \text{as}\ h\to0,\ \text{for}\ \forall s\geq 0.
\end{equation}
We now establish a Lax-Wendroff type theorem for \eqref{2D TC model}.
\begin{theorem}[Lax-Wendroff type theorem]
	\label{thm:2D-LW}
Suppose that $\{u_h\}$ is generated by a consistent TriC scheme
\eqref{2D-TriC-full-form} and satisfies the above assumptions.
If $\{u_h\}$ has uniformly bounded total variation and support and converges in $L^{1}$ to a function $u$ as $\tau,h\to0$, then $u$ is an ultra-weak solution to \eqref{2D TC model}.
\end{theorem}
\begin{proof}
Define $\phi_{ij}^s=\phi(x_i,y_j,t_s)$ with
$\phi\in C_0^\infty(\mathbb R^3)$.
As in the one-dimensional case, we omit the superscript $s$ and the
variable $t$ whenever no confusion can arise.
We multiply \eqref{2D-TriC-full-form} by
$\Delta t_sV_i^xV_j^y\phi_{ij}$ and sum over $i,j,s$.
The temporal and convective terms are handled as in
\cite{C&F_2018_shi_shu,MCOM_2026_Yu_Cheng_Liu_Shu}.
We mainly consider the dispersive terms.
Define an auxiliary linear functional:
 	\begin{equation}\label{2D-LW-lambda}
 		\begin{aligned}
 			\lambda_{ij}
 			=\frac{1}{4}\left(
 			\phi_{i-\frac{1}{2},j-\frac{1}{2}}+\phi_{i-\frac{1}{2},j+\frac{1}{2}}
 			+\phi_{i+\frac{1}{2},j-\frac{1}{2}}+\phi_{i+\frac{1}{2},j+\frac{1}{2}}
 			\right)+4\phi_{ij}
 			-\phi_{i-\frac{1}{2},j}-\phi_{i+\frac{1}{2},j}
 			-\phi_{i,j-\frac{1}{2}}-\phi_{i,j+\frac{1}{2}},
 		\end{aligned}
 	\end{equation}
This is the tensor product of \eqref{lambdaj}, and Taylor expansion gives
 	$\lambda_{ij}=\phi_{ij}+\mathcal O_\phi(h^2)$.
 	We first establish the limits with $\lambda_{ij}$ and then estimate the
 	correction corresponding to $\phi_{ij}-\lambda_{ij}$.

	Write $u_{ij}^s=u_h^s(x_i,y_j)$ and denote by
	$\mathcal S_{xxx}$ and $\mathcal S_{xxy}$ the contributions
	of the $xxx$ and $xxy$ terms to the summed identity.
	Repeated summation by parts, as in \eqref{SBP 1D}, gives
	\begin{equation}\label{eq:2D-LW-pure-limit}
		\begin{aligned}
			\mathcal S_{xxx}
			={}&-\sum_s\Delta t_s\sum_{i,j}V_j^y
			\widetilde A_{1,i+\frac{1}{2},j}^s\omega^x_{i+\frac{1}{2}}
			\phi_{xxx}(x_{i+\frac{1}{2}},y_j)
			+\mathcal O(h)\\
			={}&-\sum_s\Delta t_s\sum_{i,j}V_j^y
			A_1(u_{ij}^s)\omega^x_{i+\frac{1}{2}}
			\phi_{xxx}(x_{i+\frac{1}{2}},y_j)+o(1)\\
			={}&-\sum_s\Delta t_s\sum_{i,j}V_j^y
			\frac{\Delta x_i+\Delta x_{i+1}}2 A_1(u_{ij}^s)
			\phi_{xxx}(x_{i+\frac{1}{2}},y_j)+o(1)\\
			={}&-\int_{\mathbb{R}_{+}}\int_{\mathbb{R}^2}
			A_1(u_h)\phi_{xxx} dx dy dt+o(1)\\
			\longrightarrow{}&-\int_{\mathbb{R}_{+}}\int_{\mathbb{R}^2}
			A_1(u)\phi_{xxx} dx dy dt.
		\end{aligned}
	\end{equation}
	where $\omega^x_{i+\frac{1}{2}}
	=(\Delta x_i+\Delta x_{i+1})/2+\eta_{i+1}^x-\eta_i^x$.
	The weights $\omega^x_{i+\frac{1}{2}}$ and $\eta_i^x$ are those in
	\eqref{1Deq:omega}--\eqref{1Deq:eta-def}, with $\Delta x_j$
	replaced by $\Delta x_i$ and $\eta_i^x=\mathcal O(h)$.
	
	The first equality follows from the one-dimensional calculation
	and Taylor expansion. For the second equality, consistency gives
	\[
	|\widetilde A_{1,i+\frac{1}{2},j}^s-A_1(u_{ij}^s)|
	\leq C\sup_{(x,y)\in B_{ij}}|u_h^s(x,y)-u_{ij}^s|,
	\]
	and \eqref{eq:2D-LW-oscillation} makes the resulting error tend
	to zero. For the third equality, summation by parts gives,
	with time indices omitted,
	\begin{equation}\label{eq:2D-LW-eta-BV}
		\begin{aligned}
			&\sum_{i,j}V_j^y(\eta^x_{i+1}-\eta_i^x)
			A_1(u_{ij})\phi_{xxx}(x_{i+\frac{1}{2}},y_j)\\
			&\quad=-\sum_{i,j}V_j^y\eta_i^x
			\left[A_1(u_{ij})\phi_{xxx}(x_{i+\frac{1}{2}},y_j)
			-A_1(u_{i-1,j})\phi_{xxx}(x_{i-\frac{1}{2}},y_j)\right]
			=o(1).
		\end{aligned}
	\end{equation}
	Here we used the Lipschitz continuity of $A_1$, the smoothness
	of $\phi$, and the local BV condition
	\eqref{eq:2D-LW-oscillation}.
	The final limit uses the $L^1$ convergence of $u_h$ and the
	Lipschitz continuity of $A_1$.
	
	For the $xxy$ term, repeated summation by parts gives
	\begin{equation}\label{eq:2D-LW-mixed-SBP}
		\begin{aligned}
			&\sum_{i,j}\Bigg[
			\frac{2}{\Delta x_i+\Delta x_{i+1}}
			\left(
			\widetilde A_{2,i+1,j+\frac12}-\widetilde A_{2,i+1,j-\frac12}
			-\widetilde A_{2,i,j+\frac12}+\widetilde A_{2,i,j-\frac12}
			\right)
			\\
			&
			-\frac{2}{\Delta x_{i-1}+\Delta x_i}
			\left(
			\widetilde A_{2,i,j+\frac12}-\widetilde A_{2,i,j-\frac12}
			-\widetilde A_{2,i-1,j+\frac12}+\widetilde A_{2,i-1,j-\frac12}
			\right)
			\Bigg]\lambda_{ij}
			\\
		=	&-\sum_{i,j}\widetilde A_{2,i,j+\frac12}
			\Bigg[
			\frac{2(\lambda_{i+1,j+1}-\lambda_{i,j+1}
				-\lambda_{i+1,j}+\lambda_{ij})}{\Delta x_i+\Delta x_{i+1}}-\frac{2(\lambda_{i,j+1}-\lambda_{i-1,j+1}
				-\lambda_{ij}+\lambda_{i-1,j})}{\Delta x_{i-1}+\Delta x_i}
			\Bigg].
		\end{aligned}
	\end{equation}
	Substituting \eqref{2D-LW-lambda} and applying Taylor expansion, we obtain
	\begin{equation}\label{eq:2D-LW-mixed-Taylor}
		\begin{aligned}
			&\frac{2(\lambda_{i+1,j+1}-\lambda_{i,j+1}
				-\lambda_{i+1,j}+\lambda_{ij})}{\Delta x_i+\Delta x_{i+1}}-\frac{2(\lambda_{i,j+1}-\lambda_{i-1,j+1}
				-\lambda_{ij}+\lambda_{i-1,j})}{\Delta x_{i-1}+\Delta x_i}\\
			=&\Delta x_i\frac{\Delta y_j+\Delta y_{j+1}}2
			\phi_{xxy}(x_i,y_{j+\frac12})+\mathcal{O}(h^3).
		\end{aligned}
	\end{equation}
	The local remainder contributes $\mathcal{O}(h)$ after summation.
	Consequently, the right-hand side of \eqref{eq:2D-LW-mixed-SBP} equals
	\begin{equation}\label{eq:2D-LW-mixed-limit}
		\begin{aligned}
			&-\sum_{i,j}\Delta x_i\frac{\Delta y_j+\Delta y_{j+1}}2
			A_2(u_{ij})\phi_{xxy}(x_i,y_{j+\frac12})+o(1)=-\int_{\mathbb R^2}A_2(u_h)\phi_{xxy} dx dy+o(1),
		\end{aligned}
	\end{equation}
	by local consistency and \eqref{eq:2D-LW-oscillation}.

	The $xyy$ and $yyy$ terms follow by interchanging the coordinates.
	The convective and time terms are handled as in the previous work \cite{C&F_2018_shi_shu,MCOM_2026_Yu_Cheng_Liu_Shu}.

Finally, consider the correction involving $\phi_{ij}^s-\lambda_{ij}^s$.
By \eqref{2D-TriC-full-form}, the sum of the four dispersive contributions,
after multiplication by
$\Delta t_sV_i^xV_j^y(\phi_{ij}^s-\lambda_{ij}^s)$ and summation over
$i,j$, equals
\begin{equation}\label{eq:2D-LW-correction}
	\begin{aligned}
		&-\sum_{i,j}V_i^xV_j^y
		(\widetilde u_{ij}^{s+1}-\widetilde u_{ij}^s)
		(\phi_{ij}^s-\lambda_{ij}^s) -\Delta t_s\sum_{i,j}V_j^y
		(\widetilde f_{1,i+\frac12,j}^s-\widetilde f_{1,i-\frac12,j}^s)
		(\phi_{ij}^s-\lambda_{ij}^s)
		\\&-\Delta t_s\sum_{i,j}V_i^x
		(\widetilde f_{2,i,j+\frac12}^s-\widetilde f_{2,i,j-\frac12}^s)
		(\phi_{ij}^s-\lambda_{ij}^s).
	\end{aligned}
\end{equation}
Notice that $\phi_{ij}^s-\lambda_{ij}^s=\mathcal{O}(h^2)$ and 
\[
\left|
(\phi_{ij}^s-\lambda_{ij}^s)
-(\phi_{ij}^{s-1}-\lambda_{ij}^{s-1})
\right|
\le C h^2\Delta t_{s-1}.
\]
Summing \eqref{eq:2D-LW-correction} over $s$ and applying
summation by parts in time, the first term is
$\mathcal O(h^2)$, including the initial contribution.
For the remaining two terms, summation by parts in space,
local consistency, and \eqref{eq:2D-LW-oscillation} give
$\mathcal O(h)$.
Therefore the total correction tends to zero, which completes
the proof.
\end{proof}

\subsection{A triply conservative finite volume framework in two dimensions}\label{subsec:tric-fv 2D}
We now extend the finite volume scheme constructed in Section~\ref{subsec:tric-fv 1D} to the two-dimensional equation \eqref{2D TC model} on non-uniform Cartesian meshes. 
The construction uses tensor products of the one-dimensional quadratic B-splines.
As in the one-dimensional case, the derivatives of the B-spline will naturally determine the TriC structure. 
The quadratic B-spline in the $x$ direction is
\begin{equation}\label{eq:2D-FV-B2x}
	B_{2,i}^x(x)=
	\begin{cases}
		\displaystyle
		\frac{(x-x_{i-\frac{3}{2}})^2}
		{\Delta x_{i-1}(\Delta x_{i-1}+\Delta x_i)},
		&x\in I_{i-1}^{x},\\ 
		\displaystyle
		\frac{(x-x_{i-\frac{3}{2}})(x_{i+\frac{1}{2}}-x)}
		{(\Delta x_{i-1}+\Delta x_i)\Delta x_i}
		+\frac{(x_{i+\frac{3}{2}}-x)(x-x_{i-\frac{1}{2}})}
		{(\Delta x_i+\Delta x_{i+1})\Delta x_i},
		&x\in I_i^{x},\\
		\displaystyle
		\frac{(x_{i+\frac{3}{2}}-x)^2}
		{\Delta x_{i+1}(\Delta x_i+\Delta x_{i+1})},
		&x\in I_{i+1}^{x},\\
		0,&\text{otherwise}.
	\end{cases}
\end{equation}
Define $B_{2,j}^y(y)$ by the same expression in the $y$ coordinate,
with index $j$ and cell widths $\Delta y_j$. The corresponding
two-dimensional spline is
\begin{equation}\label{eq:2D-FV-tensor-spline}
	B_{2,ij}(x,y)=B_{2,i}^x(x)B_{2,j}^y(y).
\end{equation}
It is nonnegative and supported on the nine cells
$\bigcup_{k=i-1}^{i+1}\bigcup_{\ell=j-1}^{j+1}I_{k\ell}$.
The partition of unity and the spline masses are
\begin{equation}\label{eq:2D-FV-spline-masses}
	\begin{gathered}
		\sum_{i,j}B_{2,ij}(x,y)=1,\qquad
		\int_{\mathbb{R}^2}B_{2,ij}(x,y) dx dy=V_i^xV_j^y.
	\end{gathered}
\end{equation}
We also introduce the linear B-spline
\begin{equation}\label{eq:2D-FV-B1x}
	B_{1,i+\frac{1}{2}}^x(x)=
	\begin{cases}
		(x-x_{i-\frac{1}{2}})/\Delta x_i,&x\in I_i^{x},\\
		(x_{i+\frac{3}{2}}-x)/\Delta x_{i+1},&x\in I_{i+1}^{x},\\
		0,&\text{otherwise},
	\end{cases}
\end{equation}
and define $B_{1,j+\frac{1}{2}}^y$ analogously. As in the one-dimensional case, we have
\begin{align}
	\partial_{x}(B_{2,i}^x)&=
	\frac{2}{\Delta x_{i-1}+\Delta x_i}B_{1,i-\frac{1}{2}}^x
	-\frac{2}{\Delta x_i+\Delta x_{i+1}}B_{1,i+\frac{1}{2}}^x,
	\label{eq:2D-FV-spline-first}\\
	\partial_{xx}(B_{2,i}^x)&=
	\frac{2}{\Delta x_{i-1}+\Delta x_i}
	\left(\frac{\mathbf1_{I_{i-1}^x}}{\Delta x_{i-1}}
	-\frac{\mathbf1_{I_i^x}}{\Delta x_i}\right)
	-\frac{2}{\Delta x_i+\Delta x_{i+1}}
	\left(\frac{\mathbf1_{I_i^x}}{\Delta x_i}
	-\frac{\mathbf1_{I_{i+1}^x}}{\Delta x_{i+1}}\right),
	\label{eq:2D-FV-spline-second}
\end{align}
where $\mathbf{1}_{I_i^x}$ is the characteristic function of the interval $I_{i}^{x}$.
Define the triple cell average using the normalized tensor-product spline moment:
\begin{equation}\label{eq:2D-FV-triple-average}
	\begin{aligned}
		\bar{\bar{\bar u}}_{ij}(t)
		&\triangleq\frac{1}{V_i^xV_j^y}
		\int_{\mathbb{R}^2}u(x,y,t)B_{2,i}^x(x)B_{2,j}^y(y) dx dy.
	\end{aligned}
\end{equation}
Then, for a sufficiently smooth solution, multiplying \eqref{2D TC model} by $B_{2,ij}$ and formally integrating by parts, we have:
\begin{equation}\label{eq:2D-FV-tested-balance}
	\begin{aligned}
		V_{i}^{x} V_{j}^{y} \frac{d}{dt} \bar{\bar{\bar u}}_{ij}(t)
		={}&\int_{\mathbb{R}^2}
		\bigl[f_1(u)\partial_{x}(B_{2,i}^x)B_{2,j}^y
		+f_2(u) B_{2,i}^x \partial_{y}(B_{2,j}^y)\bigr] dx dy\\
		&+\left\langle A_1(u),\partial_{xxx}(B_{2,i}^x) B_{2,j}^y\right\rangle
		+\left\langle A_2(u), \partial_{xx}(B_{2,i}^x) \partial_{y}(B_{2,j}^y)\right\rangle\\
		&+\left\langle A_3(u), \partial_{x}(B_{2,i}^x) \partial_{yy}(B_{2,j}^y) \right\rangle
		+\left\langle A_4(u),B_{2,i}^x \partial_{yyy}(B_{2,j}^y) \right\rangle.
	\end{aligned}
\end{equation}
The brackets denote inner product in the distributional sense. Compact support of the splines
eliminates the boundary terms. This identity \eqref{eq:2D-FV-tested-balance} is used to motivate the discretization.
Using \eqref{eq:2D-FV-spline-first} and \eqref{eq:2D-FV-spline-second}, we define the following fluxes:
\begin{align}
	F_{1,i+\frac{1}{2},j}
	&=\frac{2}{(\Delta x_i+\Delta x_{i+1})V_j^y}
	\int_{\mathbb{R}^2} f_1(u)B_{1,i+\frac{1}{2}}^x B_{2,j}^y dx dy,
	\label{eq:2D-FV-exact-f1}\\
	F_{2,i,j+\frac{1}{2}}
	&=\frac{2}{V_i^x(\Delta y_j+\Delta y_{j+1})}
	\int_{\mathbb{R}^2} f_2(u)B_{2,i}^x B_{1,j+\frac{1}{2}}^y dx dy,
	\label{eq:2D-FV-exact-f2}\\
	G_{1,i+\frac{1}{2},j}
	&=\frac{1}{V_j^y}\int_{\mathbb{R}}
	A_1(u(x_{i+\frac{1}{2}},y,t))B_{2,j}^y(y) dy,
	\label{eq:2D-FV-exact-A1}\\
	G_{4,i,j+\frac{1}{2}}
	&=\frac{1}{V_i^x}\int_{\mathbb{R}}
	A_4(u(x,y_{j+\frac{1}{2}},t))B_{2,i}^x(x) dx,
	\label{eq:2D-FV-exact-A4}\\
	G_{2,i,j+\frac{1}{2}}
	&=\frac{2}{\Delta x_i(\Delta y_j+\Delta y_{j+1})}
	\int_{I_i^{x}}\int_{I_j^y\cup I_{j+1}^y}
	A_2(u)B_{1,j+\frac{1}{2}}^y(y) dy dx,
	\label{eq:2D-FV-exact-A2}\\
	G_{3,i+\frac{1}{2},j}
	&=\frac{2}{(\Delta x_i+\Delta x_{i+1})\Delta y_j}
	\int_{I_i^{x}\cup I_{i+1}^{x}}\int_{I_j^y}
	A_3(u)B_{1,i+\frac{1}{2}}^x(x) dy dx.
	\label{eq:2D-FV-exact-A3}
\end{align}
Each flux recovers the constant state if the $u(x,y,t)$ is a constant. The equation \eqref{eq:2D-FV-tested-balance} is a semi-discrete scheme. Applying the forward Euler temporal discretization to \eqref{eq:2D-FV-tested-balance}, we have the following fully discrete scheme:
\begin{equation}\label{eq:2D-FV-forward-euler}
	\frac{\bar{\bar{\bar{u}}}_{ij}^{s+1}-\bar{\bar{\bar{u}}}_{ij}^s}{\Delta t_s}
	+\frac{\mathcal{F}_{i+\frac{1}{2},j}^{x,s}-\mathcal{F}_{i-\frac{1}{2},j}^{x,s}}{V_i^x}
	+\frac{\mathcal{F}_{i,j+\frac{1}{2}}^{y,s}-\mathcal{F}_{i,j-\frac{1}{2}}^{y,s}}{V_j^y}=0.
\end{equation}
where $\bar{\bar{\bar{u}}}_{ij}^s$ is the evolved triple cell average at $t=t_{s}$ and the normalized fluxes are defined as: 
\begin{equation}\label{eq:2D-FV-combined-x}
	\begin{aligned}
		\mathcal{F}_{i+\frac{1}{2},j}^{x,s}={}&
		\widetilde{f}_{1,i+\frac{1}{2},j}^s +\frac{2}{\Delta x_i+\Delta x_{i+1}}
		\left[
		\frac{\widetilde{A}_{1,i+\frac{3}{2},j}^s-\widetilde{A}_{1,i+\frac{1}{2},j}^s}{\Delta x_{i+1}}
		-\frac{\widetilde{A}_{1,i+\frac{1}{2},j}^s-\widetilde{A}_{1,i-\frac{1}{2},j}^s}{\Delta x_i}
		\right]\\
		&+\frac{2}{\Delta x_i+\Delta x_{i+1}}
		\left[\frac{\widetilde{A}_{2,i+1,j+\frac{1}{2}}^s-\widetilde{A}_{2,i+1,j-\frac{1}{2}}^s}{V_j^y}
		-\frac{\widetilde{A}_{2,i,j+\frac{1}{2}}^s-\widetilde{A}_{2,i,j-\frac{1}{2}}^s}{V_j^y}
		\right],
	\end{aligned}
\end{equation}
\begin{equation}\label{eq:2D-FV-combined-y}
	\begin{aligned}
		\mathcal{F}_{i,j+\frac{1}{2}}^{y,s}={}&
		\widetilde{f}_{2,i,j+\frac{1}{2}}^s+\frac{2}{\Delta y_j+\Delta y_{j+1}}
		\left[
		\frac{\widetilde{A}_{4,i,j+\frac{3}{2}}^s-\widetilde{A}_{4,i,j+\frac{1}{2}}^s}{\Delta y_{j+1}}
		-\frac{\widetilde{A}_{4,i,j+\frac{1}{2}}^s-\widetilde{A}_{4,i,j-\frac{1}{2}}^s}{\Delta y_j}
		\right]\\
		&+\frac{2}{\Delta y_j+\Delta y_{j+1}}
		\left[
		\frac{\widetilde{A}_{3,i+\frac{1}{2},j+1}^s-\widetilde{A}_{3,i-\frac{1}{2},j+1}^s}{V_i^x}
		-\frac{\widetilde{A}_{3,i+\frac{1}{2},j}^s-\widetilde{A}_{3,i-\frac{1}{2},j}^s}{V_i^x}
		\right].
	\end{aligned}
\end{equation}
The tilde quantities in \eqref{eq:2D-FV-combined-x} and
\eqref{eq:2D-FV-combined-y} are numerical approximations to the
normalized flux integrals in
\eqref{eq:2D-FV-exact-f1}-\eqref{eq:2D-FV-exact-A3}, evaluated at
$t=t_s$:
\[
\begin{aligned}
	\widetilde{f}_{1,i+\frac{1}{2},j}^s&\approx F_{1,i+\frac{1}{2},j}(t_s),
	&\qquad \widetilde{f}_{2,i,j+\frac{1}{2}}^s&\approx F_{2,i,j+\frac{1}{2}}(t_s),\\
	\widetilde{A}_{1,i+\frac{1}{2},j}^s&\approx G_{1,i+\frac{1}{2},j}(t_s),
	&\qquad \widetilde{A}_{4,i,j+\frac{1}{2}}^s&\approx G_{4,i,j+\frac{1}{2}}(t_s),\\
	\widetilde{A}_{2,i,j+\frac{1}{2}}^s&\approx G_{2,i,j+\frac{1}{2}}(t_s),
	&\qquad \widetilde{A}_{3,i+\frac{1}{2},j}^s&\approx G_{3,i+\frac{1}{2},j}(t_s).
\end{aligned}
\]
The fluxes in \eqref{eq:2D-FV-exact-f1}-\eqref{eq:2D-FV-exact-A3} are all normalized integrals, thus they will reduce to the physical fluxes for a constant solution. Numerical quadrature with enough accuracy is used to evaluate these integrals. Furthermore, the reconstruction of the point values is also needed.
The reconstruction and quadrature used to approximate these integrals
are specified below.

\subsubsection{Implementation of numerical quadrature and numerical fluxes}
Let $\{(\xi_a,\omega_a)\}_{a=1}^{N_q}$ be a Gauss--Legendre rule
on $[-\frac{1}{2},\frac{1}{2}]$, normalized by
$\sum_{a=1}^{N_q}\omega_a=1$. With $x_k$ and $y_\ell$ denoting
cell midpoints, set
\begin{equation}\label{eq:2D-FV-quadrature-points}
	x_{k,a}=x_k+\Delta x_k\xi_a,\qquad
	y_{\ell,b}=y_\ell+\Delta y_\ell\xi_b.
\end{equation}
Thus the tensor-product rule on $I_{k\ell}$ has physical weights
$\Delta x_k\Delta y_\ell\omega_a\omega_b$.

For stability, we will introduce the upwinding mechanism in both convective and dispersive fluxes, therefore, we will use limit values from different directions. Denote the limit point value associated with the $x$ direction by $u_{i+\frac{1}{2},j}^{x,-,s}(x,y)$ and
$u_{i+\frac{1}{2},j}^{x,+,s}(x,y)$, and those associated with the $y$ direction by $u_{i,j+\frac{1}{2}}^{y,-,s}(x,y)$ and $u_{i,j+\frac{1}{2}}^{y,+,s}(x,y)$.
The detailed implementation of the evolved data is specified in the next subsection.

For the convective terms, there are two integrals \eqref{eq:2D-FV-exact-f1}-\eqref{eq:2D-FV-exact-f2} needed to be evaluated. Define the Lax-Friedrichs flux
\begin{equation}\label{eq:2D-FV-LF}
	\widehat{f}_{r}(u_{L},u_{R};\alpha)
	=\frac{f_r(u_{L})+f_r(u_{R})}{2}-\frac{\alpha}{2}(u_{R}-u_{L}), 
	\quad r=1,2,\quad \alpha\geq\sup_{z\in\mathcal{I}}|f_r'(z)|,
\end{equation}
and denote
\begin{equation}\label{eq:2D-FV-pointwise-convection}
	\begin{aligned}
		f_{1,i+\frac{1}{2},\ell}^{\mathrm{LF},s}(x,y)
		&=\widehat f_1\bigl(
		u_{i+\frac{1}{2},\ell}^{x,-,s}(x,y),
		u_{i+\frac{1}{2},\ell}^{x,+,s}(x,y);
		\alpha_1\bigr),\\
		f_{2,k,j+\frac{1}{2}}^{\mathrm{LF},s}(x,y)
		&=\widehat f_2\bigl(
		u_{k,j+\frac{1}{2}}^{y,-,s}(x,y),
		u_{k,j+\frac{1}{2}}^{y,+,s}(x,y);
		\alpha_2\bigr).
	\end{aligned}
\end{equation}
Then the generalized numerical convective fluxes are
\begin{equation}\label{eq:2D-FV-quadrature-f1}
	\begin{aligned}
		\widetilde{f}_{1,i+\frac{1}{2},j}^s
		={}&\frac{2}{(\Delta x_i+\Delta x_{i+1})V_j^y}
		\sum_{k=i}^{i+1}\sum_{\ell=j-1}^{j+1}
		\Delta x_k\Delta y_\ell\sum_{a,b=1}^{N_q}\omega_a\omega_b \times B_{1,i+\frac{1}{2}}^x(x_{k,a})B_{2,j}^y(y_{\ell,b})
		f_{1,i+\frac{1}{2},\ell}^{\mathrm{LF},s}(x_{k,a},y_{\ell,b}),
	\end{aligned}
\end{equation}
\begin{equation}\label{eq:2D-FV-quadrature-f2}
	\begin{aligned}
		\widetilde{f}_{2,i,j+\frac{1}{2}}^s
		={}&\frac{2}{V_i^x(\Delta y_j+\Delta y_{j+1})}
		\sum_{k=i-1}^{i+1}\sum_{\ell=j}^{j+1}
		\Delta x_k\Delta y_\ell\sum_{a,b=1}^{N_q}\omega_a\omega_b \times B_{2,i}^x(x_{k,a})B_{1,j+\frac{1}{2}}^y(y_{\ell,b})
		f_{2,k,j+\frac{1}{2}}^{\mathrm{LF},s}(x_{k,a},y_{\ell,b}).
	\end{aligned}
\end{equation}

For the dispersive terms, similar to the one-dimensional case, the so-called ``reverse-monotone'' numerical flux $\hat{A}(\cdot,\cdot)$ is needed. This construction is for stability and also employed in \cite{JCP_2004_Levy_Shu_Yan}. 
As in the one-dimensional case, we also adopt the reverse Godunov flux, defined by
\begin{equation}\label{eq:2D-FV-dispersive-point-flux}
	\widehat A_r^{\mathrm G}(v,w)=
	\begin{cases}
		\displaystyle\max_{z\in[v,w]}A_r(z),&v\leq w,\\[1mm]
		\displaystyle\min_{z\in[w,v]}A_r(z),&v>w,
	\end{cases}
	\qquad r=1,2,3,4.
\end{equation}
Next, we divide the dispersive terms into ``pure-terms'' \eqref{eq:2D-FV-exact-A1}-\eqref{eq:2D-FV-exact-A4} and ``mixed-terms'' \eqref{eq:2D-FV-exact-A2}-\eqref{eq:2D-FV-exact-A3}. 
The pure dispersive traces are defined by
\begin{equation}\label{eq:2D-FV-dispersive-traces}
	\begin{aligned}
		a_{1,i+\frac{1}{2},\ell}^s(y)
		&=\widehat A_1^{\mathrm G}\bigl(
		u_{i+\frac{1}{2},\ell}^{x,-,s}(x_{i+\frac{1}{2}},y),
		u_{i+\frac{1}{2},\ell}^{x,+,s}(x_{i+\frac{1}{2}},y)\bigr),\\
		a_{4,k,j+\frac{1}{2}}^s(x)
		&=\widehat A_4^{\mathrm G}\bigl(
		u_{k,j+\frac{1}{2}}^{y,-,s}(x,y_{j+\frac{1}{2}}),
		u_{k,j+\frac{1}{2}}^{y,+,s}(x,y_{j+\frac{1}{2}})\bigr).
	\end{aligned}
\end{equation}
The corresponding integrals \eqref{eq:2D-FV-exact-A1}-\eqref{eq:2D-FV-exact-A4} are approximated by
\begin{align}
	\widetilde{A}_{1,i+\frac{1}{2},j}^s
	&=\frac{1}{V_j^y}\sum_{\ell=j-1}^{j+1}\Delta y_\ell
	\sum_{b=1}^{N_q}\omega_b B_{2,j}^y(y_{\ell,b})
	\widehat A_1^{\mathrm G}\bigl(
	u_{i+\frac{1}{2},\ell}^{x,-,s}(x_{i+\frac{1}{2}},y_{\ell,b}),
	u_{i+\frac{1}{2},\ell}^{x,+,s}(x_{i+\frac{1}{2}},y_{\ell,b})\bigr),
	\label{eq:2D-FV-quadrature-A1}\\
	\widetilde{A}_{4,i,j+\frac{1}{2}}^s
	&=\frac{1}{V_i^x}\sum_{k=i-1}^{i+1}\Delta x_k
	\sum_{a=1}^{N_q}\omega_a B_{2,i}^x(x_{k,a})
	\widehat{A}_4^{\mathrm G}\bigl(
	u_{k,j+\frac{1}{2}}^{y,-,s}(x_{k,a},y_{j+\frac{1}{2}}),
	u_{k,j+\frac{1}{2}}^{y,+,s}(x_{k,a},y_{j+\frac{1}{2}})\bigr).
	\label{eq:2D-FV-quadrature-A4}
\end{align}
For the mixed dispersive terms, the same dispersive traces are defined by
\begin{equation}\label{eq:2D-FV-pointwise-mixed}
	\begin{aligned}
		a_{2,i,j+\frac{1}{2}}^s(x,y)
		&=\widehat A_2^{\mathrm G}\bigl(
		u_{i,j+\frac{1}{2}}^{y,-,s}(x,y),
		u_{i,j+\frac{1}{2}}^{y,+,s}(x,y)\bigr),\\
		a_{3,i+\frac{1}{2},j}^s(x,y)
		&=\widehat A_3^{\mathrm G}\bigl(
		u_{i+\frac{1}{2},j}^{x,-,s}(x,y),
		u_{i+\frac{1}{2},j}^{x,+,s}(x,y)\bigr).
	\end{aligned}
\end{equation}
The corresponding integrals \eqref{eq:2D-FV-exact-A2}-\eqref{eq:2D-FV-exact-A3} are approximated by
\begin{equation}\label{eq:2D-FV-quadrature-A2}
	\begin{aligned}
		\widetilde{A}_{2,i,j+\frac{1}{2}}^s
		=\frac{2}{\Delta x_i(\Delta y_j+\Delta y_{j+1})}
		\sum_{\ell=j}^{j+1}\Delta x_i\Delta y_\ell
		\sum_{a,b=1}^{N_q}\omega_a\omega_b \times B_{1,j+\frac{1}{2}}^y(y_{\ell,b})
		a_{2,i,j+\frac{1}{2}}^s(x_{i,a},y_{\ell,b}),
	\end{aligned}
\end{equation}
\begin{equation}\label{eq:2D-FV-quadrature-A3}
	\begin{aligned}
		\widetilde{A}_{3,i+\frac{1}{2},j}^s
		=\frac{2}{(\Delta x_i+\Delta x_{i+1})\Delta y_j}
		\sum_{k=i}^{i+1}\Delta x_k\Delta y_j
		\sum_{a,b=1}^{N_q}\omega_a\omega_b \times B_{1,i+\frac{1}{2}}^x(x_{k,a})
		a_{3,i+\frac{1}{2},j}^s(x_{k,a},y_{j,b}).
	\end{aligned}
\end{equation}
Notice that the point values used in the numerical fluxes above are not directly available from the evolved triple cell average and must therefore be supplied by reconstruction.
Actually, the following point values are needed:
\begin{equation}\label{eq:2D-FV-required-pointwise-states}
	\begin{aligned}
		&u_{i+\frac{1}{2},\ell}^{x,\pm,s}(x_{k,a},y_{\ell,b}),
		&&k=i,i+1,\  \ell=j-1,j,j+1,\quad u_{k,j+\frac{1}{2}}^{y,\pm,s}(x_{k,a},y_{\ell,b}),
		\quad k=i-1,i,i+1,\ \ell=j,j+1,\\
		&u_{i+\frac{1}{2},\ell}^{x,\pm,s}(x_{i+\frac{1}{2}},y_{\ell,b}),
		&&\ell=j-1,j,j+1,\quad u_{k,j+\frac{1}{2}}^{y,\pm,s}(x_{k,a},y_{j+\frac{1}{2}}),\quad k=i-1,i,i+1,
	\end{aligned}
\end{equation}
where $a,b=1,\ldots,N_q$, and the superscript $\pm$ denotes the two limit values from different directions.

\subsubsection{Triple cell averages and linear polynomial reconstruction}
Now we introduce how to get high-order point values \eqref{eq:2D-FV-required-pointwise-states} from two-dimensional triple cell averages  \eqref{eq:2D-FV-triple-average}. 
Consider the cell $I_{ij}$ and its corresponding $5\times5$ stencil \begin{equation}\label{eq:2D-FV-reconstruction-stencil}
	\mathcal{S}_{ij}
	=\{I_{i+r,j+t}:r,t=-2,-1,0,1,2\}.
\end{equation}
The reconstruction is to find a polynomial $u_{ij}^s\in\mathbb{Q}^4=\mathbb{P}^4 \otimes \mathbb{P}^4$ satisfying the following conditions:
\begin{equation}\label{eq:2D-FV-moment-matching}
	\frac{1}{V_{i+r}^xV_{j+\ell }^y}
	\int_{\mathbb{R}^2}
	u_{ij}^s(x,y)
	B_{2,i+r}^x(x)B_{2,j+\ell}^y(y) dx dy
	=\bar{\bar{\bar{u}}}_{i+r,j+\ell}^s,
	\qquad r,\ell=-2,\ldots,2.
\end{equation}
The resulting reconstruction is linear in the prescribed triple cell averages.
All point values used in the numerical fluxes are obtained by evaluating the corresponding cell polynomials, extended to adjacent cells when necessary.

\subsection{Examples of triply conservative schemes}
In addition to the TriC-FV scheme, as in the one-dimensional case, we present several existing high-order TriC schemes in the literature \cite{MCS_2023_Ahmat_Qiu, MCOM_2008_ultraweak, JCP_2004_Levy_Shu_Yan}. These works mainly focus on one-dimensional problems and present the corresponding schemes in one dimension. Although the extensions to two-dimensional Cartesian meshes are straightforward via tensor products, we restrict our discussion to schemes having subsequent developments, particularly those for which explicit two-dimensional formulations are available.

\subsubsection{TriC finite difference scheme}
In \cite{MCS_2023_Ahmat_Qiu}, the authors present the detailed construction in one dimension and remark that the method extends to two dimensions in a dimension-by-dimension manner. The paper also reports two-dimensional numerical results.
However, the corresponding two-dimensional formula is not written out explicitly.
As the example is not the main focus of our paper, here we do not examine its two-dimensional TriC property and consistency, and focus instead on other schemes with explicit two-dimensional formulas.

\subsubsection{TriC ultra-weak discontinuous Galerkin methods}
In \cite{SISC_2022_Huang_Cheng}, the authors developed the UWDG method and proposed a new formulation to solve the Zakharov-Kuznetsov equation:
\begin{equation}
	u_{t} + f(u)_{x} + u_{xxx} + u_{xyy} = 0.
\end{equation}
Thus $f_1=f$, $A_1(u)=A_3(u)=u$, and $f_2=A_2=A_4=0$.
On a uniform Cartesian mesh $\{I_{ij}\}$ of $\Omega$, let $V_h^k=\{v\in L^2(\Omega): v|_{I_{ij}}\in\mathbb{Q}^k(I_{ij})\},\ k\geq2$,
where $\mathbb{Q}^{k}$ is the tensor-product space of the one-dimensional space $\mathbb{P}^{k}$.
Denote $(a,b)_{ij}=\int_{I_{ij}}ab dx dy$ and
$(a,b)=\sum_{i,j}(a,b)_{ij}$.
For a piecewise smooth function $w$, denote its traces by $w_{i+\frac{1}{2},y}^{\pm}=w(x_{i+\frac{1}{2}}^{\pm},y)$, $w_{x,j+\frac{1}{2}}^{\pm}=w(x,y_{j+\frac{1}{2}}^{\pm})$, $w_{i+\frac{1}{2},j+\frac{1}{2}}^{\sigma_{1},\sigma_{2}}=w(x_{i+\frac{1}{2}}^{\sigma_{1}},y_{j+\frac{1}{2}}^{\sigma_{2}}),\ \sigma_{1},\sigma_{2}\in\{-,+\}$.
The semi-discrete UWDG method is to find a numerical solution $u_h \in V_h^k$
such that, for every $v\in V_h^k$ and every cell $I_{ij}$, the following equation holds:
\begin{equation}\label{eq:2D-UWDG-full-form}
	\begin{aligned}
		&(\partial_tu_h,v)_{ij}\\
		=&(f(u_h),v_x)_{ij}
		+(u_h,v_{xxx}+v_{xyy})_{ij}
-\int_{I_j^y}
		\Bigl[
		(\widehat f+\check u_{xx})_{i+\frac{1}{2},y}
		v_{i+\frac{1}{2},y}^{-}
		-
		(\widehat f+\check u_{xx})_{i-\frac{1}{2},y}
		v_{i-\frac{1}{2},y}^{+}
		\Bigr] dy
		\\
		+&\int_{I_j^y}
		\Bigl[
		(\widetilde{u}_x)_{i+\frac{1}{2},y}
		(v_x)_{i+\frac{1}{2},y}^{-}
		-
		(\widetilde{u}_x)_{i-\frac{1}{2},y}
		(v_x)_{i-\frac{1}{2},y}^{+}
		\Bigr] dy
		-\int_{I_j^y}
		\Bigl[
		\widehat u_{i+\frac{1}{2},y}
		(v_{xx})_{i+\frac{1}{2},y}^{-}
		-
		\widehat u_{i-\frac{1}{2},y}
		(v_{xx})_{i-\frac{1}{2},y}^{+}
		\Bigr] dy
		\\
		-&\int_{I_j^y}
		\Bigl[
		\widetilde{u}_{i+\frac{1}{2},y}
		(v_{yy})_{i+\frac{1}{2},y}^{-}
		-
		\widetilde{u}_{i-\frac{1}{2},y}
		(v_{yy})_{i-\frac{1}{2},y}^{+}
		\Bigr] dy
		+\int_{I_i^{x}}
		\Bigl[
		(\check u_y)_{x,j+\frac{1}{2}}
		(v_x)_{x,j+\frac{1}{2}}^{-}
		-
		(\check u_y)_{x,j-\frac{1}{2}}
		(v_x)_{x,j-\frac{1}{2}}^{+}
		\Bigr] dx
		\\
		-&\int_{I_i^{x}}
		\Bigl[
		\widehat u_{x,j+\frac{1}{2}}
		(v_{xy})_{x,j+\frac{1}{2}}^{-}
		-
		\widehat u_{x,j-\frac{1}{2}}
		(v_{xy})_{x,j-\frac{1}{2}}^{+}
		\Bigr] dx
		-(\check u_y)_{i+\frac{1}{2},j+\frac{1}{2}}
		v_{i+\frac{1}{2},j+\frac{1}{2}}^{-,-}
		+(\check u_y)_{i-\frac{1}{2},j+\frac{1}{2}}
		v_{i-\frac{1}{2},j+\frac{1}{2}}^{+,-}
		\\
		+&(\check u_y)_{i+\frac{1}{2},j-\frac{1}{2}}
		v_{i+\frac{1}{2},j-\frac{1}{2}}^{-,+}
		-(\check u_y)_{i-\frac{1}{2},j-\frac{1}{2}}
		v_{i-\frac{1}{2},j-\frac{1}{2}}^{+,+}
		+\widetilde{u}_{i+\frac{1}{2},j+\frac{1}{2}}
		(v_y)_{i+\frac{1}{2},j+\frac{1}{2}}^{-,-}
		-\widetilde{u}_{i-\frac{1}{2},j+\frac{1}{2}}
		(v_y)_{i-\frac{1}{2},j+\frac{1}{2}}^{+,-}
		\\
		-&\widetilde{u}_{i+\frac{1}{2},j-\frac{1}{2}}
		(v_y)_{i+\frac{1}{2},j-\frac{1}{2}}^{-,+}
		+\widetilde{u}_{i-\frac{1}{2},j-\frac{1}{2}}
		(v_y)_{i-\frac{1}{2},j-\frac{1}{2}}^{+,+}.
	\end{aligned}
\end{equation}
The fluxes are chosen as
\begin{equation}\label{eq:2D-UWDG-face-traces}
	\begin{aligned}
		&(\check u_{xx})_{i+\frac{1}{2},y}
		=(\partial_{xx} u_{h})_{i+\frac{1}{2},y}^{+},\quad
		(\widetilde{u}_x)_{i+\frac{1}{2},y}
		=(\partial_{x}u_{h})_{i+\frac{1}{2},y}^{+},\quad
		\widehat u_{i+\frac{1}{2},y}
		=(u_h)_{i+\frac{1}{2},y}^{-},
		\\
		&\widetilde{u}_{i+\frac{1}{2},y}
		=(u_h)_{i+\frac{1}{2},y}^{+},\quad 
		(\check u_y)_{x,j+\frac{1}{2}}
		=(\partial_{y}u_{h})_{x,j+\frac{1}{2}}^{+},\quad
		\widehat u_{x,j+\frac{1}{2}}
		=(u_h)_{x,j+\frac{1}{2}}^{-}.
	\end{aligned}
\end{equation}
The convective flux is the Lax-Friedrichs flux \eqref{eq: 1D LF flux}.
The vertex fluxes are
\begin{equation}\label{eq:2D-UWDG-corner-traces}
	(\check u_y)_{i+\frac{1}{2},j+\frac{1}{2}}
	=(u_{h,y})_{i+\frac{1}{2},j+\frac{1}{2}}^{-,+},
	\qquad
	\widetilde{u}_{i+\frac{1}{2},j+\frac{1}{2}}
	=(u_h)_{i+\frac{1}{2},j+\frac{1}{2}}^{+,-}.
\end{equation}

Now we prove that the scheme \eqref{eq:2D-UWDG-full-form} equipped with the forward Euler method in time can be rewritten in the TriC form \eqref{2D-TriC-full-form}.
Take $\psi_{ij}=B_{2,i}^xB_{2,j}^y\in V_h^k,\ \widetilde{u}_{ij}
	=\frac{(u_h,\psi_{ij})}{V_i^xV_j^y}$. 
The test function is admissible because $k\geq2$. Since the spline function
is $C^1$ and piecewise quadratic, the convective trace terms, derivative-trace terms, all mixed trace terms, and both vertex terms vanish by summing \eqref{eq:2D-UWDG-full-form} over all cells. The remaining equation of the scheme is
\begin{equation}\label{eq:2D-UWDG-spline-balance}
	\begin{aligned}
		\frac{d}{dt}(u_h,\psi_{ij})
		=&\int_\Omega f(u_h)(B_{2,i}^x)'B_{2,j}^y +  u_h(B_{2,i}^x)'(B_{2,j}^y)'' dx dy+\sum_p[(B_{2,i}^x)'']_{x_{p+\frac{1}{2}}}
		\int \widehat u_{p+\frac{1}{2},y}B_{2,j}^y(y) dy,
	\end{aligned}
\end{equation}
where $[\phi]_{x_{j+\frac{1}{2}}} = \phi(x_{j+\frac{1}{2}}^{+}) - \phi(x_{j+\frac{1}{2}}^{-})$.
Define the generalized fluxes by
\begin{subequations}\label{eq:2D-examples-DG-fluxes}
	\begin{align}
		\widetilde{f}_{1,i+\frac{1}{2},j}
		&=\frac{2}{(\Delta x_i+\Delta x_{i+1})V_j^y}
		\int_\Omega f(u_h)B_{1,i+\frac{1}{2}}^xB_{2,j}^y dx dy,
		\label{eq:2D-UWDG-generalized-f1}\\
		\widetilde{A}_{1,i+\frac{1}{2},j}
		&=\frac{1}{V_j^y}\int_{I^{y}_{j-1}\cup I^{y}_{j} \cup I_{j+1}^y}
		\widehat u_{i+\frac{1}{2},y}B_{2,j}^y(y) dy,
		\label{eq:2D-UWDG-generalized-A1}\\
		\widetilde{A}_{3,i+\frac{1}{2},j}
		&=\frac{2}{(\Delta x_i+\Delta x_{i+1})\Delta y_j}
		\int_{I_i^{x}\cup I_{i+1}^{x}}\int_{I_j^y}
		u_h(x,y)B_{1,i+\frac{1}{2}}^x(x) dy dx.
		\label{eq:2D-UWDG-generalized-A3}
	\end{align}
\end{subequations}
Set $\widetilde{f}_2=\widetilde{A}_2=\widetilde{A}_4=0$ where we omit the subscript for simplicity.
Substitution of the spline derivative identities \eqref{eq:2D-FV-spline-first}-\eqref{eq:2D-FV-spline-second}
into \eqref{eq:2D-UWDG-spline-balance} and applying forward Euler in time gives the TriC form \eqref{2D-TriC-full-form}.

Next we consider the consistency \eqref{eq:2D-consistency-constants}-\eqref{eq:2D-consistency-mass} of the UWDG scheme.
It is clear that it holds the constant-preserving property and the positive integer $c$ is taken as $5$. Furthermore, the generalized fluxes and conserved quantity are Lipschitz continuous.
Hence Definition~\ref{def:2D-consistency} holds.

\subsubsection{Local discontinuous Galerkin schemes}
\label{subsubsec:2D-examples-LDGI}
We consider the two Cartesian LDG schemes proposed in Xu and Shu~\cite{CMAME_2007_Xu_Shu}. 
We use the same meshes, DG spaces and notations as before. In \cite{CMAME_2007_Xu_Shu}, they proposed 2D LDG schemes for solving the 2D KdV-type equation and the Zakharov-Kuznetsov equation.
\paragraph{The two-dimensional KdV-type equation.}
For the 2D KdV-type equation $u_t+f(u)_x+g(u)_y+u_{xxx}+u_{yyy}=0$,
the first-order system is
\begin{equation}\label{eq:2D-LDG2007-KdV-system}
	\begin{aligned}
		&u_t+f(u)_x+g(u)_y+p_x+q_y=0,\qquad p=v_x,\ v=u_x,\ q=w_y,\ w=u_y.
	\end{aligned}
\end{equation}
The LDG scheme seeks $u_h,p_h,v_h,q_h,w_h\in V_h^k (k\geq2)$ such that, for every
$\phi,\eta,\zeta,\xi,\chi\in V_h^k$,
\begin{equation}\label{eq:2D-LDG2007-KdV-weak}
	\begin{aligned}
		((u_h)_t,\phi)
		={}&(f(u_h)+p_h,\phi_x)+(g(u_h)+q_h,\phi_y)+\sum_{e\in\mathcal{E}_x}\int_e(\widehat f+\widehat p)[\phi]_x d\sigma 
		+\sum_{e\in\mathcal{E}_y}\int_e(\widehat g+\widehat q)[\phi]_y d\sigma,\\
		(p_h,\eta)
		={}&-(v_h,\eta_x)
		-\sum_{e\in\mathcal{E}_x}\int_e\widehat v[\eta]_x d\sigma,\ (v_h,\zeta)
		={}-(u_h,\zeta_x)
		-\sum_{e\in\mathcal{E}_x}\int_e\widehat u^{ x}[\zeta]_x d\sigma,\\
		(q_h,\xi)
		={}&-(w_h,\xi_y)
		-\sum_{e\in\mathcal{E}_y}\int_e\widehat w[\xi]_y d\sigma,\ (w_h,\chi)
		={}-(u_h,\chi_y)
		-\sum_{e\in\mathcal{E}_y}\int_e\widehat u^{ y}[\chi]_y d\sigma,
	\end{aligned}
\end{equation}
where $\mathcal E_x$ and $\mathcal E_y$ denote the sets
of vertical and horizontal interfaces, respectively,
each counted once, and $[w]_d=w^{+,d}-w^{-,d}$, $d=x,y$. The fluxes are defined as:
\begin{equation}\label{eq:2D-LDG2007-KdV-traces}
	\begin{aligned}
		&\widehat p=p_h^{+,x},\quad
		\widehat v=v_h^{+,x},\quad \widehat u^{ x}=u_h^{-,x},
		\quad &&e\in\mathcal{E}_x,\\
		&\widehat q=q_h^{+,y},\quad
		\widehat w=w_h^{+,y},\quad \widehat u^{ y}=u_h^{-,y},
		&&e\in\mathcal{E}_y.
	\end{aligned}
\end{equation}
Here $\widehat f=\widehat f(u_h^{-,x},u_h^{+,x})$ and
$\widehat g=\widehat g(u_h^{-,y},u_h^{+,y})$ are consistent monotone
fluxes.
As in the UWDG case, 
set $\psi=\psi_{ij}=B_{2,i}^xB_{2,j}^y$ and choose
$\phi=\psi$, $\eta=\psi_x$, $\zeta=-\psi_{xx}$,
$\xi=\psi_y$, and $\chi=-\psi_{yy}$ in \eqref{eq:2D-LDG2007-KdV-weak} and utilize the forward Euler method to discretize the temporal derivative, we can get the TriC form. Furthermore, as the fluxes are chosen as the classical fluxes, the scheme is consistent in the sense of \eqref{eq:2D-consistency-constants}-\eqref{eq:2D-consistency-mass}.

\paragraph{The Zakharov--Kuznetsov equation.}
For $f(u)=3u^2$, the equation
$u_t+f(u)_x+u_{xxx}+u_{xyy}=0$ is written as
\begin{equation}\label{eq:2D-LDG2007-ZK-system}
	u_t+f(u)_x+p_x+r_y=0,\qquad
	p=q_x,\quad q=u_x,\quad r=s_x,\quad s=u_y.
\end{equation}
For this equation, we have $f_{1}(u)=f(u),f_{2}(u)=A_{2}(u)=A_{4}(u)=0,\ A_{1}(u)=A_{3}(u)=u $.  The LDG scheme seeks $u_h,p_h,q_h,r_h,s_h\in V_h^k$ such that, for every
$\phi,\eta,\zeta,\xi,\chi\in V_h^k$,
\begin{equation}\label{eq:2D-LDG2007-ZK-weak}
	\begin{aligned}
		((u_h)_t,\phi)
		={}&(f(u_h)+p_h,\phi_x)+(r_h,\phi_y)+\sum_{e\in\mathcal{E}_x}\int_e(\widehat f+\widehat p)[\phi]_x ds
		+\sum_{e\in\mathcal{E}_y}\int_e\widehat r[\phi]_y d\sigma,\\
		(p_h,\eta)
		={}&-(q_h,\eta_x)
		-\sum_{e\in\mathcal{E}_x}\int_e\widehat q[\eta]_x d\sigma,\	(q_h,\zeta)
		={}-(u_h,\zeta_x)
		-\sum_{e\in\mathcal{E}_x}\int_e\widehat u^{ x}[\zeta]_x d\sigma,\\
		(r_h,\xi)
		={}&-(s_h,\xi_x)
		-\sum_{e\in\mathcal{E}_x}\int_e\widehat s[\xi]_x d\sigma,\ (s_h,\chi)
		={}-(u_h,\chi_y)
		-\sum_{e\in\mathcal{E}_y}\int_e \widehat{u}^{y}[\chi]_y d\sigma.
	\end{aligned}
\end{equation}
The dispersive fluxes are
\begin{equation}\label{eq:2D-LDG2007-ZK-traces}
	\begin{aligned}
		&\widehat p=p_h^{+,x},\quad \widehat q=q_h^{+,x},\quad
		\widehat s=s_h^{+,x},\quad \widehat u^{ x}=u_h^{-,x},
		&&e\in\mathcal{E}_x,\\
		&\widehat{u}^{y}=u_h^{-,y},\qquad
		\widehat r=r_h^{+,y}-\alpha_0[u_h]_y,\qquad \alpha_0>0,
		&&e\in\mathcal{E}_y.
	\end{aligned}
\end{equation}
The convective flux $\widehat f(u_h^{-,x},u_h^{+,x})$ again satisfies
the classical consistency and monotonicity.

Set $\psi=\psi_{ij}$ and choose
$\phi=\psi$, $\eta=\psi_x$, $\zeta=\psi_{xx}$,
$\xi=\psi_y$, and $\chi=\psi_{xy}$ in
\eqref{eq:2D-LDG2007-ZK-weak}.
For simplicity, we also omit the subscript in the following.
Applying forward Euler in time gives the TriC form with
the generalized fluxes in \eqref{eq:2D-examples-DG-fluxes},
where $\widehat u$ in \eqref{eq:2D-UWDG-generalized-A1}
is $\widehat u^{ x}$, and
$\widetilde f_2=\widetilde A_2=\widetilde A_4=0$.

\section{Numerical experiments}\label{Sec4}
In this section, we present the numerical results of the TriC-FV scheme both in one dimension and two dimensions. The results include the accuracy tests, comparisons between the TriC scheme and a non-TriC scheme, dynamics of compactons. 
In presenting the numerical results, if an exact solution is available, we refer to it as ``Exact''; otherwise, we use the result computed by the LDG scheme as the reference solution and refer to it as ``Ref'' or ``LDG''.
When comparing the differences between the TriC scheme and a non-TriC scheme, we refer to them as ``TriC-FV'' and ``Non-TriC'', respectively.
In implementation, the SSPRK3 \eqref{SSPRK3} \cite{SIREV_2001_Gottlieb_Shu_Tadmor} temporal discretization is utilized.
 Let $\mathcal{L}_h(\mathbf U)$ denote the spatial discretization operator.
The three stages are
\begin{equation}\label{SSPRK3}
	\begin{aligned}
		\mathbf U^{(1)}&=\mathbf U^s+\Delta t_s\mathcal{L}_h(\mathbf U^s),\\
		\mathbf U^{(2)}&=\frac34\mathbf U^s
		+\frac{1}{4}\left[\mathbf U^{(1)}
		+\Delta t_s\mathcal{L}_h(\mathbf U^{(1)})\right],\\
		\mathbf U^{s+1}&=\frac{1}{3}\mathbf U^s
		+\frac23\left[\mathbf U^{(2)}
		+\Delta t_s\mathcal{L}_h(\mathbf U^{(2)})\right].
	\end{aligned}
\end{equation}
In one dimension, the time step is taken as:
\begin{equation}\label{eq:num-1d-common-CFL}
	\Delta t_n\leq
	\frac{\nu}{a_n\max\{h_{\min}^{-1},h_{\min}^{-5/3}\}
		+b_nh_{\min}^{-3}},
	\qquad h_{\min}=\min_i\Delta x_i.
\end{equation}
Here $a_n=\max |f'(u_h)|$, $b_n = \max|A'(u_h)|$.
In two dimensions, the same rule becomes
\begin{equation}\label{eq:num-unified-MCS-step}
	\Delta t_n\leq
	\frac{\nu}{
		\begin{aligned}
			&a_{1,n}\max\{h_x^{-1},h_x^{-5/3}\}
			+a_{2,n}\max\{h_y^{-1},h_y^{-5/3}\}\\[-1mm]
			&\quad+b_{1,n}h_x^{-3}+b_{2,n}h_x^{-2}h_y^{-1}
			+b_{3,n}h_x^{-1}h_y^{-2}+b_{4,n}h_y^{-3}
	\end{aligned}},
	\qquad
	h_x=\min_i\Delta x_i,\quad h_y=\min_j\Delta y_j,
\end{equation}
where $a_{i,n}$ and $b_{i,n}$ are the corresponding maxima of convection and dispersion coefficients.
Here we take $\nu = 0.18$ in both \eqref{eq:num-1d-common-CFL} and \eqref{eq:num-unified-MCS-step} for all experiments.

\subsection{One-dimensional numerical results}\label{subsec:num-1d}
\subsubsection{Accuracy for smooth solutions}\label{subsubsec:num-1d-accuracy}
In this part, we test the fifth-order accuracy of the TriC-FV scheme for both linear and nonlinear problems on both uniform and non-uniform meshes.
All tests are taken as periodic on $\Omega=[0,2\pi]$. The mesh is set as
\begin{equation}\label{eq:num-mesh-map}
	x_{k+\frac{1}{2}}=x_L+L\left(s_k+
	\frac{\varepsilon}{2\pi}\sin(2\pi s_k)\right),\qquad
	s_k=\frac{k}{N},\quad k=0,\ldots,N,
\end{equation}
where $L=x_R-x_L$ and $\varepsilon=0$ for the uniform mesh or $\varepsilon=0.2$ for the non-uniform mesh. 

\paragraph{Example 1: linear accuracy test}
The linear problem is
\begin{equation}\label{eq:num-1d-linear}
	u_t+u_x+u_{xxx}=0,\qquad
	u(x,t)=1+0.2\sin x+0.1\cos(2x+6t),
\end{equation}
We take the final time $T = 0.05$.
Table~\ref{tab:num-1d-linear} shows the optimal fifth-order accuracy.

\paragraph{Example 2: nonlinear accuracy test}
For the nonlinear problem, the equation is 
\begin{equation}\label{eq:num-1d-nonlinear}
	\begin{aligned}
		u_t+(u^2)_x+(u^2)_{xxx}&=s(x,t),\qquad
		u(x,t)=1+0.2\sin(x-t),
	\end{aligned}
\end{equation}
with $s(x,t)=-0.2\cos(x-t)-0.12\sin(2x-2t)$. 
Table~\ref{tab:num-1d-nonlinear} presents the errors in different norms at the final time $T=0.1$. 
The results also verify the fifth-order accuracy of the TriC-FV scheme. 
\begin{table}[!htbp]
	\centering\small
	\setlength{\tabcolsep}{5pt}
	\caption{Smooth linear accuracy test \eqref{eq:num-1d-linear}.}\label{tab:num-1d-linear}
	\begin{tabular}{|cr|rr|rr|rr|}
		\hline
		$\varepsilon$ & $N$ & $L^1$-error & order & $L^2$-error & order & $L^\infty$-error & order \\ 
		\hline
0.0 & 20 & 8.73e-04 & -- & 4.01e-04 & -- & 3.47e-04 & -- \\
0.0 & 40 & 2.78e-05 & 4.97 & 1.28e-05 & 4.97 & 1.08e-05 & 5.01 \\
0.0 & 80 & 8.63e-07 & 5.01 & 3.97e-07 & 5.01 & 3.31e-07 & 5.03 \\
\hline
0.2 & 20 & 1.16e-03 & -- & 6.04e-04 & -- & 6.90e-04 & -- \\
0.2 & 40 & 4.05e-05 & 4.85 & 2.16e-05 & 4.81 & 2.60e-05 & 4.73 \\
0.2 & 80 & 1.28e-06 & 4.98 & 6.86e-07 & 4.98 & 8.08e-07 & 5.01 \\
		\hline
	\end{tabular}
	
	\caption{Smooth nonlinear accuracy test \eqref{eq:num-1d-nonlinear}.\label{tab:num-1d-nonlinear}}
	\begin{tabular}{|cr|rr|rr|rr|}
		\hline
		$\varepsilon$ & $N$ & $L^1$-error & order & $L^2$-error & order & $L^\infty$-error & order \\
		\hline
		0.0 & 20 & 5.21e-05 & -- & 2.44e-05 & -- & 2.17e-05 & -- \\
		0.0 & 40 & 1.61e-06 & 5.01 & 7.59e-07 & 5.01 & 6.81e-07 & 4.99 \\
		0.0 & 80 & 5.00e-08 & 5.01 & 2.35e-08 & 5.01 & 2.12e-08 & 5.00 \\
		\hline
		0.2 & 20 & 8.26e-05 & -- & 4.50e-05 & -- & 5.52e-05 & -- \\
		0.2 & 40 & 2.83e-06 & 4.87 & 1.53e-06 & 4.88 & 1.89e-06 & 4.87 \\
		0.2 & 80 & 8.97e-08 & 4.98 & 4.85e-08 & 4.98 & 5.99e-08 & 4.98 \\
		\hline
	\end{tabular}
\end{table}


\subsubsection{Comparison between TriC and non-TriC schemes}\label{compare_1D}
In this part, we compare the performance between the TriC-FV scheme and a non-TriC FD scheme on degenerate convection-dispersion equations, to show the validity of the TriC theory. Furthermore, we also utilize the numerical solution obtained by the $\mathbb{P}^{2}$ LDG-I scheme proposed in \cite{JCP_2004_Levy_Shu_Yan} as the reference solution.
The dispersive part of the non-TriC FD scheme is defined as:
\begin{equation}\label{eq:num-front-chain}
	\frac{v_j^{s+1}-v_j^s}{\Delta t_s}
	+D_+^2\left(3(v^{s+1})^2D_-v^{s+1}\right)_j=0,
	\qquad
	D_\pm v_j=\pm\frac{v_{j\pm1}-v_j}{h}.
\end{equation}
It is consistent for smooth solutions, both conservative and doubly conservative, but lacks the TriC structure. In the following, we use ``Non-TriC'' to denote the specific scheme \eqref{eq:num-front-chain}. We will show that in the cases presented in this part, the scheme \eqref{eq:num-front-chain} may fail to converge to a correct solution.

\begin{figure}[htbp]
	\centering
	\subfigure[Travelling front at $T=1.5$.]{%
		\includegraphics[width=.48\textwidth]{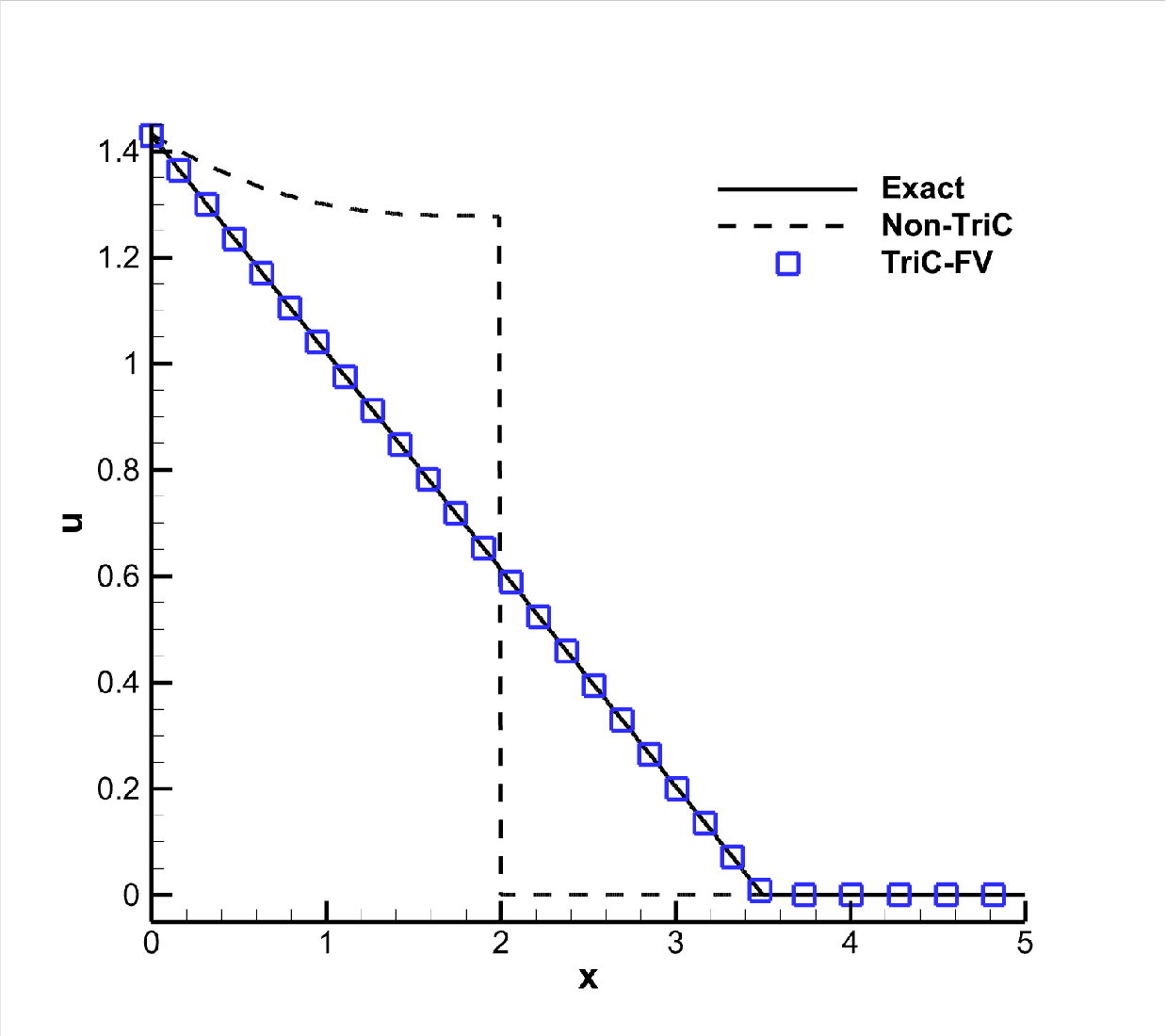}%
		\label{fig:num-front}%
	}
	\quad
	\subfigure[Initial discontinuous problem at $T=0.1$.]{%
		\includegraphics[width=.48\textwidth]{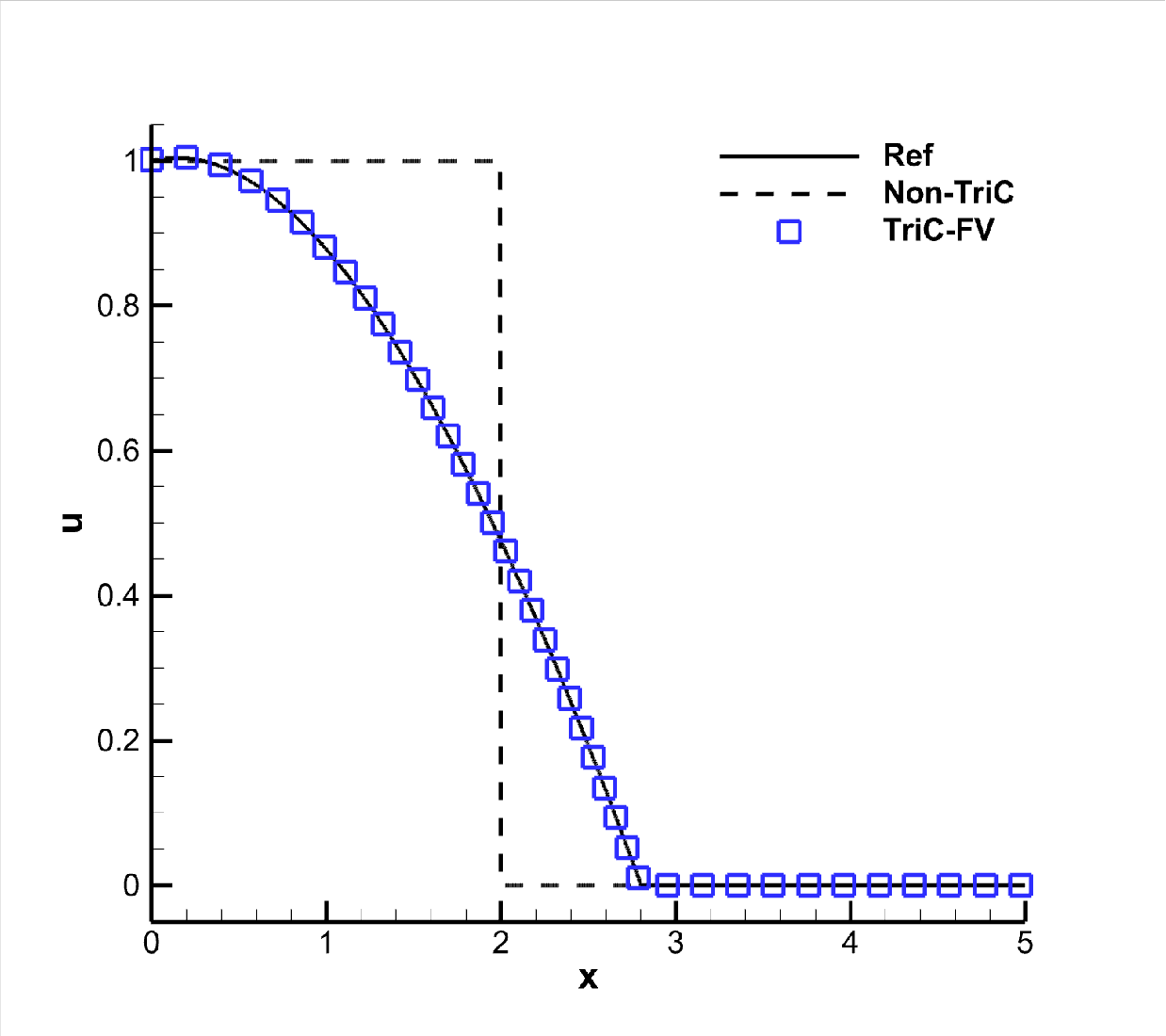}%
		\label{fig:num-step-1d-profiles}%
	}
	\caption{Comparison of the numerical solution of $K(0,3)$ equation in Example 3 and Example 4 with the TriC-FV solution and the non-TriC solution.}
\end{figure}

\paragraph{Example 3: the $K(0,3)$ equation with an exact traveling front solution}
\label{subsubsec:num-1d-front}
Consider the exact travelling front proposed in \cite{StAM_2024_Anco_Gandarias}
\begin{equation}\label{eq:num-front}
	u_t+(u^3)_{xxx}=0,\qquad x\in(0,5),
\end{equation}
with the exact solution $u(x,t)=\frac{(2+t-x)_+}{\sqrt{6}}$, where $(z)_+=\max(z,0)$.
At $T=1.5$, the front moves from $x=2$ to $x=3.5$. 
We use Picard iteration to solve this nonlinear system \eqref{eq:num-front-chain}. The Picard iteration starts from $v^{n+1,0}=v^n$ and solves
\begin{equation}\label{eq:num-chain-picard}
	\frac{v_j^{s+1,k+1}-v_j^s}{\Delta t_s}
	+D_+^2\left(q^{s+1,k}D_-v^{s+1,k+1}\right)_j=0,
	\qquad q_j^{s+1,k}=3(v_j^{s+1,k})^2.
\end{equation}
Notice that for \eqref{eq:num-front-chain}-\eqref{eq:num-chain-picard}, the Picard iteration retains $v_j^s=0$ at nodes $x_j\geq2$ and this phenomenon will persist under the mesh refinement. Therefore the front of the numerical solution calculated by the non-TriC scheme will stay at $x=2$.
The TriC-FV scheme is run on $N=100$ non-uniform cells, with $\varepsilon=0.2$ in \eqref{eq:num-mesh-map}. 
Figure~\ref{fig:num-front} shows that TriC-FV captures the exact front, while the non-TriC scheme front stays at $x=2$. 
Computation with mesh refinements confirms the numerical convergence of the numerical solutions.

\paragraph{Example 4: the $K(0,3)$ equation with discontinuous initial condition}\label{subsubsec:num-1d-front1}
We again consider the equation \eqref{eq:num-front}, with the following initial and boundary conditions:
\begin{equation}
	u(x,0)=\mathbf1_{\{x<2\}},\quad u(0,t) = 1,\ u(5,t) =  \partial_{x} u (5,t) = 0.
\end{equation}
For such an initial condition, we can easily prove that the solution of \eqref{eq:num-front-chain}-\eqref{eq:num-chain-picard} converges to itself as the mesh refines, and $u(x,t)=\mathbf1_{\{x<2\}}$ does not satisfy the ultra-weak formulation \eqref{eq:ultra-weak-formulation} of the original equation \eqref{eq:num-front}. This phenomenon is also verified by the numerical result presented in Fig.~\ref{fig:num-step-1d-profiles}.
Computation with mesh refinements confirms the numerical convergence of the numerical solutions.

\subsubsection{Compactons}
In this part, we simulate the compactons, i.e., the solitons with compact support. We also utilize this example to compare the fifth-order TriC scheme and the first-order TriC scheme to test the high-resolution of high-order schemes.
\paragraph{Propagation of a $K(3,3)$ compacton}
\label{subsubsec:num-1d-compacton}
The Rosenau--Hyman equation \cite{PRL_1993_Rosenau_Hyman} is
\begin{equation}\label{eq:num-K33}
	u_t+(u^3)_x+(u^3)_{xxx}=0.
\end{equation}
Its exact compacton is
\begin{equation}\label{eq:num-K33-compacton}
	U(z)=
	\begin{cases}
		\sqrt{3/2}\cos(z/3),& |z|\leq3\pi/2,\\
		0,& |z|>3\pi/2,
	\end{cases}
	\qquad u(x,t)=U(x-x_0-t).
\end{equation}
We take $x_0=-5$, the physical domain $[-20,20]$, and $T=2$.
TriC-FV uses non-uniform meshes.

\begin{figure}[htbp]
	\centering
	\subfigure[Numerical solution of the $K(3,3)$ equation.]{%
		\includegraphics[width=.32\textwidth]{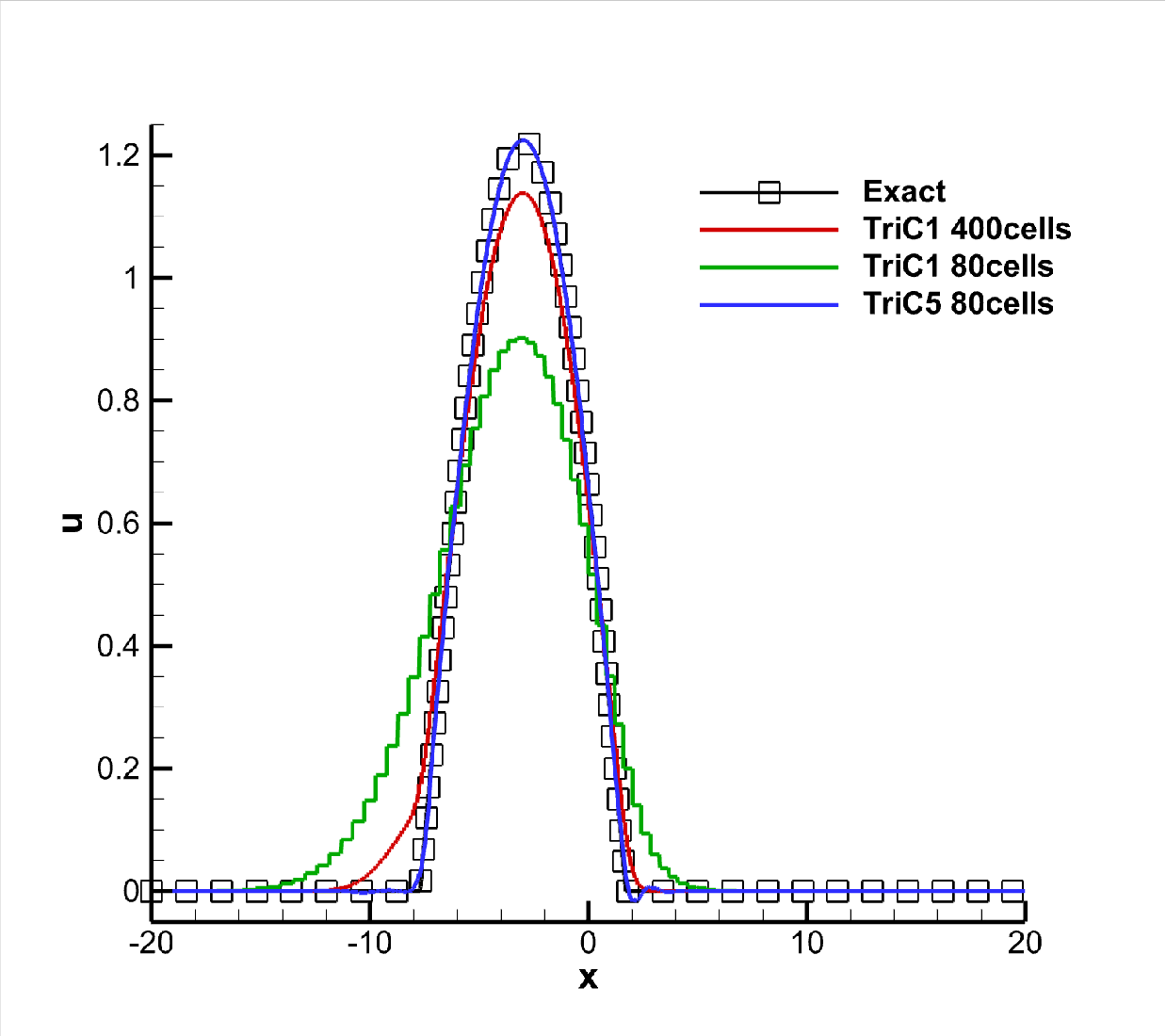}%
		\label{fig:1D_K33_compare}%
	}
	\quad
	\subfigure[Zoomed-in View of the Peak Region.]{%
		\includegraphics[width=.32\textwidth]{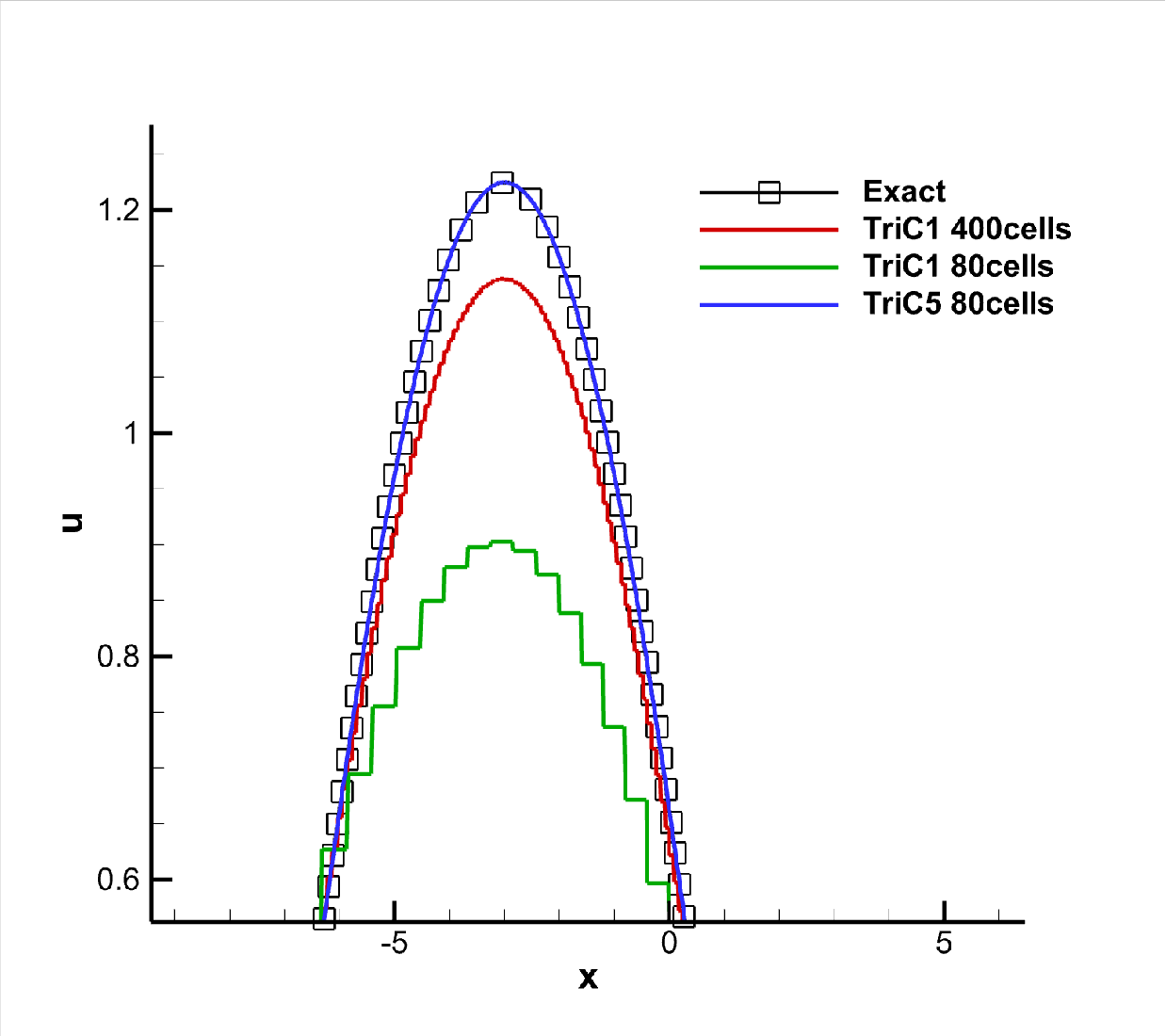}%
		\label{fig:1D_K33_peak detail}
	}
	\subfigure[Zoomed-in View of the Support Region.]{%
	\includegraphics[width=.32\textwidth]{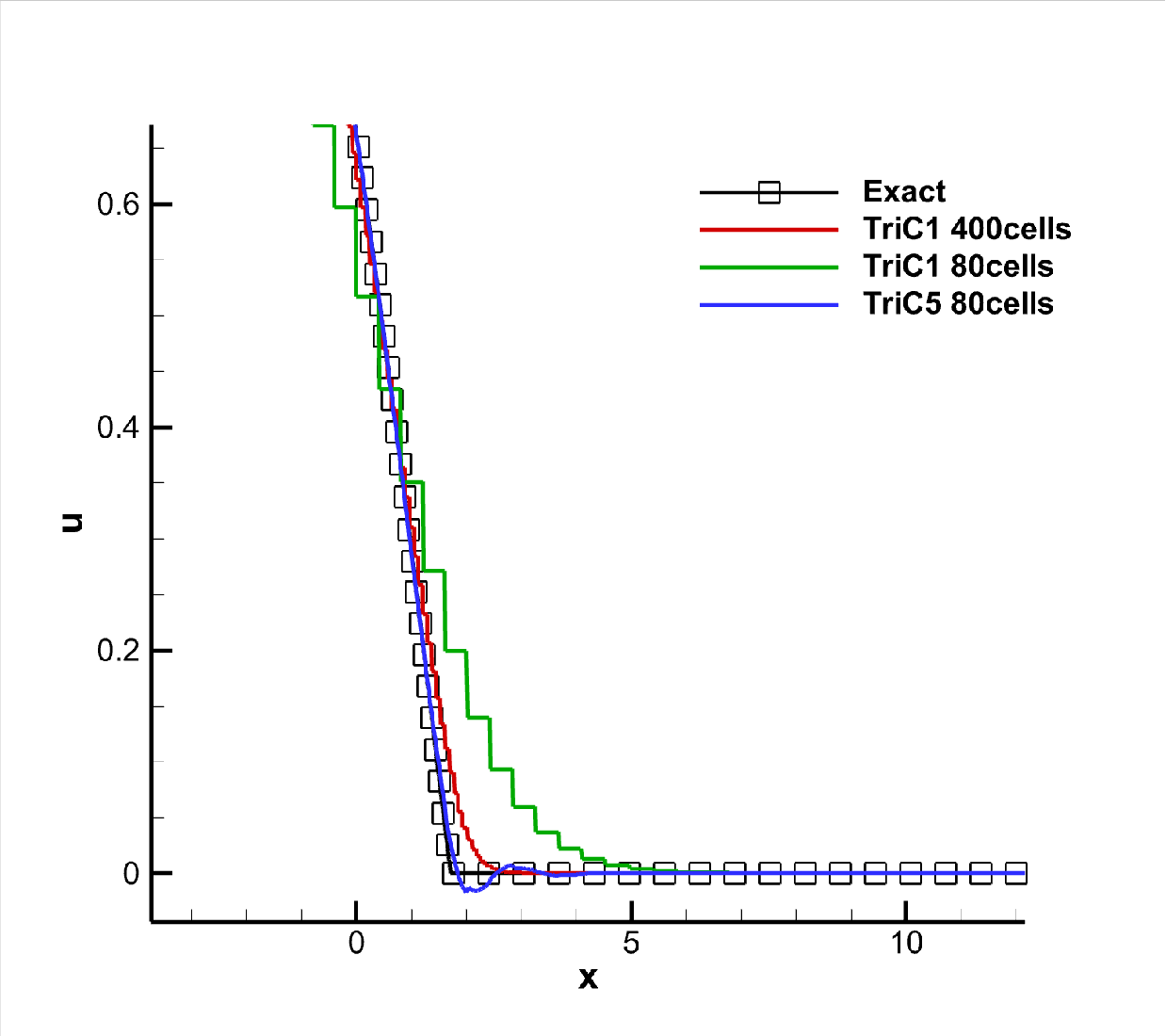}%
	\label{fig:1D_K33_support detail}%
	}
	\caption{Comparison of the high-order solution and low-order solution of \eqref{eq:num-K33-compacton} in Example 5.}
	\label{fig:2D comparison}
\end{figure}

\paragraph{Example 5: low- and high-order TriC-FV schemes}
\label{subsubsec:num-1d-low-high}
For the same compacton, we compare the fifth-order TriC-FV scheme (which reconstructs a fourth-degree polynomial in each cell) and the first-order (which reconstructs a constant in each cell) TriC-FV scheme.
Both schemes use the same triple means, numerical fluxes, and SSPRK3 temporal discretization.
We compare the example under both the same freedom and the same meshes.
Figure~\ref{fig:2D comparison} therefore includes two comparisons: equal evolved unknowns on $80$ cells, and equal freedoms with degree-four reconstruction on
$80$ cells and constant reconstruction on $400$ cells.
All meshes are non-uniform with $\varepsilon=0.2$.
The high-order reconstruction better resolves the peak and
support endpoints, including when the constant reconstruction
uses five times as many evolved unknowns.

\subsection{Two-dimensional numerical results}\label{subsec:num-2d}
In this part, we test the TriC-FV scheme for two-dimensional equations
with pure and mixed third derivatives. All results computed by the TriC-FV scheme are evaluated on the non-uniform Cartesian mesh
\begin{equation}\label{eq:num-2d-mapped-mesh}
	\begin{aligned}
		x_{i-\frac{1}{2}}&=x_L+(x_R-x_L)
		\left[\frac{i}{N_x}+\frac{\varepsilon_x}{2\pi}
		\sin\left(\frac{2\pi i}{N_x}\right)\right],
		&&0\leq i\leq N_x,\\
		y_{j-\frac{1}{2}}&=y_L+(y_R-y_L)
		\left[\frac{j}{N_y}+\frac{\varepsilon_y}{2\pi}
		\sin\left(\frac{2\pi j}{N_y}\right)\right],
		&&0\leq j\leq N_y,
	\end{aligned}
\end{equation}
with $\varepsilon_x=\varepsilon_y=0.15$.
The FD and LDG-I comparisons use uniform grids. 

\subsubsection{Accuracy for smooth solutions}\label{subsubsec:num-2d-accuracy}
\paragraph{Example 6: linear accuracy test}
The linear problem is
\begin{equation}\label{eq:num-2d-linear}
	u_t+0.6u_x+0.3u_y+u_{xxx}+0.4u_{xxy}
	+0.7u_{xyy}+0.5u_{yyy}=0,
\end{equation}
with the exact solution
\begin{equation}\label{eq:num-2d-linear-exact}
	u(x,y,t)=\sin(x+y+1.7t)+0.1\cos(2x-y+6.4t).
\end{equation}
We take $\Omega=[0,2\pi]^2$ and the final time $T=0.1$. Table~\ref{tab:num-2d-linear} shows the fifth order accuracy.

\paragraph{Example 7: nonlinear accuracy test}
For the nonlinear problem, we consider
\begin{equation}\label{eq:num-2d-nonlinear}
	\begin{aligned}
		&u_t+(u^2)_x+(u^3/3)_{xxx}+(u^3/3)_{xyy}=s(x,y,t),
	\end{aligned}
\end{equation}
with the exact solution $u(x,y,t)=1+0.2\sin(x+y-1.3t)+0.1\cos(2x-y-1.1t)$ and the corresponding source $s(x,y,t)$. We take $\Omega=[0,2\pi]^2$ and the final time $T=0.05$. 
Table~\ref{tab:num-2d-nonlinear} gives the optimal rate of the convergence. 

\begin{table}[htbp]
	\centering\small
	\setlength{\tabcolsep}{5pt}
	\caption{Smooth linear accuracy test \eqref{eq:num-2d-linear}.}
	\label{tab:num-2d-linear}
	\begin{tabular}{|c|rr|rr|rr|}
		\hline
		Mesh & $L^1$-error & order & $L^2$-error & order & $L^\infty$-error & order \\
		\hline
12$\times$12 & 9.70e-02 & -- & 1.85e-02 & -- & 1.00e-02 & -- \\
18$\times$18 & 1.57e-02 & 4.50 & 3.07e-03 & 4.43 & 1.81e-03 & 4.23 \\
24$\times$24 & 4.01e-03 & 4.74 & 7.96e-04 & 4.69 & 4.59e-04 & 4.77 \\
36$\times$36 & 5.55e-04 & 4.88 & 1.11e-04 & 4.85 & 6.52e-05 & 4.81 \\
48$\times$48 & 1.34e-04 & 4.95 & 2.69e-05 & 4.94 & 1.56e-05 & 4.97 \\
		\hline
	\end{tabular}
		\setlength{\tabcolsep}{5pt}
	\caption{Smooth nonlinear accuracy test \eqref{eq:num-2d-nonlinear}.}
	\label{tab:num-2d-nonlinear}
	\begin{tabular}{|c|rr|rr|rr|}
		\hline
		Mesh & $L^1$-error & order & $L^2$-error & order & $L^\infty$-error & order \\
		\hline
		16$\times$16 & 2.16e-02 & -- & 4.29e-03 & -- & 2.59e-03 & -- \\
		24$\times$24 & 3.19e-03 & 4.72 & 6.49e-04 & 4.66 & 3.75e-04 & 4.76 \\
		32$\times$32 & 7.86e-04 & 4.87 & 1.61e-04 & 4.84 & 9.37e-05 & 4.82 \\
		48$\times$48 & 1.06e-04 & 4.95 & 2.18e-05 & 4.94 & 1.29e-05 & 4.89 \\
		\hline
	\end{tabular}
\end{table}

\subsubsection{Comparison between TriC and non-TriC schemes}\label{compare_2d}
In this part we compare the TriC-FV scheme with a non-TriC finite difference scheme for degenerate dispersive equations. The LDG method provides a reference solution.
These tests examine whether a scheme can approach an incorrect limit under refinement and verify the validity of TriC schemes.

\begin{figure}[htbp]
	\centering
	\subfigure[Exact.]{%
		\includegraphics[width=.32\textwidth]{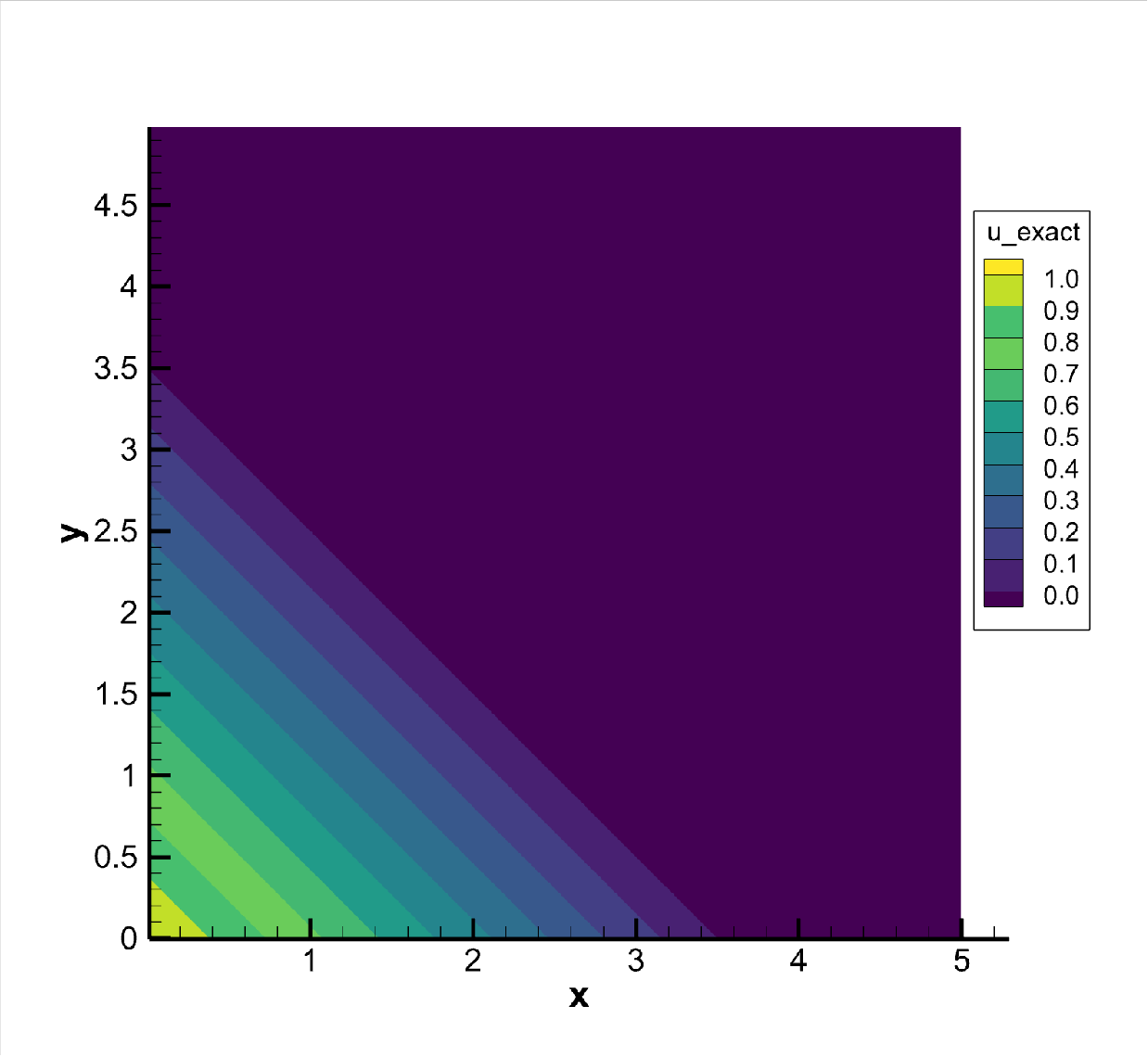}%
		\label{fig:2D_mixed_Exact}%
	}
	\quad
	\subfigure[TriC-FV solution.]{%
		\includegraphics[width=.32\textwidth]{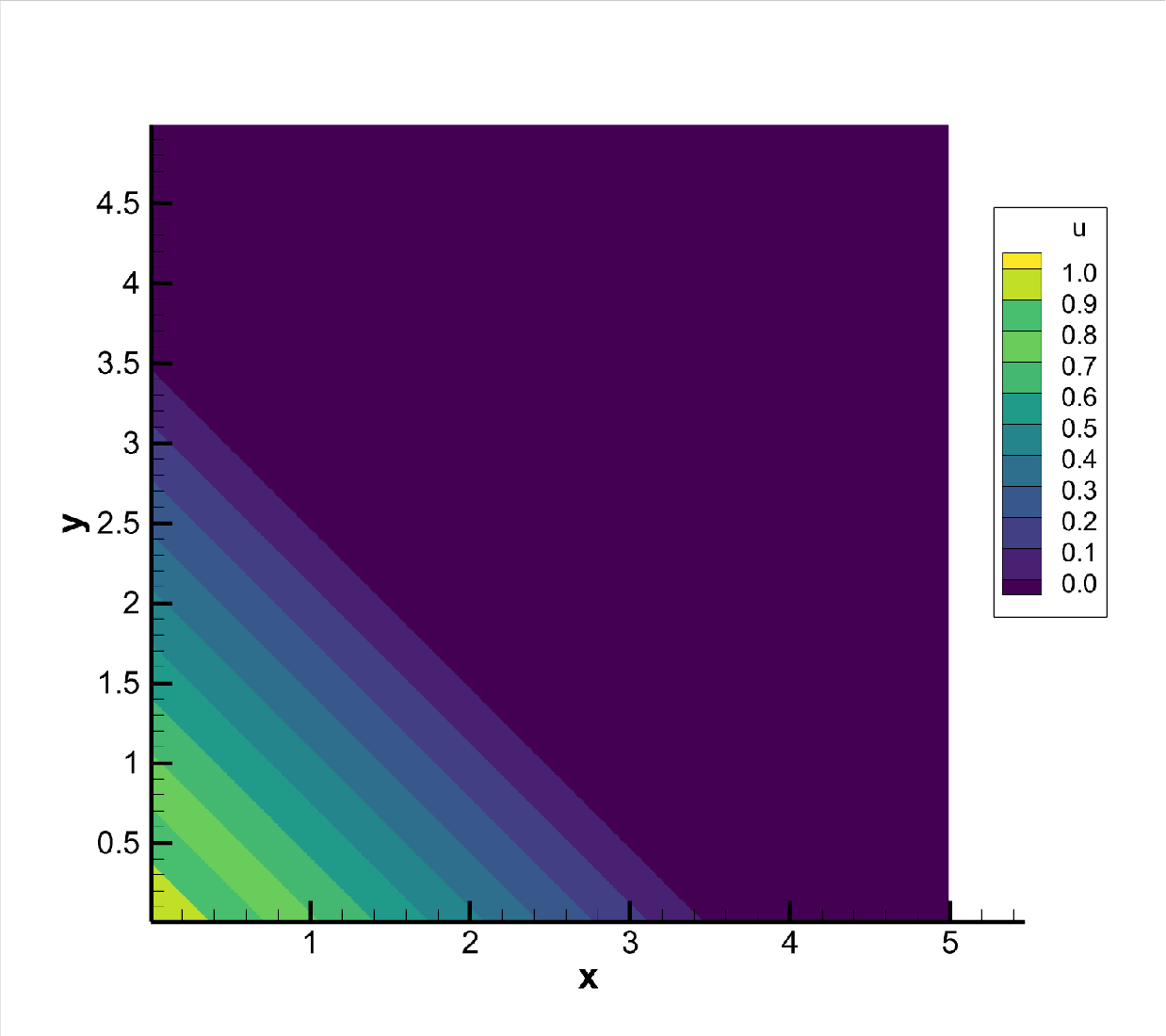}%
		\label{fig:2D_mixed_2d_TriC-FV}
	}
	\subfigure[Non-TriC solution.]{%
		\includegraphics[width=.32\textwidth]{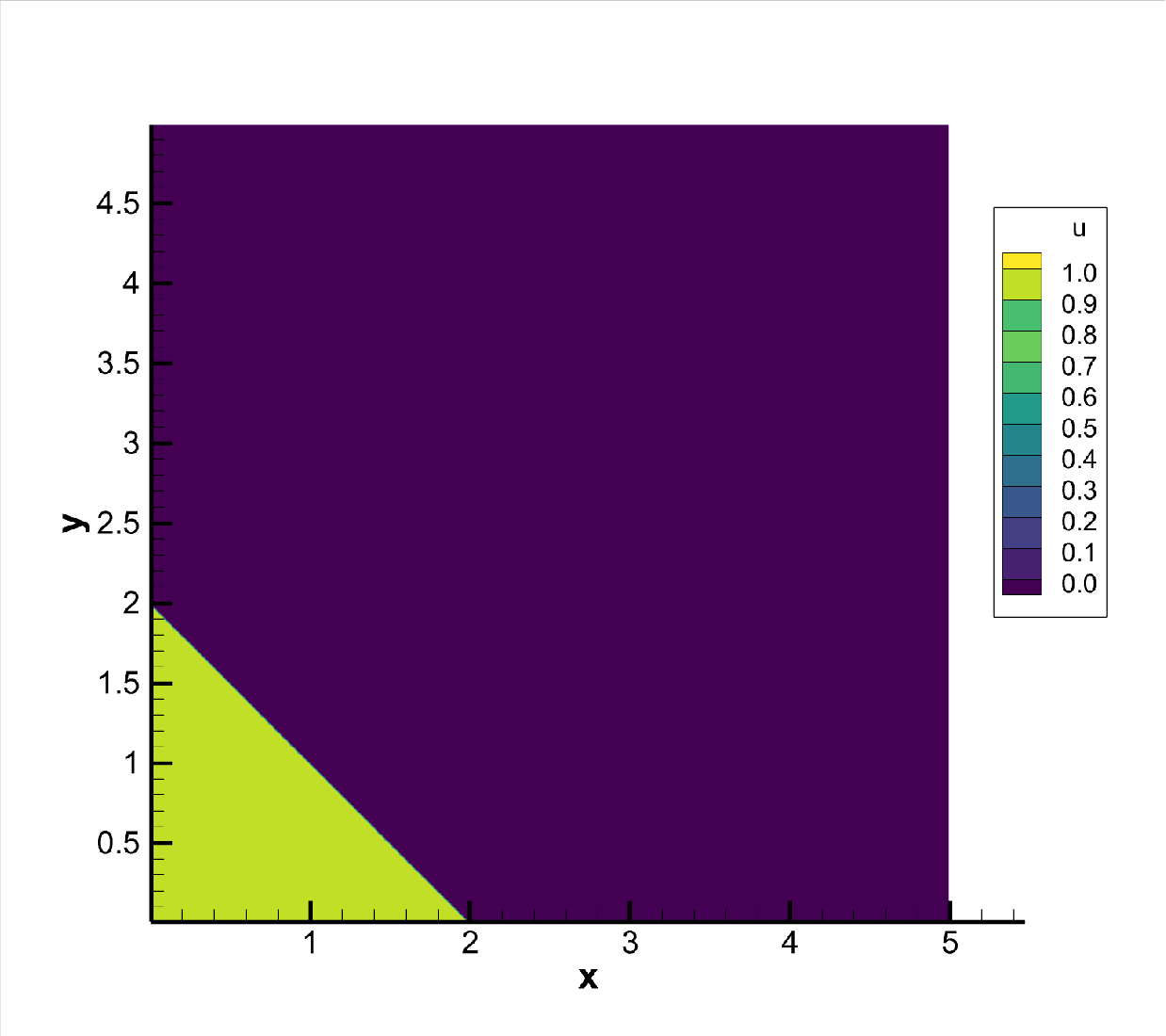}%
		\label{fig:2D_mixed_2d_NonTriC}%
	}	
	\caption{Comparison of contours of the results of \eqref{eq:num-front-2d-pde} generated by the TriC-FV scheme and the non-TriC scheme in Example 8.} 
	\label{fig:2D_mixed_compare}
\end{figure}

\paragraph{Example 8: an exact travelling front}\label{subsubsec:num-2d-front}
Here we test an anisotropic problem which has a mixed derivative. 
The equation is close to the one-dimensional problem:
\begin{equation}\label{eq:num-front-2d-pde}
	u_t+(u^3)_{xxx}+(u^3)_{xyy}=0,
	\qquad (x,y)\in(0,5)^2,
\end{equation}
with the exact solution $u(x,y,t)=\frac{(2+t-x-y)_+}{\sqrt{12}}$.
The non-TriC scheme is
\begin{equation}\label{eq:num-front-2d-chain}
	\begin{aligned}
		\frac{v_{ij}^{s+1}-v_{ij}^s}{\Delta t_n}+D_{+,x}^{2}\left(3(v^{s+1})^2D_{-,x}v^{s+1}\right)_{ij}+D_{+,x}D_{+,y}\left(3(v^{s+1})^2D_{-,y}v^{s+1}\right)_{ij}=0,
	\end{aligned}
\end{equation}
where $D_{\pm,x}$ denotes the $x$-directional forward or backward difference operator, and $D_{\pm,y}$ is defined analogously.
The non-TriC scheme is run on $160^{2}$ cells and TriC-FV scheme uses $80^2$ non-uniform cells.
Figures~\ref{fig:2D_mixed_compare} and \ref{fig:mixed_section} show that TriC-FV scheme can capture the correct moving front, while the non-TriC solution remains zero in the propagation region.
Computation with mesh refinements confirms the numerical convergence of the solutions.

\subsubsection{Compactons}
In this part we utilize the TriC-FV scheme to simulate the two-dimensional compacton proposed in \cite{PRL_2007_Rosenau_Hyman_Staley}:
\begin{equation}\label{eq:num-PRL-model}
	u_t+(u^2)_x+\frac{1}{3}\partial_x\Delta(u^3)=0,
	\qquad \Delta=\partial_{xx}+\partial_{yy}.
\end{equation}
It belongs to \eqref{2D TC model} with $f_1=u^2$,
$A_1=A_3=u^3/3$, and $f_2=A_2=A_4=0$.

\begin{figure}[htbp]
	\centering
	\subfigure[Comparison of \eqref{eq:num-front-2d-pde} along the cut $y=0.5$ in Example 8]{%
		\includegraphics[width=.48\textwidth]{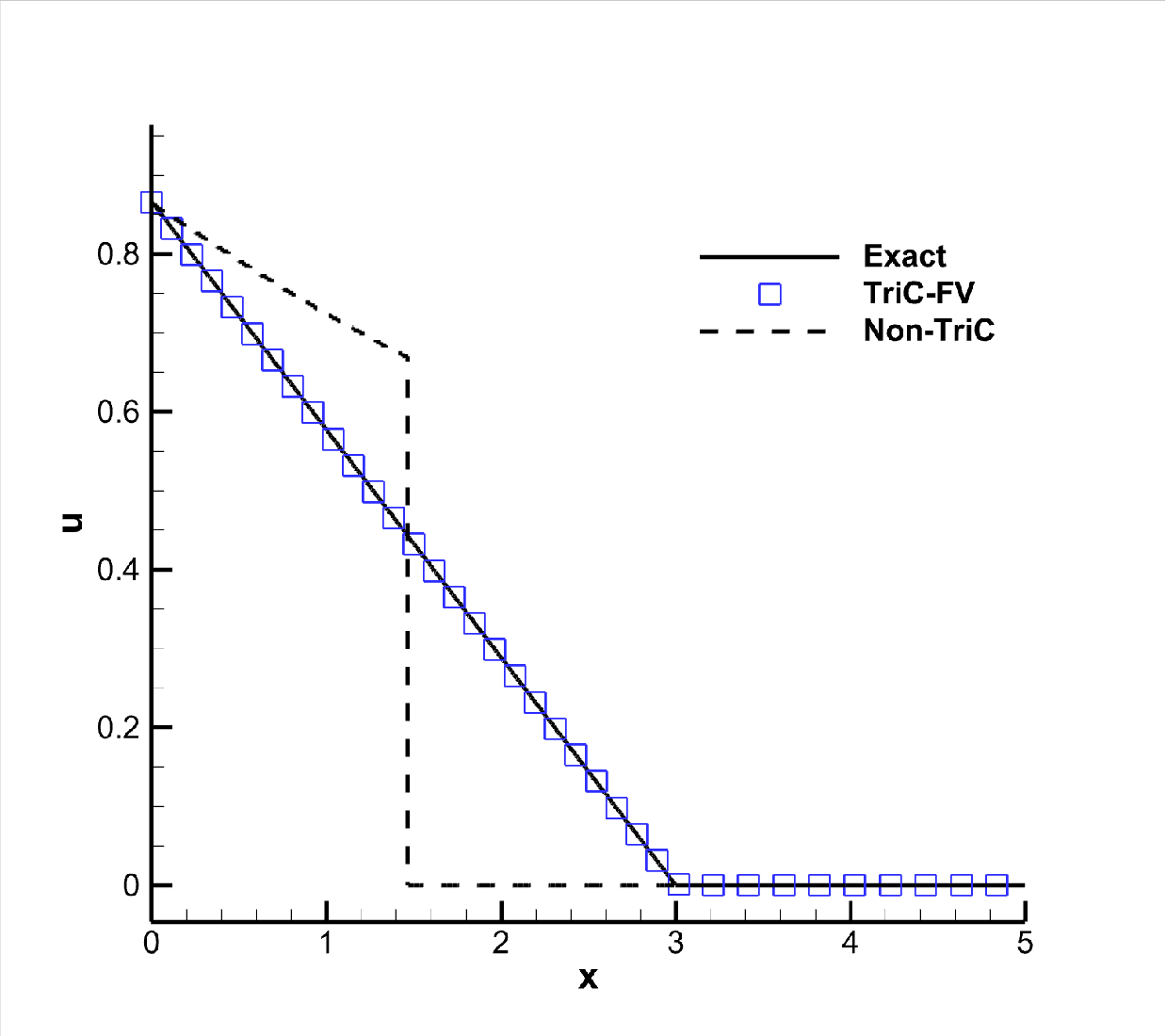}%
		\label{fig:mixed_section}%
	}
	\quad
	\subfigure[Comparison of \eqref{eq:num-parabolic-compacton} along the cut $y=0$ in Example 9.]{%
		\includegraphics[width=.48\textwidth]{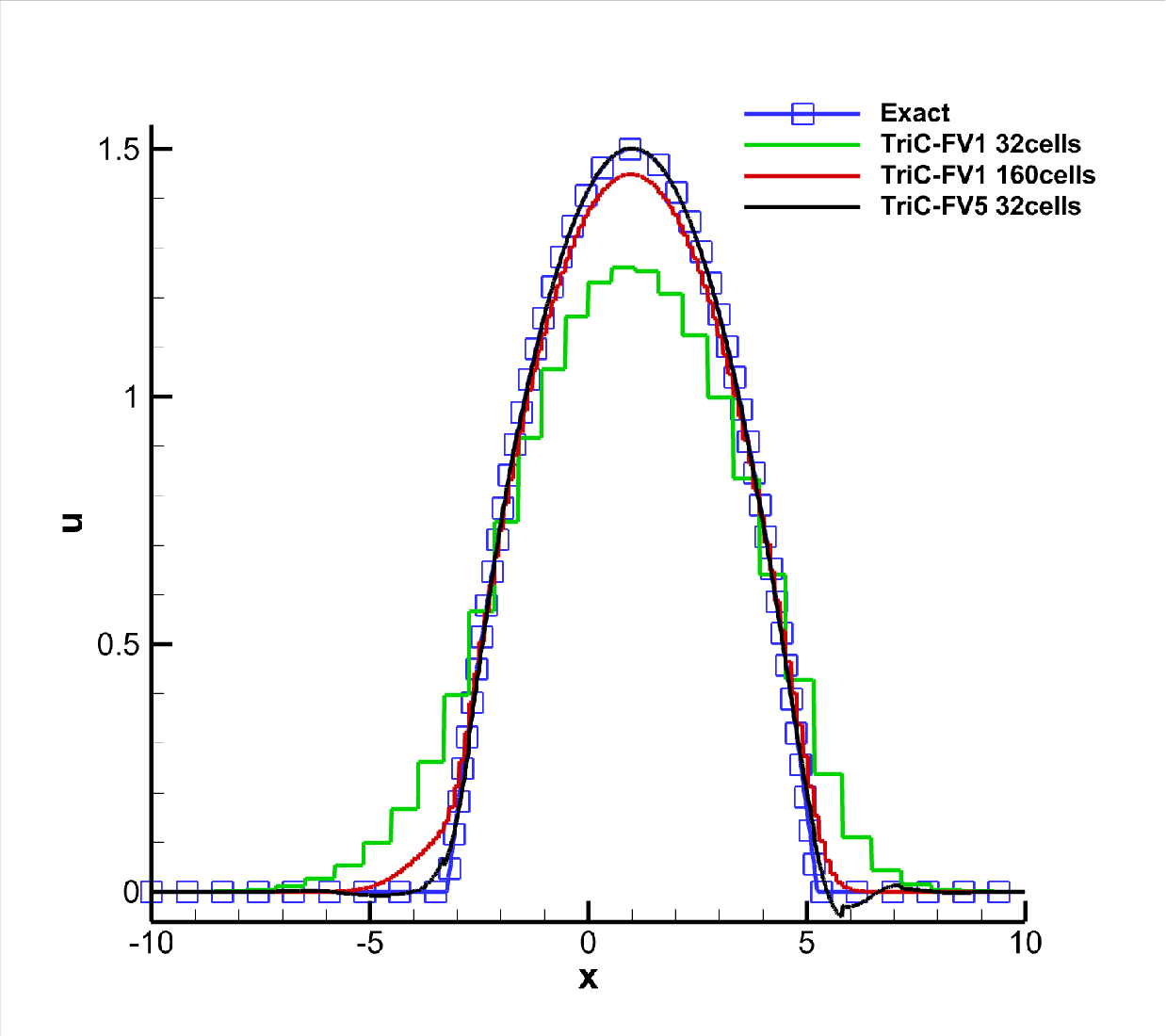}%
		\label{fig:2D_mixed_2d_TriC-FV }
	}
	\caption{Comparisons of the problems \eqref{eq:num-front-2d-pde} and \eqref{eq:num-parabolic-compacton} along the cuts in one dimension.} 
	\label{fig:2D_section_compare}
\end{figure}

\begin{figure}[htbp]
	\centering
	\subfigure[Exact solution]{%
		\includegraphics[width=.48\textwidth]{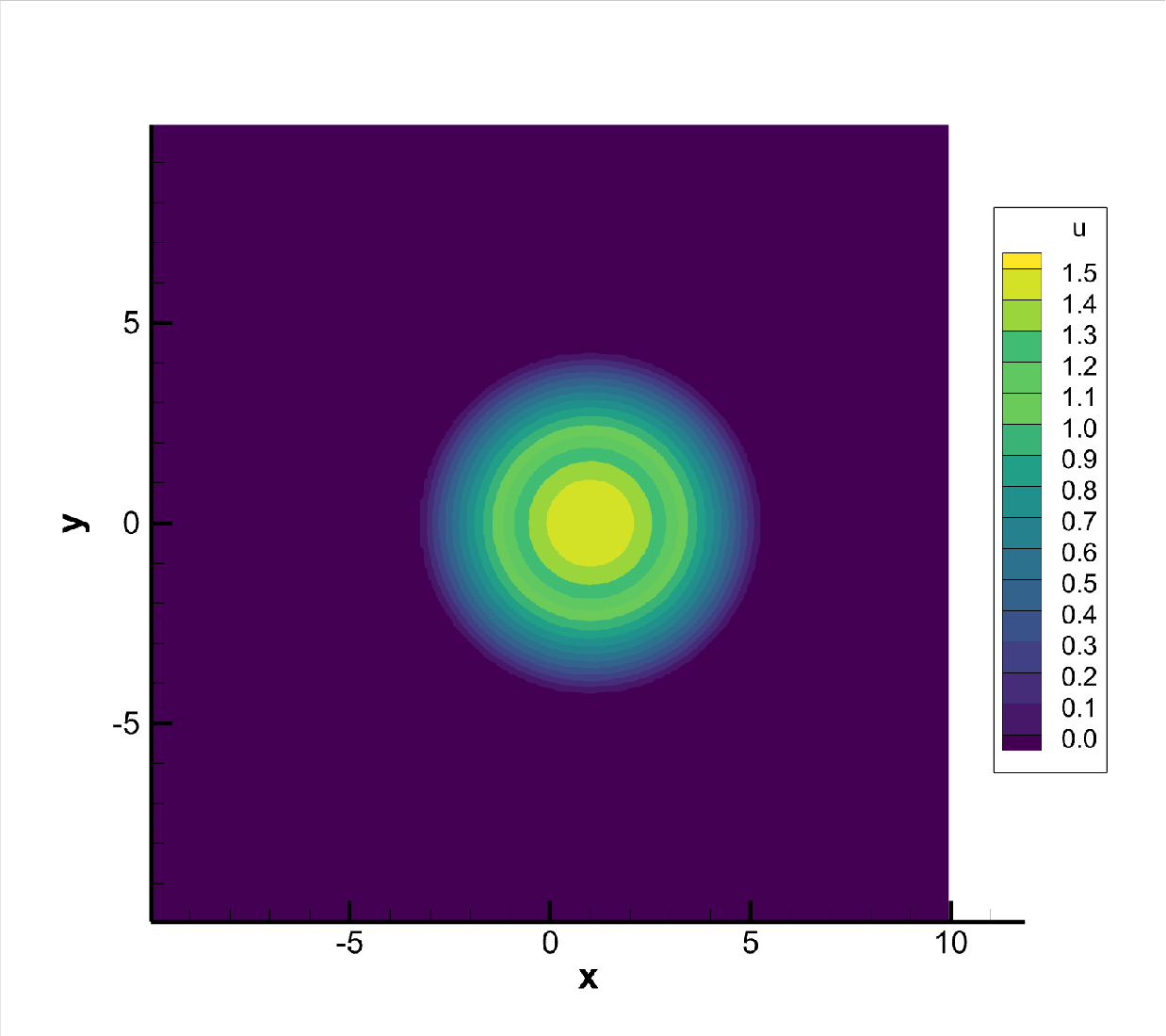}%
		\label{fig:2D_compacton_2d_exact}%
	}
	\subfigure[TriC-FV5 scheme on $32^{2}$ cells.]{%
		\includegraphics[width=.48\textwidth]{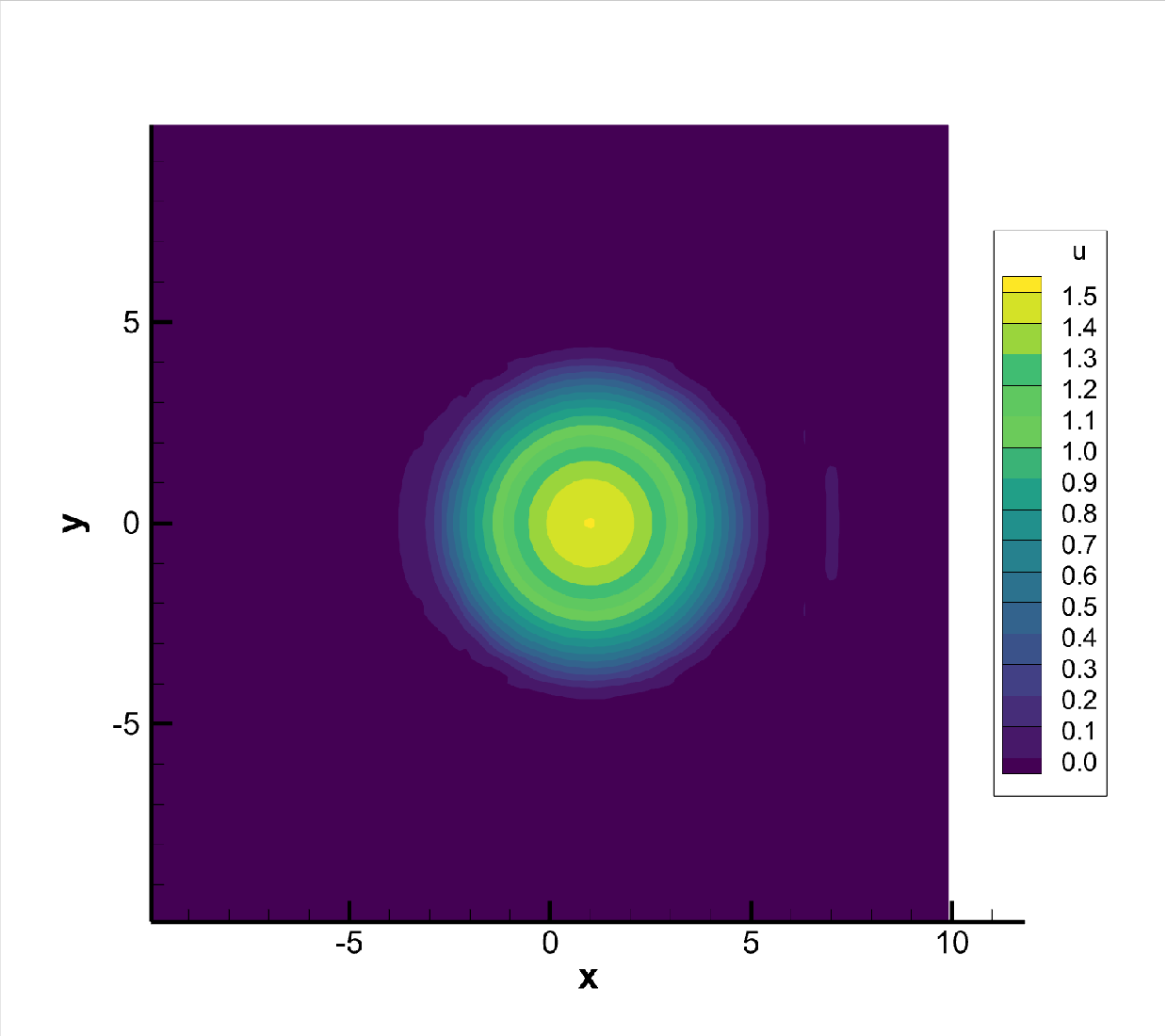}%
		\label{fig:2D_compacton_2d_FV5_32_32}
	}\\
	\subfigure[TriC-FV1 scheme on $32^{2}$ cells.]{%
		\includegraphics[width=.48\textwidth]{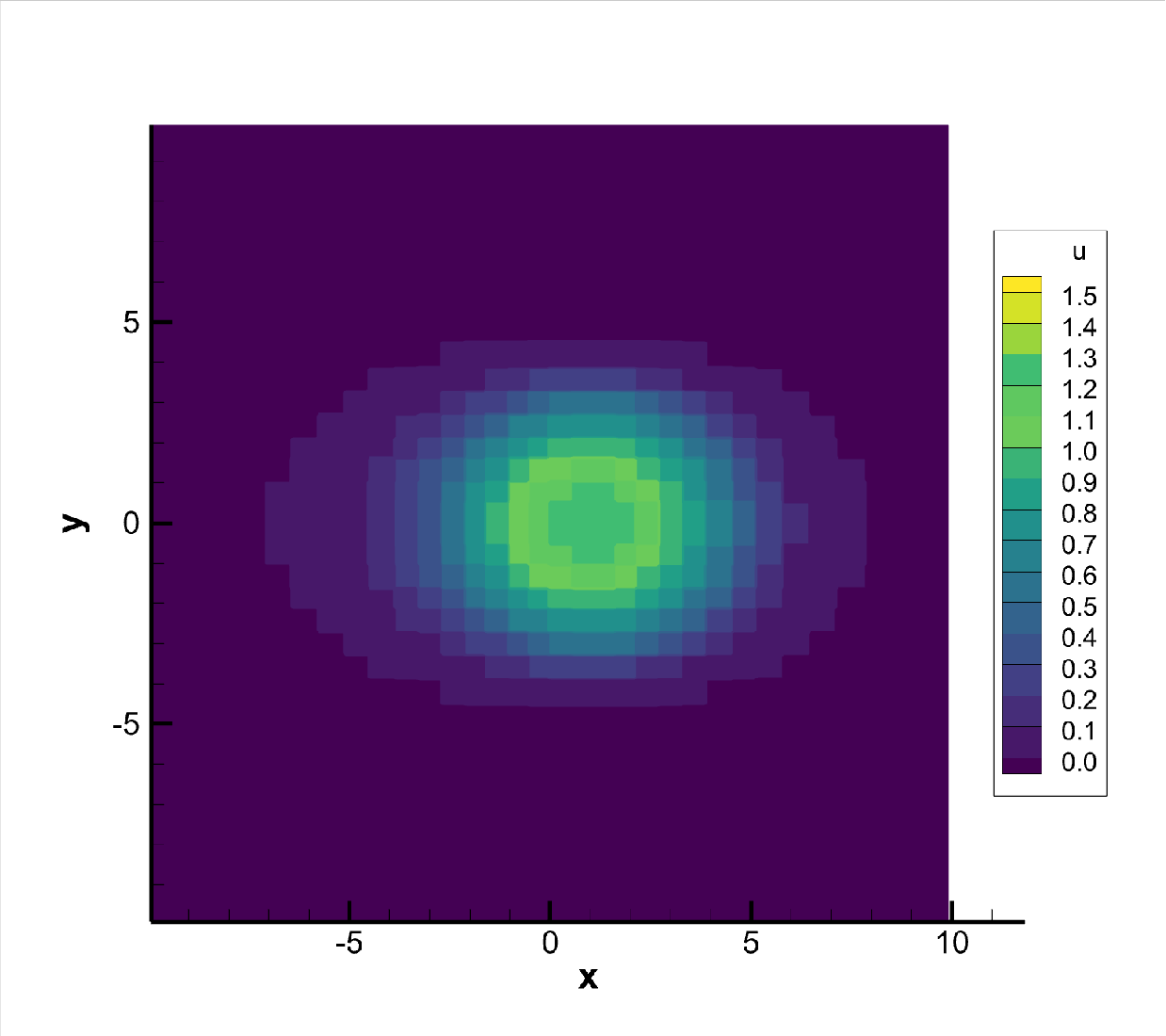}%
		\label{fig:2D_compacton_2d_FV1_32_32}%
	}
	\subfigure[TriC-FV1 scheme on $160^{2}$ cells.]{%
		\includegraphics[width=.48\textwidth]{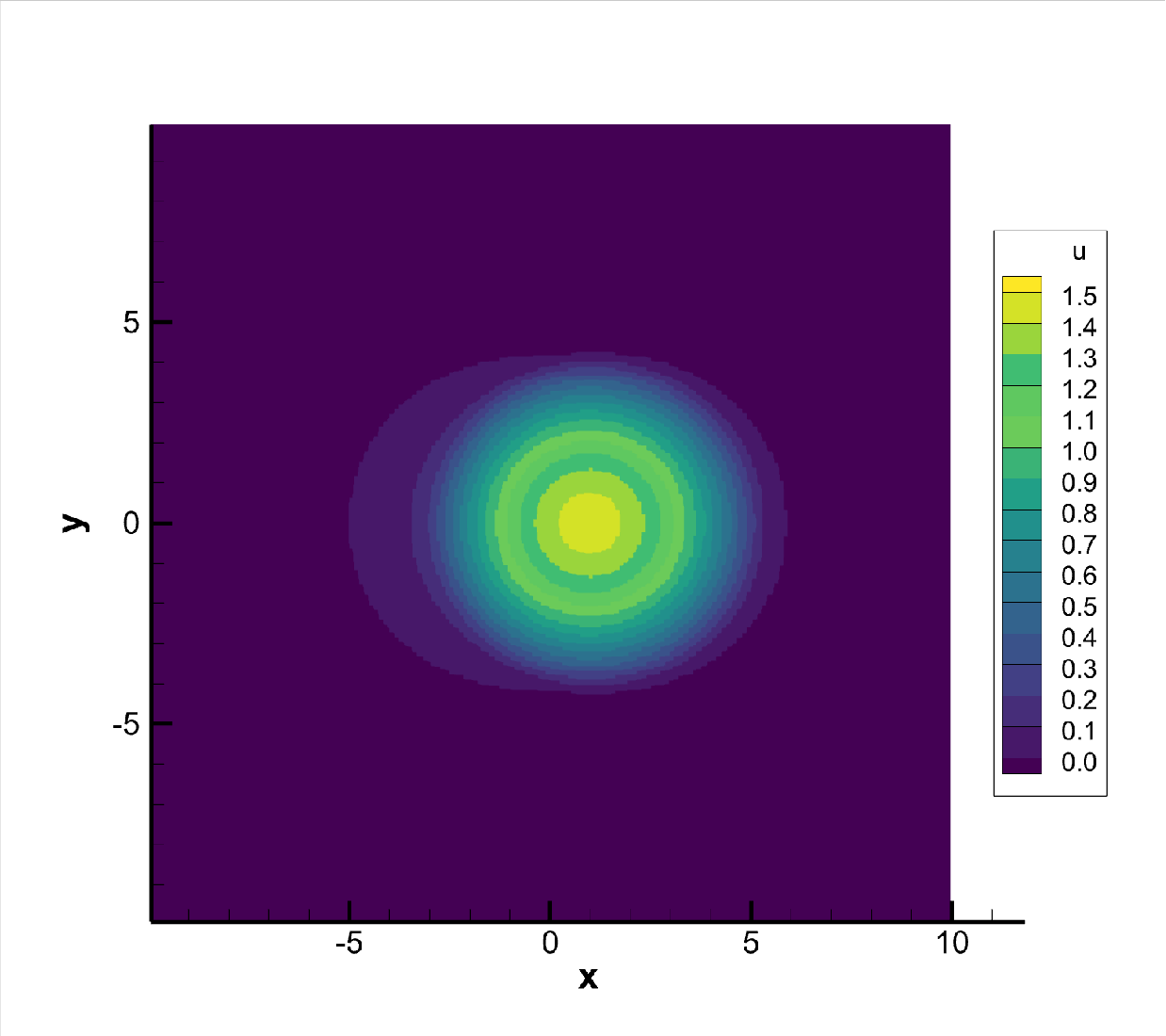}%
		\label{fig:2D_compacton_2d_FV1_160_160}%
	}
	\caption{Comparison of contours of the numerical solutions of the problem \eqref{eq:num-parabolic-compacton} generated by the high-order and low-order TriC-FV schemes in Example 9.} 
	\label{fig:2D_compacton}
\end{figure}

\paragraph{Example 9: comparison of the high-order and low-order schemes}\label{subsubsec:num-2d-compacton}
The equation~\eqref{eq:num-PRL-model} has the exact solution if we take the corresponding initial and boundary conditions:
\begin{equation}\label{eq:num-parabolic-compacton}
	u(x,y,t)=\left[
	\frac{3}{2}-\frac{(x-x_0-t)^2+(y-y_0)^2}{12}
	\right]_+.
\end{equation}
We take $(x_0,y_0)=(0,0)$, $\Omega=[-10,10]^2$, and the final time
$T=1$. The fifth-order TriC-FV (TriC-FV5) scheme uses $32^2$ non-uniform cells and the first-order TriC (TriC-FV1) scheme with the piecewise constant reconstruction uses $32^{2}$ cells and $160^{2}$ cells, at $T=1$.
Figure~\ref{fig:2D_compacton} compares the fifth-order result on
$32^2$ cells with constant reconstructions on $32^2$ and $160^2$
cells, corresponding to equal cells and equal freedoms. 
Their one-dimensional comparison along the cut $y=0$ is presented in Figure~\ref{fig:2D_section_compare}.
The results show that even for this non-smooth solution, high-order results better resolve the solution. 

\begin{figure}[htbp]
	\centering
	\subfigure[$T=0$ and $T=4$]{%
		\includegraphics[width=.48\textwidth]{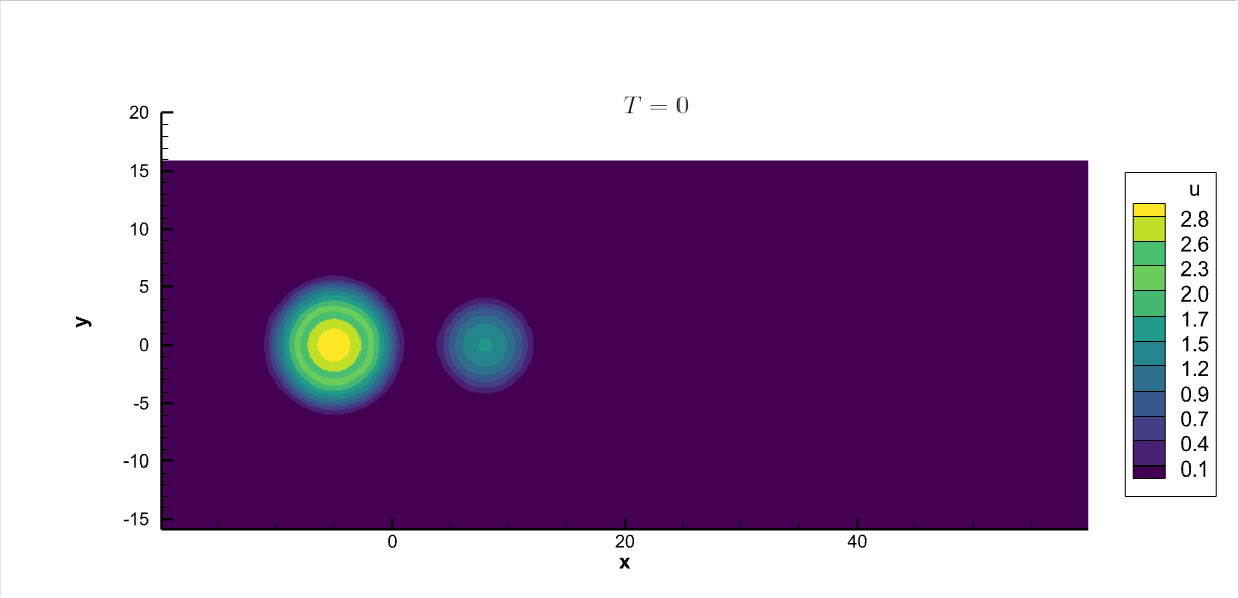}\quad
		\includegraphics[width=.48\textwidth]{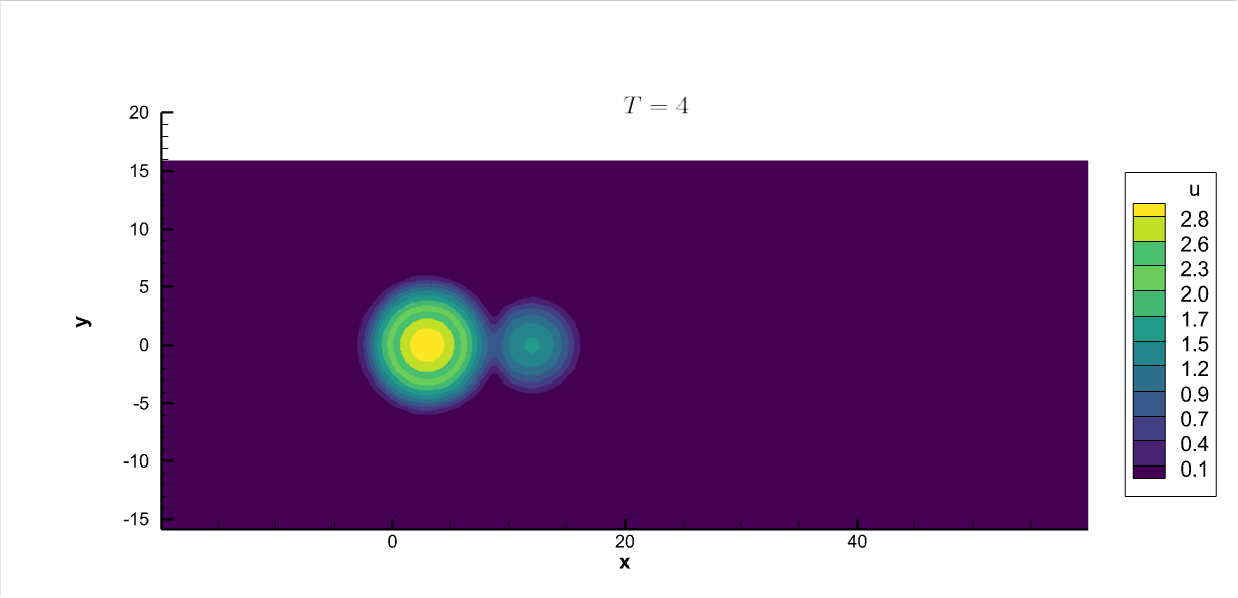}%
		\label{fig:2D_two_compactons_TriC-FV1}%
	}\\
	\subfigure[$T=8$ and $T=12$]{%
		\includegraphics[width=.48\textwidth]{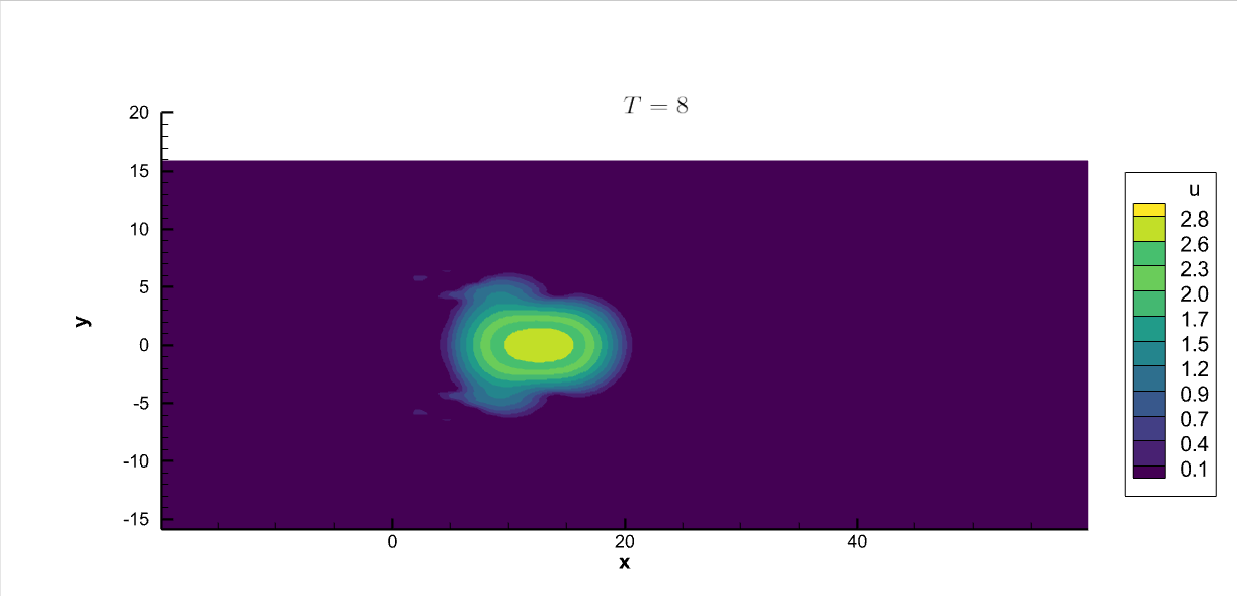}\quad
		\includegraphics[width=.48\textwidth]{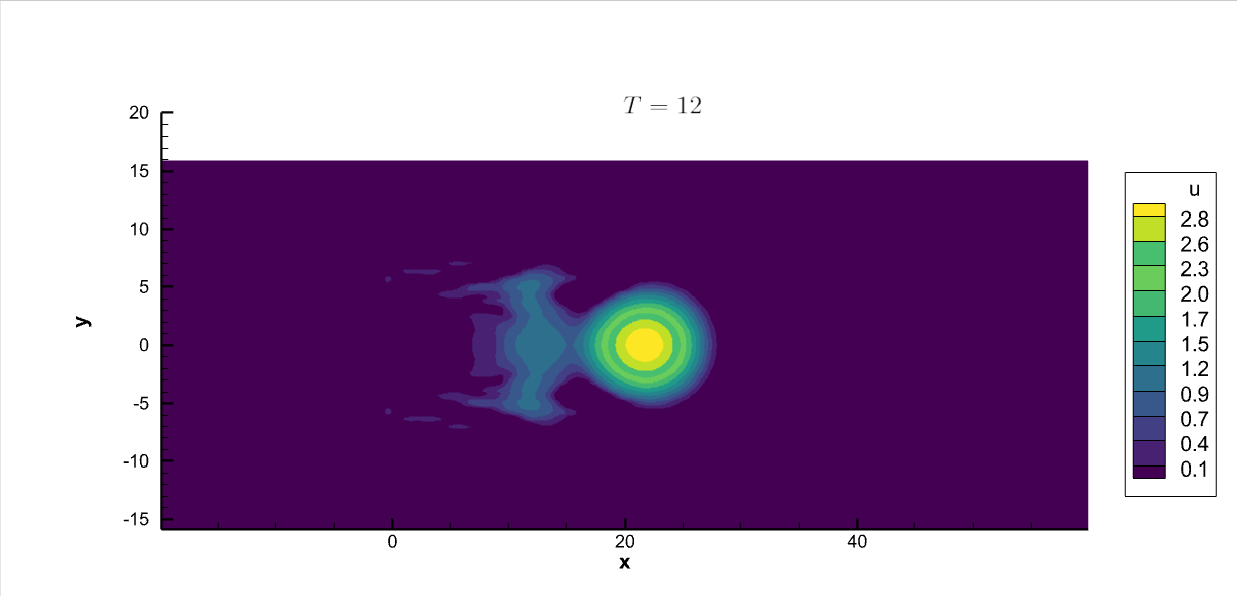}%
		\label{fig:2D_two_compactons_TriC-FV2}
	}\\ 
	\subfigure[$T=16$ and $T=20$]{%
		\includegraphics[width=.48\textwidth]{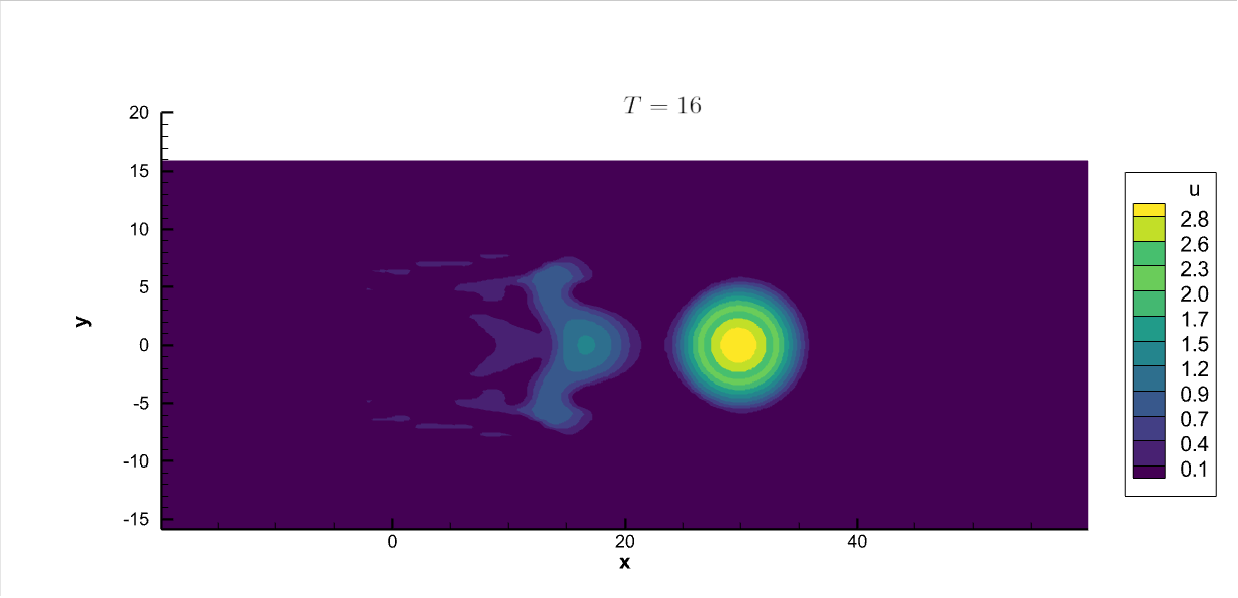}\quad
		\includegraphics[width=.48\textwidth]{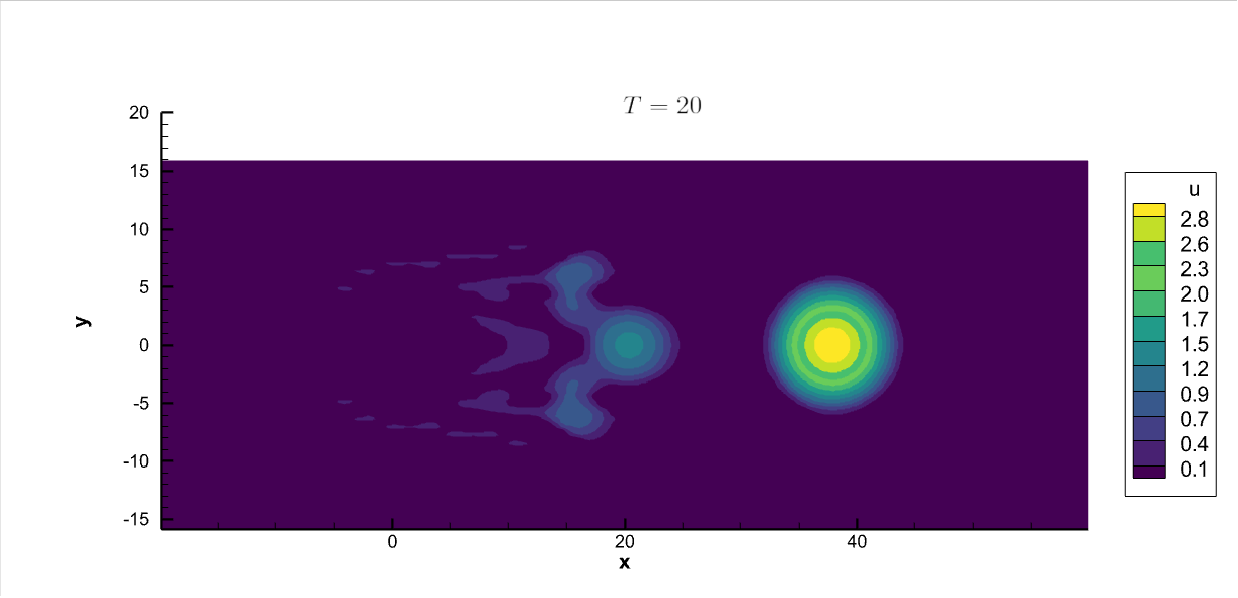}%
		\label{fig:2D_two_compactons_TriC-FV3}%
	}
	\caption{Two-compacton interaction of \eqref{eq:num-two-compactons} in Example 10 on $[-20,60]\times[-16,16]$ at different times.} 
	\label{fig:num-interaction}
\end{figure}

\paragraph{Example 10: interaction of two compactons}\label{subsubsec:num-interaction}
For \eqref{eq:num-PRL-model}, the initial condition is taken as
\begin{equation}\label{eq:num-two-compactons}
	u_0(x,y)
	=\left(3-\frac{(x+5)^2+y^2}{12}\right)_+
	+\left(\frac{3}{2}-\frac{(x-8)^2+y^2}{12}\right)_+,
\end{equation}
with two isolated compactons having speeds $2$ and $1$. This example can model the collision of two compactons. 
We take the physical domain $[-20,60]\times[-16,16]$ and compute through $T=20$ to follow the interaction and subsequent separation.
The results of the TriC-FV scheme on $96\times48$ non-uniform cells are presented in Figure~\ref{fig:num-interaction}. 

\begin{figure}[htbp]
	\centering
	\begin{minipage}{0.48\textwidth}
		\centering
		\includegraphics[width=\linewidth]{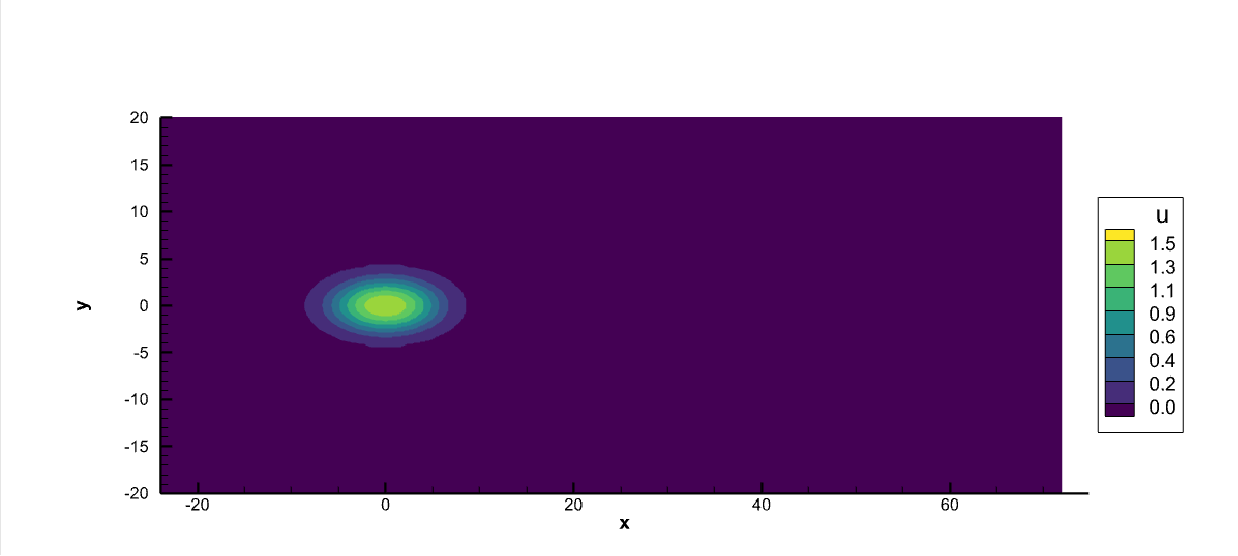}
		
		\small $T=0$
	\end{minipage}
	\hfill
	\begin{minipage}{0.48\textwidth}
		\centering
		\includegraphics[width=\linewidth]{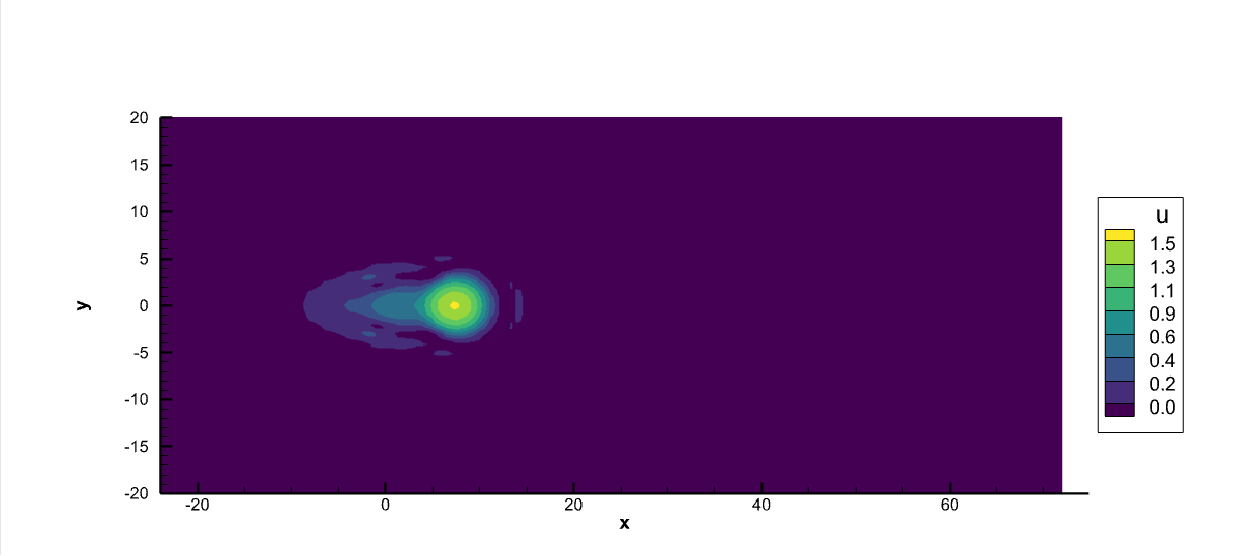}
		
		\small $T=6$
	\end{minipage}
	
	\begin{minipage}{0.48\textwidth}
		\centering
		\includegraphics[width=\linewidth]{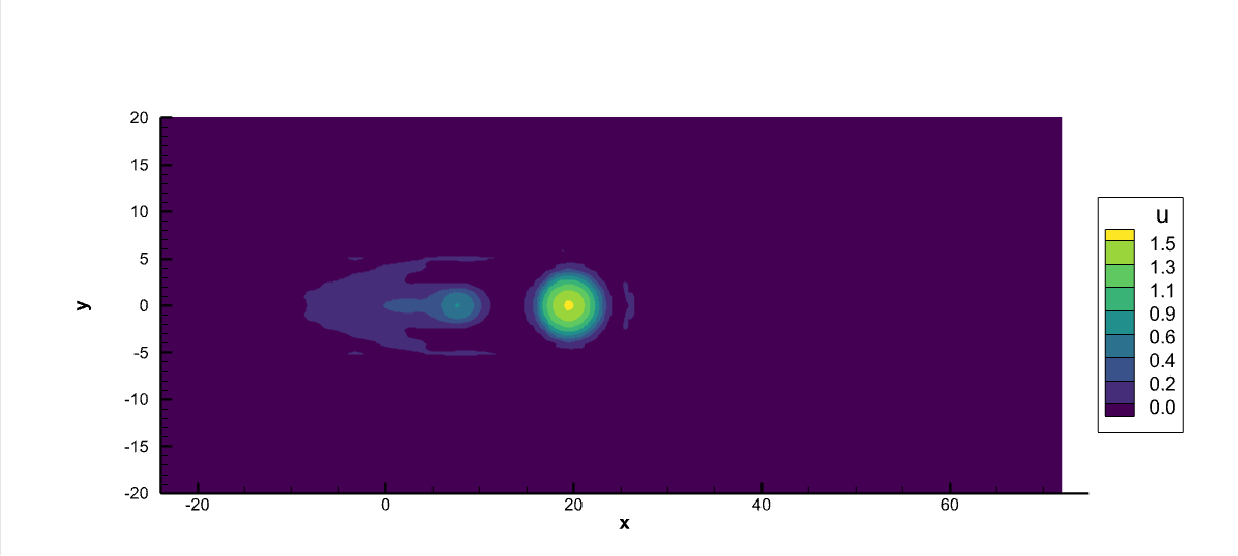}
		
		\small $T=18$
	\end{minipage}
	
	\caption{Compacton of \eqref{eq:num-prl-C203-pulse-data} emerging from a compactly supported initial data in Example 11.}
	\label{fig:num-pulse-t18-evolution}
\end{figure}

\paragraph{Example 11: formation from an elongated pulse}\label{subsubsec:num-elongated-pulse}
We next consider the formation experiment for $C_2(2,0+3)$ in
\cite{PRL_2007_Rosenau_Hyman_Staley}.
This experiment could illustrate compactons gradually emerging from a compactly supported initial data.
We take the following initial condition for \eqref{eq:num-PRL-model}:
\begin{equation}\label{eq:num-prl-C203-pulse-data}
	\begin{aligned}
		\rho(x,y)&=\left[\left(\frac{x-x_c}{R_x}\right)^2
		+\left(\frac{y-y_c}{R_y}\right)^2\right]^{1/2},\quad 
		u_0(x,y)=
		\begin{cases}
			U_0\cos^2\!\left(\dfrac{\pi\rho(x,y)}2\right),&\rho\leq1,\\
			0,&\rho>1,
		\end{cases}
	\end{aligned}
\end{equation}
with $U_0=1.5,\ R_x=9,\ R_y=4.5,\ (x_c,y_c)=(0,0)$.
We take $\Omega=[-24,72]\times[-24,24]$ and the final time $T=18$.
The TriC-FV scheme uses $128\times64$ non-uniform cells.
Figure~\ref{fig:num-pulse-t18-evolution} shows compactons being continuously emitted to the right.

\section{Conclusion}\label{Sec5}
To solve the problem that a numerical solution generated by a scheme that is both conservative and doubly conservative may still converge to a non-ultra-weak solution of the degenerate convection-dispersion equation, we have proposed a triply conservative theory for numerically solving the degenerate convection-dispersion equation both in one and two space dimensions. 
We further developed a Lax-Wendroff type theorem for the degenerate convection-dispersion equations to provide a new theoretical guarantee to numerical schemes for such equations.
We then constructed a new high-order triply conservative finite volume scheme on non-uniform meshes based on a deep connection between the triple cell average and the quadratic B-spline. We also presented some high-order triply conservative schemes proposed in the literature.
Numerical experiments both in one dimension and two dimensions verified the triply conservative theory and illustrated the better resolution of the high-order TriC-FV scheme. 
We also used the new scheme to simulate the two-dimensional anisotropic problems proposed in the literature \cite{PRL_2007_Rosenau_Hyman_Staley} and observed the same phenomenon.

\bibliographystyle{amsplain_initials}

\begin{thebibliography}{99}
	
	\bibitem{MCS_2023_Ahmat_Qiu}
	M.~Ahmat and J.~Qiu, \emph{Direct WENO scheme for dispersion-type equations},
	Mathematics and Computers in Simulation \textbf{204} (2023), 216--229.
	
	\bibitem{StAM_2024_Anco_Gandarias}
	S.~C. Anco and M.~Gandarias, \emph{Weak compactons of nonlinearly dispersive
		KdV and KP equations}, Studies in Applied Mathematics \textbf{153} (2024),
	no.~4, e12777.
	
	\bibitem{Book_Brenner_Scott}
	S.~C. Brenner and L.~R. Scott, \emph{The mathematical theory of finite element
		methods}, Springer, 2008.
	
	\bibitem{MCOM_2008_ultraweak}
	Y.~Cheng and C.-W. Shu, \emph{A discontinuous Galerkin finite element method
		for time dependent partial differential equations with higher order
		derivatives}, Mathematics of Computation \textbf{77} (2008), 699--730.
	
	\bibitem{JCP_2001_Chertock_Levy}
	A.~Chertock and D.~Levy, \emph{Particle methods for dispersive equations},
	Journal of Computational Physics \textbf{171} (2001), no.~2, 708--730.
	
	\bibitem{Book_de_Boor}
	C.~De~Boor, \emph{A practical guide to splines}, vol.~27, Springer, New York,
	1978.
	
	\bibitem{JCP_1995_deFrutos_Lopez_Sanz}
	J.~deFrutos, M.~L{\'o}pez-Marcos, and J.~M. Sanz-Serna, \emph{A finite
		difference scheme for the $K(2,2)$ compacton equation}, Journal of
	Computational Physics \textbf{120} (1995), no.~2, 248--252.
	
	\bibitem{CPAM_2019_Germain_Benjamin_Jeremy}
	P.~Germain, B.~Harrop-Griffiths, and J.~L. Marzuola, \emph{Existence and
		uniqueness of solutions for a quasilinear KdV equation with degenerate
		dispersion}, Communications on Pure and Applied Mathematics \textbf{72}
	(2019), no.~11, 2449--2484.
	
	\bibitem{SIREV_2001_Gottlieb_Shu_Tadmor}
	{S. Gottlieb}, {C.-W. Shu}, and {E. Tadmor}, \emph{Strong stability-preserving
		high-order time discretization methods}, SIAM Review \textbf{43} (2001),
	no.~1, 89--112.
	
	\bibitem{JCP_1999_Holden_Karlsen_Risebro}
	H.~Holden, K.~H. Karlsen, and N.~H. Risebro, \emph{Operator splitting methods
		for generalized Korteweg-de Vries equations}, Journal of Computational
	Physics \textbf{153} (1999), no.~1, 203--222.
	
	\bibitem{SISC_2022_Huang_Cheng}
	J.~Huang, Y.~Liu, Y.~Liu, Z.~Tao, and Y.~Cheng, \emph{A class of adaptive
		multiresolution ultra-weak discontinuous Galerkin methods for some nonlinear
		dispersive wave equations}, SIAM Journal on Scientific Computing \textbf{44}
	(2022), no.~2, A745--A769.

	
	\bibitem{KdV_original}
	D.~J. Korteweg and G.~De~Vries, \emph{On the change of form of long waves
		advancing in a rectangular canal, and on a new type of long stationary
		waves}, The London, Edinburgh, and Dublin Philosophical Magazine and Journal
	of Science \textbf{39} (1895), no.~240, 422--443.
	
	\bibitem{Lax_Wendroff_1960}
	P.~Lax and B.~Wendroff, \emph{Systems of conservation laws}, Communications on
	Pure and Applied Mathematics \textbf{13} (1960), 217--237.
	
	\bibitem{Book_Leveque}
	R.~J. LeVeque, \emph{Numerical methods for conservation laws}, Birkh\"auser
	Basel, 1992.
	
	\bibitem{JCP_2004_Levy_Shu_Yan}
	D.~Levy, C.-W. Shu, and J.~Yan, \emph{Local discontinuous Galerkin methods for
		nonlinear dispersive equations}, Journal of Computational Physics
	\textbf{196} (2004), no.~2, 751--772.
	
	\bibitem{PRL_1993_Rosenau_Hyman}
	P.~Rosenau and J.~M. Hyman, \emph{Compactons: solitons with finite wavelength},
	Physical Review Letters \textbf{70} (1993), no.~5, 564--567.
	
	\bibitem{PRL_2007_Rosenau_Hyman_Staley}
	P.~Rosenau, J.~M. Hyman, and M.~Staley, \emph{Multidimensional compactons},
	Physical Review Letters \textbf{98} (2007), no.~2, 024101.
	
	\bibitem{MCS_2007_Rus_Villatoro}
	F.~Rus and F.~R. Villatoro, \emph{Pad{\'e} numerical method for the
		Rosenau-Hyman compacton equation}, Mathematics and Computers in Simulation
	\textbf{76} (2007), no.~1-3, 188--192.

	
	\bibitem{C&F_2018_shi_shu}
	C.~Shi and C.-W. Shu, \emph{On local conservation of numerical methods for
		conservation laws}, Computers \& Fluids \textbf{169} (2018), 3--9.
	
	\bibitem{SIREV_2009_Shu}
	{C.-W. Shu}, \emph{High order weighted essentially nonoscillatory schemes for
		convection dominated problems}, SIAM Review \textbf{51} (2009), no.~1,
	82--126.
	
	\bibitem{JCP_1989_Shu_Osher}
	{C.-W. Shu} and {S. Osher}, \emph{Efficient implementation of essentially
		non-oscillatory shock-capturing schemes, II}, Journal of Computational
	Physics \textbf{83} (1989), no.~1, 32--78.
	
	\bibitem{CMAME_2007_Xu_Shu}
	Y.~Xu and C.-W. Shu, \emph{Error estimates of the semi-discrete local
		discontinuous Galerkin method for nonlinear convection--diffusion and KdV
		equations}, Computer Methods in Applied Mechanics and Engineering
	\textbf{196} (2007), no.~37-40, 3805--3822.
	
	\bibitem{MCOM_2026_Yu_Cheng_Liu_Shu}
	K.~Yu, J.~Cheng, Y.~Liu, and C.-W. Shu, \emph{A Lax-Wendroff type theorem of
		doubly conservative schemes for degenerate convection-diffusion equations},
	Mathematics of Computation \textbf{https://doi.org/10.1090/mcom/4167} (2026).
	
	\bibitem{arxiv_2026_Yu_Cheng_Shu}
	K.~Yu, J.~Cheng, and C.-W. Shu, \emph{A high-order doubly conservative weighted
		essentially non-oscillatory finite volume scheme for degenerate
		convection-diffusion equations on unstructured meshes}, submitted.
	
	\bibitem{ZETF_1974_ZK-equation}
	V.~Zakharov and E.~Kuznetsov, \emph{On three-dimensional solitons}, Zhurnal
	Eksperimental'noi i Teoreticheskoi Fiziki \textbf{66} (1974), 594--597.
	
	\bibitem{SISC_2012_ZLS_fvm}
	{X. Zhang}, {Y. Liu}, and {C.-W. Shu}, \emph{Maximum-principle-satisfying high
		order finite volume weighted essentially nonoscillatory schemes for
		convection-diffusion equations}, SIAM Journal on Scientific Computing
	\textbf{34} (2012), no.~2, A627--A658.
	
	\bibitem{CAMC_2026_Zheng_Xu}
	W.~Zheng and Y.~Xu, \emph{Invariants preserving time-implicit local
		discontinuous Galerkin schemes for high-order nonlinear wave equations},
	Communications on Applied Mathematics and Computation \textbf{8} (2026),
	no.~1, 149--176.
	
	\bibitem{JCP_2018_Zhu_Shu_MR-WENO}
	{J. Zhu} and {C.-W. Shu}, \emph{A new type of multi-resolution WENO schemes
		with increasingly higher order of accuracy}, Journal of Computational Physics
	\textbf{375} (2018), 659--683.
	
	\bibitem{Nonlinearity_2018_Zilburg_Rosenau}
	A.~Zilburg and P.~Rosenau, \emph{Loss of regularity in the $K(m,n)$ equations},
	Nonlinearity \textbf{31} (2018), no.~6, 2651--2665.
	
\end{thebibliography}
%
%
\providecommand{\bysame}{\leavevmode\hbox to3em{\hrulefill}\thinspace}
\providecommand{\MR}{\relax\ifhmode\unskip\space\fi MR }
\providecommand{\MRhref}[2]{%
	\href{http://www.ams.org/mathscinet-getitem?mr=#1}{#2}
}
\providecommand{\href}[2]{#2}

\end{document}